\documentclass[num-refs]{wiley-article}

\usepackage{graphicx}
\usepackage{amsmath}
\usepackage{epstopdf} 

\usepackage[english]{babel}
\usepackage{amsmath,amssymb,amsxtra,latexsym}
\usepackage{graphicx}
\usepackage{color}
\usepackage{epsfig}
\usepackage{amsmath}%
\usepackage{amsfonts}%
\usepackage{amssymb}
\usepackage{float}
\usepackage{placeins}
\usepackage{algorithm}
\usepackage[noend]{algpseudocode}
\usepackage{multicol}
\usepackage{afterpage}

\newcommand\norm[1]{\left\lVert#1\right\rVert}
\newcommand{\sgn}{\text{sgn}}
\newcommand{\MM}{\text{MM}}
\newcommand{\am}[1]{#1}

\numberwithin{equation}{section}
\papertype{Original Manuscript}

\title{Convergence of a La\-gran\-gian--Eu\-le\-rian scheme via weak asymptotic analysis for one-dimensional hyperbolic problems}

\author[1\authfn{1}]{Eduardo Abreu}
\author[2]{Arthur Espírito Santo}
\author[3\authfn{2}]{Wanderson Lambert}
\author[4]{John Pérez}

\affil[1]{IMECC -	Universidade Estadual de Campinas  - Campinas, SP - Brazil}
\affil[2]{IME - Universidade Federal do Rio Grande do Sul - Porto Alegre - RS - Brazil.}
\affil[3]{Universidade Federal de Alfenas - Po\c{c}os de Caldas, MG - Brazil.}
\affil[4]{Instituto Tecnológico Metropolitano, Medell\'{\i}n, Colombia.}

\corraddress{Arthur Espírito Santo, IME - Universidade Federal do Rio Grande do Sul - Porto Alegre - RS - Brazil.}
\corremail{arthurmes@ufrgs.br}

\fundinginfo{FAPESP (2019/20991-8)\authfn{1}\authfn{2}; CNPq 306385/2019-8)\authfn{1}; PETROBRAS (2015/00398-0)\authfn{1}.}

\runningauthor{E. Abreu, A. Espírito Santo, W. Lambert, J. Pérez}

\begin{document}

\maketitle

\begin{abstract}
In this paper, we study both convergence and bounded variation properties of a new fully discrete conservative La\-gran\-gian--Eu\-le\-rian scheme to the entropy solution in the sense of Kruzhkov (scalar case) by using a weak asymptotic analysis. We discuss theoretical deve\-lopments on the conception of {\it no-flow curves} for hyperbolic problems within scientific computing. The resulting algorithms have been proven to be effective to study nonlinear wave formations and rarefaction interactions. We present experiments to a study based on the use of the Wasserstein distance to show the effectiveness of the {\it no-flow curves} approach in the cases of shock interaction with an entropy wave related to the inviscid Burgers' model problem and to a $2 \times 2$ nonlocal traffic flow symmetric system of type Keyfitz--Kranzer.

\end{abstract}
\newpage
\section{Introduction}

In this work, we present a fully discrete Lagrangian--Eulerian 
scheme, and its numerical analysis by a weak asymptotic method, 
for the computational treatment of first-order hyperbolic 
partial differential equations given by
\begin{equation} \label{noLi} \frac{\partial u}{\partial t}
+ \frac{\partial H(u)}{\partial x} = 0, 
\quad x \in \mathbb{R}, \quad t \in\mathbb{R}^+, \quad \quad \quad u(x,0) = u_0(x),
\end{equation} 
where $u=u(x,t):\mathbb{R}\times \mathbb{R}^+\to\Omega\subset \mathbb{R}$
and $H:\Omega \to \mathbb{R}$. 

To give a brief summary about the previous works that lead us 
to obtain the new method developed in this paper, we  cite first  
\cite{douglas2000locally}. In that article, an innovative locally 
conservative scheme was developed for the scalar parabolic 
convection-dominated convection--diffusion model problem of 
two-phase, immiscible, incompressible flow in porous media, the 
\textit{Locally Conservative Eulerian--Lagrangian Method} (LCELM). 
This scheme was related to the divergence form of the parabolic flow 
equation, where the use of a space-time divergence form allowed 
the localization of the transport and the desired local conservation 
property could also be localized. Such a scheme was very competitive 
from a computational point of view for scalar nonlinear transport 
problems in heterogeneous porous media. To the best of our knowledge, 
the LCELM procedure was the first work in the literature to introduce 
a {\it space-time local} conservation (by physical and geometric 
arguments), a region where the mass conservation takes place 
(locally), which was coined as {\it integral tube} along with the 
so-called integral curves for scalar parabolic convection-diffusion 
models in porous media transport problems (see Eqs. (5.4a)-(5.4b) 
in \cite{douglas2000locally}). Moreover, so-called integral 
curves were associated with points on the boundary (usually vertices) 
of the finite elements in that framework.

In this work, the fully discrete La\-gran\-gian--Eu\-le\-rian formulation 
is based on the new and substantial improvement interpretation of 
the {\it integral tube}, which is now subject to condition 
$O\left(\frac{H(u)}{u}\right) \propto \left[\frac{\Delta x}{\Delta t}\right] \to 0$ 
(called \textit{no-flow curves} \cite{EAJP19}), where the quantities 
$u$ and $H(u)$ come from the scalar initial value problem (\ref{noLi}) 
and under a suitable stability estimate that is very effective in 
computing practice; as a result, we also obtain a {\it weak CFL 
condition} that does not depend on the derivative of the flux 
function $H(u)$, but only on the introduced no-flow curves. This simple and interesting technique is the key 
ingredient of the La\-gran\-gian--Eu\-le\-rian framework that provides 
information about the local wave propagation speed. This issue 
is not discussed in \cite{douglas2000locally}. In addition, the 
{\it no-flow curves} reveal to be also a desingularization 
analysis tool for the construction of computationally stable 
schemes \cite{AMPB21,ALPS1719,EAJP19}. 
In \cite{EAJP19}, the authors also discussed a new reinterpretation 
of the La\-gran\-gian--Eu\-le\-rian no-flow curves as an anti-diffusion term 
into the viscosity coefficient defined by the quantity 
$O\left(\frac{H(u)}{u}\right)$, with a distinct identification to 
the Finite Volume framework. Here we are interested in the design 
of novel fully-discrete schemes based on the new concept no-flow 
curves subject to $O\left(\frac{H(u)}{u}\right) \propto 
\left[\frac{\Delta x}{\Delta t}\right] \to 0$ which is substantially 
different in theory foundations from the previous and relevant 
LCELM method.

\begin{figure}[t]
	    \centering
	    \includegraphics[width=.8\textwidth]{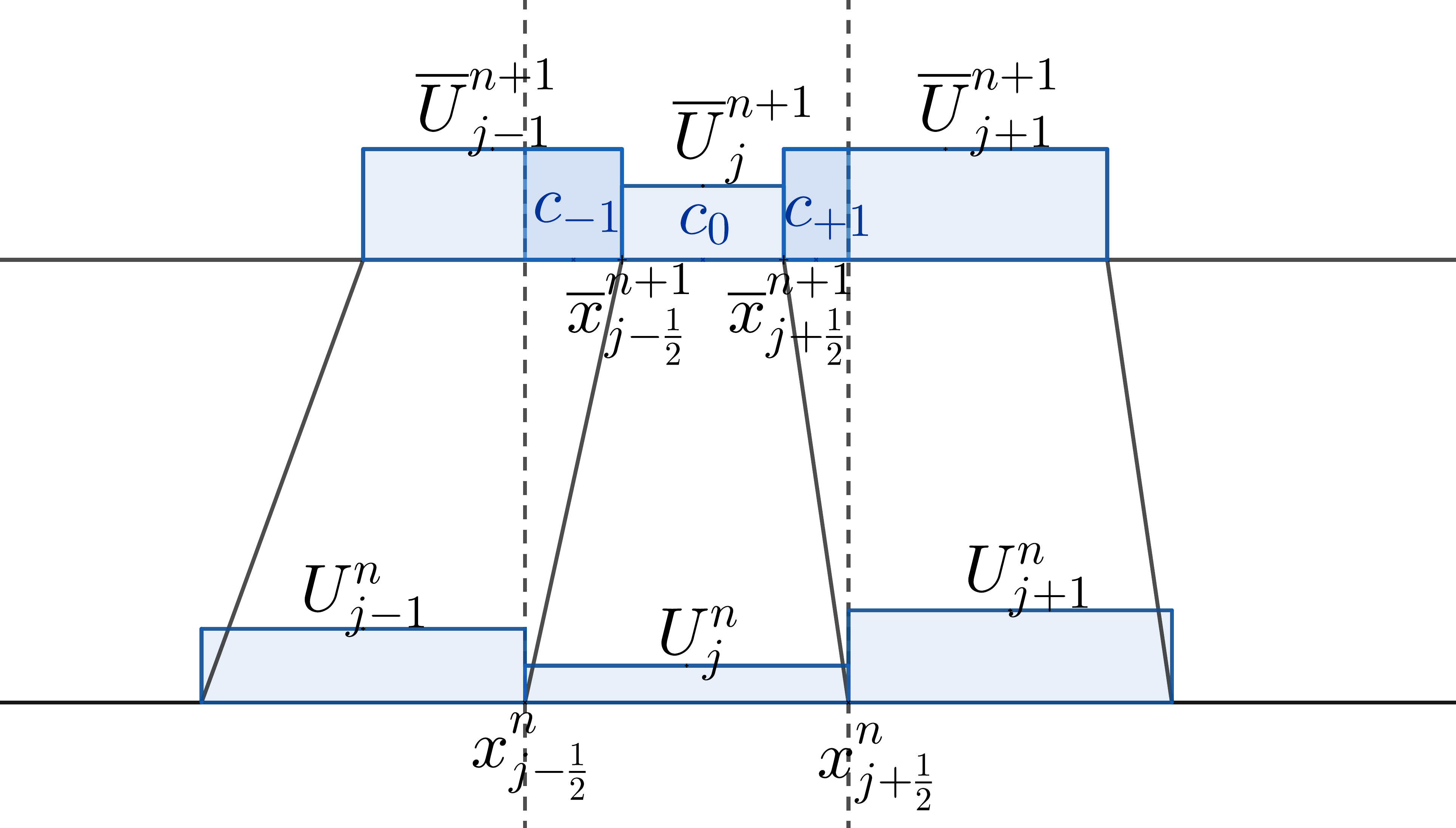}
	    \caption{Projection to original mesh.}
	    \label{figura02}
	\end{figure}

The idea of how the fully discrete explicit locally conservative 
La\-gran\-gian--Eu\-le\-rian scheme works is outlined in the following 
three basic steps, for each time step, from $t^n$ to $t^{n+1}$. 
First, by using a space-time divergence form of model (\ref{noLi}) 
we obtain an exact construction of the set of ODE modeling the 
so-called no-flow curves (see Figure \ref{figura01}). Second, 
by integrating the underlying conservation law (in divergence 
form) over the special region in the space-time domain (see 
region $D_j^{n,n+1}(x,t)$ in Figure \ref{figura01}), where the 
conservation of the mass flux takes place, we get the 
La\-gran\-gian--Eu\-le\-rian conservation relations. Third, by combining 
steps one and two, we perform an evolution (Lagrangian) step 
in time, from a primal grid to a nonstaggered grid and turn back 
by means of the final (Eulerian) projection step with a piecewise 
linear reconstruction (see Figure \ref{figura02}). A key hallmark 
of such method is the dynamic tracking forward of no-flow curves 
subject to $O\left(\frac{H(u)}{u}\right) \propto \left[\frac{\Delta x}{\Delta t}\right] \to 0$, per time step. The scheme is free of local Riemann problem solutions 
and does not use adaptive space-time discretization. This is a 
considerable improvement compared to the classical backward 
tracking over time of the characteristic curves over each time 
step interval, which is based on the strong form of the problem 
and that are not reversible for systems in general.

To analyze the properties of the {explicit} La\-gran\-gian--Eu\-le\-rian 
scheme, we use recent improvements on the weak asymptotic 
analysis (see \cite{ACP16,ACP17,MC12}), which was first
defined in a distinct framework by \cite{DOS03}. The weak 
asymptotic solutions are consistent with the traditional 
solutions in one-dimensional and multi-D (see also 
\cite{ACP16,ACP17,ALPS1719}).  The weak asymptotic analysis 
has been used to study the existence of solutions for scalar 
equations and systems of hyperbolic equations, giving solutions 
a new meaning (see 
\cite{ACP16,ACP17,MC13,VDVS05a,VDDM,DOS03,VDVS05b,BNORVS11,GAO10,EPVS06,VMS06} 
and the references therein). An interesting aspect of this 
theory is that it makes it possible to prove the existence 
(and, for the scalar case, the uniqueness) of a solution by 
means of numerical methods. 

In a previous work \cite{ALPS1719}, we defined a weak 
asymptotic solution for a scalar equation (\ref{noLi}) and 
outlined the proof of stability of the numerical method used. 
Weak asymptotic methods aim to investigate the nonlinear 
phenomena that appear in evolutionary equations (see 
\cite{ACP16,ACP17,MC12}). In addition, they have been used 
in explicit calculations by several authors to study wave 
interactions inside the solutions to Riemann or Cauchy 
problems when these solutions involve nonclassical products 
of Heaviside functions, $\delta$-Dirac distributions, and 
even their derivatives. In \cite{ACP16}, the authors show 
how one can construct families of continuous functions 
which asymptotically satisfy scalar equations with discontinuous 
nonlinearities and systems having irregular solutions. Through 
a weak asymptotic analysis, it has been proven, for scalar 
equations, that the initial value problem is well posed in 
the $L^1$ sense for the approximate solutions constructed. 
It was possible to demonstrate that the family of solutions 
obtained from the method satisfies Kruzhkov entropy, which 
provided a better understanding of the mathematical 
computations for the construction of accurate numerical 
schemes satisfying classical and entropy (Kruzhkov) solutions. 
It turns out that, beside the applications in nontrivial 
models of hyperbolic conservation laws, the convergent 
La\-gran\-gian--Eu\-le\-rian scheme via the weak asymptotic method 
treat such models with computational efficiency. From this 
theory, we obtain explicit estimates for the stability 
condition, and these estimates  are numerically implemented 
and give very good numerical solutions as we see in the 
Section \ref{NumExp}.

The proposed weak asymptotic analysis in this work, together 
with the La\-gran\-gian--Eu\-le\-rian approach 
\cite{APS17a,ALPS1719,AMPB21,EAJP19}, fits in properly with 
the classical theory while improving the mathematical 
computations for the construction of new accurate numerical 
schemes, given by Proposition \ref{prop1}, where convergence, 
existence, and stability are established. 
In the  classical treatment for proving convergence of a 
numerical scheme, first it is proved that the approximate 
solutions, generated by the numerical scheme, are of finite 
total variation (or bounded variation), actually, for scalar 
equations are proved that the total variation are not increasing. 
This condition gives the compact embedding of $BV$ functions  
in $L^1$; one can use the Helly's selection theorem to 
show pointwise convergence of the generated sequence of 
solutions. After, it is proved that the generated sequence 
is a weak solution for the scalar equations and, finally, 
it is proved that the solution satisfies an entropy 
criterion, the most commom used is the Kruzhkov entropy 
condition, see \cite{GDL}. The treatment used in the weak 
asymptotic analysis is very similar to the previous steps. 
To prove the main properties of the numerical scheme it is 
proved that the numerical scheme generates a solution with 
bounded total variation, also it is used a similar argument 
of Helly's selection theorem, see Appendix \ref{precompact}. 
Using the asymptotic theory, the solution to the proposed 
scheme is obtained as a family of functions, 
$\{u(\cdot,t,\epsilon)\}_\epsilon$, bounded in 
$\mathbb{L}^1(\mathbb{S}^1)$ uniformly in parameter 
$\epsilon$ of the model (\ref{noLi})). The main difference 
is the kind of limit that is taken in solution in Equation 
$(\ref{weaksol})$. For instance, our approach allows us 
to deal with the reconstruction (accomplished by means of 
robust choices of slope limiters) of variable $u$ of 
Eq. $(\ref{noLi})$ in the (numerical) flux terms. We also 
give sufficient conditions for a total variation 
nonincreasing (Section \ref{tvdsection}) through a suitable 
semidiscrete formulation of Eq. $(\ref{noLi})$. In addition, 
we obtain the maximum principle and the entropy (Kruzhkov) 
solution to model $(\ref{noLi})$ thanks to a suitable 
interpretation of the approximate solutions provided by 
the analysis presented in Section \ref{maximumsec}. 
The weak asymptotic theory handles reconstructions easily 
and opens new possibilities for the design of new methods 
for hyperbolic problems under a weak CFL-type stability 
depending on the no-flow curves associated with the novel 
Lagrangian--Eulerian approach instead of using estimates 
on the eigenvalues, the weak asymptotic theory fits very well with the Lagragian-Eulerian scheme and it is a very promissor technique to be used to prove convergence for improved numerical methods of Lagrangian-Eulerian class.

The rest of this paper is structured as follows. Section 
\ref{sec:lefvm} introduces the La\-gran\-gian--Eu\-le\-rian scheme. 
Section \ref{convergencee} presents the convergence of the 
numerical method via the weak asymptotic analysis. We propose 
a scheme to find a solution to $(\ref{noLi})$ and prove that 
the resulting solution converges to the weak solution of 
$(\ref{noLi})$. We demonstrate that the numerical scheme 
obtained in Section \ref{sec:lefvm} is compatible with that 
used in the weak asymptotic analysis. Section \ref{tvdsection} 
proves that our scheme has a total variation nonincreasing 
property for the solution. To complete our analysis, Section 
$\ref{maximumsec}$ demonstrates that the obtained solution 
satisfies the maximum principle and Kruzhkov entropy 
solution. Appendix \ref{lip} shows evidence that the 
reconstructions used in this study are Lipschitz continuous. 
Section \ref{NumExp} shows and discusses a set of 
representative computational results, including numerical 
studies with the W1 distance. Section \ref{ConRem} 
summarizes our concluding remarks.

\section{La\-gran\-gian--Eu\-le\-rian scheme}\label{sec:lefvm}

To construct the fully discrete La\-gran\-gian--Eu\-le\-rian 
scheme, we first consider the scalar one-dimensional (1D) 
Cauchy problem (\ref{noLi}) in its divergent form 
\begin{equation}
\nabla\cdot \left[\begin{array}{c}
H(u) \\ u
\end{array}\right] = 0,
\quad t > 0, \quad  x\in \mathbb{R}, 
\qquad u(x,0)=u_{0}(x), \quad   x\in \mathbb{R}, \quad \text{ where } \nabla=\left(\frac{\partial}{\partial x},\frac{\partial}{\partial t}\right). 
\label{div}
\end{equation}
For the sake of clarity and completeness, we sketch in what follows the main three steps of the La\-gran\-gian--Eu\-le\-rian approach, namely, the construction of the set of ODE modeling the no-flow curves (Section \ref{subnflow}),
the La\-gran\-gian--Eu\-le\-rian conservation relations (Section \ref{secConRel}) and the final projection step by piecewise linear reconstruction (Section \ref{projStep}), along with an evolution algorithm of the fully discrete nonstaggered Lagrangian-Eulerian method.
As in the La\-gran\-gian--Eu\-le\-rian schemes present in \cite{APS17a,ALPS1719,AMPB21,EAJP19}, local conservation is obtained by integrating the conservation law over the region in the space-time domain where the conservation of mass flux takes place.

	\subsection{\textbf{La\-gran\-gian--Eu\-le\-rian no-flow curves}}\label{subnflow}
	Consider the La\-gran\-gian--Eu\-le\-rian control volumes
	\begin{equation} 
	{ D_j^{n,n+1} }= \{(x,t) \,\, / \,\, \sigma^n_{j-\frac{1}{2}}(t) 
\leq x \leq \sigma^n_{j+\frac{1}{2}}(t), \,\, t^n \leq t \leq t^{n+1} 
	\}, 
	\label{ler}
	\end{equation}
	where $\sigma_{j\pm\frac{1}{2}}^n(t)$
	are the {\it La\-gran\-gian--Eu\-le\-rian no-flow curves}
such that $\sigma_{j\pm\frac{1}{2}}^n(t^n) = x_{j\pm\frac{1}{2}}^n$
and subject to $O\left(\frac{H(u)}{u}\right) \propto \left[\frac{\Delta x}{\Delta t}\right] \to 0$ \cite{ALPS1719,AMPB21,EAJP19} along with $u$ and $H(u)$ given from (\ref{div}).
	These curves correspond to the lateral boundaries
	of domain ${ D_j^{n,n+1} }$ in (\ref{ler}), and we define 
	$\bar{x}_{j\pm\frac{1}{2}}^{n} := 
	\sigma_{j\pm\frac{1}{2}}^n(t^{n+1}) $ as their 
	endpoints in time $t^{n+1}$. For each control volume, $h_j^n$ is defined as $h_j^n = x_{j+\frac{1}{2}}^n-x_{j-\frac{1}{2}}^n$.

\begin{figure}[ht]
		\centering 
		\includegraphics[scale=0.2]{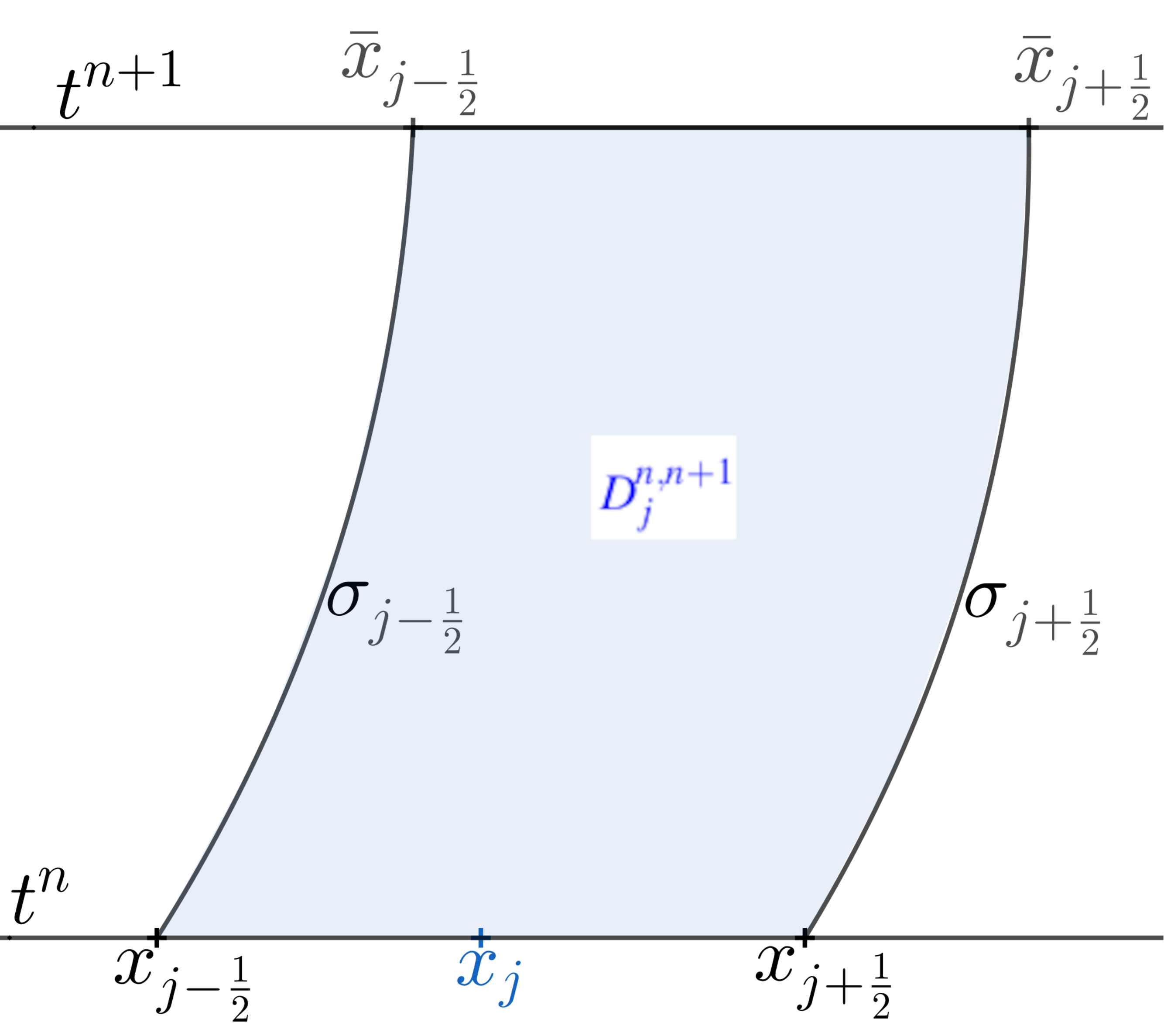}\qquad
		\includegraphics[scale=0.2]{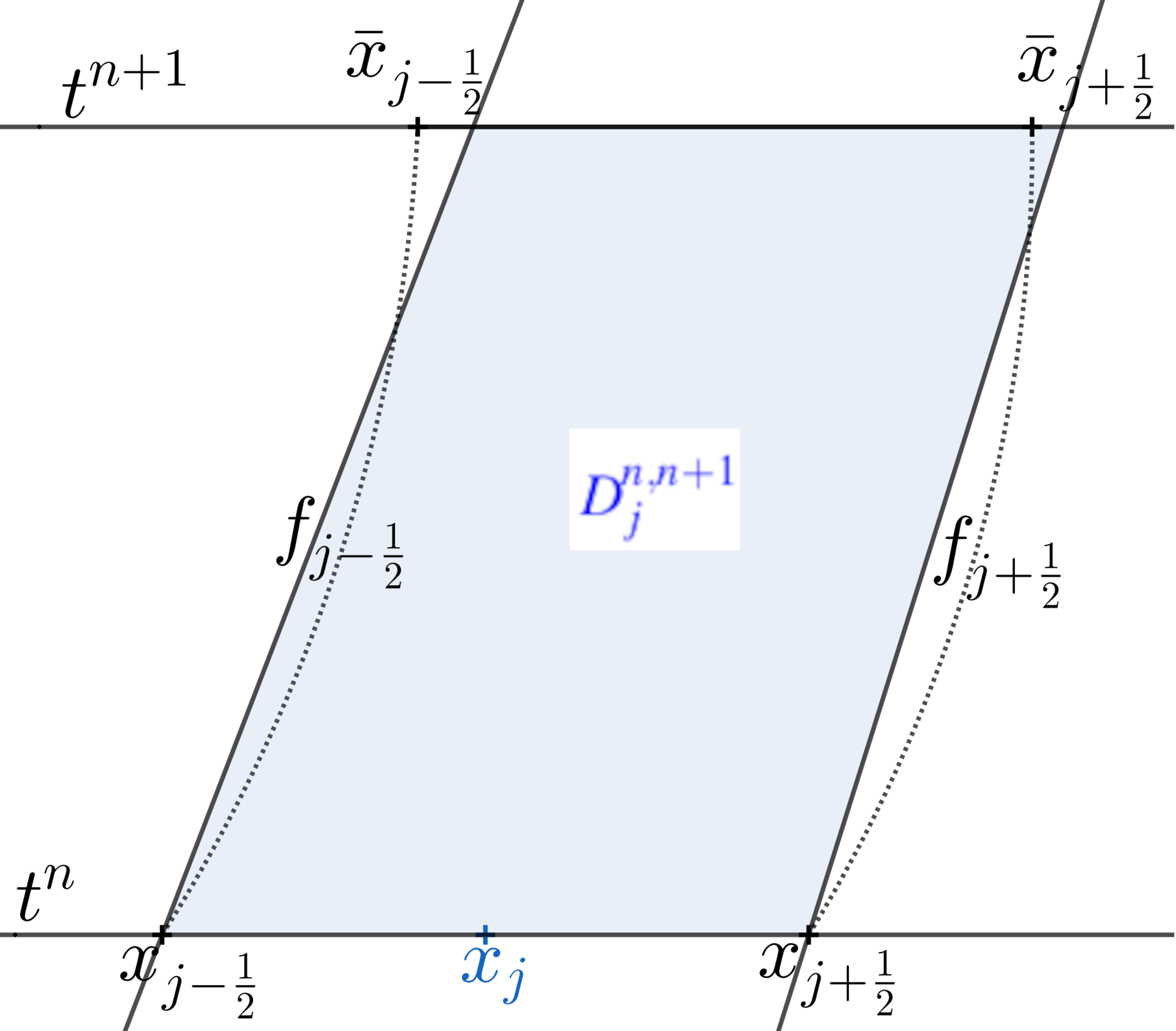}
		\caption{Left: No flow region ${ D_j^{n,n+1} }$. Right: Linear approximation for ${ D_j^{n,n+1} }$.}
		\label{figura01}
	\end{figure}

	The numerical scheme is expected to satisfy certain 
	type of mass conservation (due to the inherent 
	nature of the conservation law)
	from time $t^n$ in the space domain 
	$[x_{j-\frac{1}{2}}^n,x_{j+\frac{1}{2}}^n]$ to 
	time $t^{n+1}$ in the space domain  
	$[\bar{x}_{j-\frac{1}{2}}^{n+1},\bar{x}_{j+\frac{1}{2}}^{n+1}]$;
        see also \cite{douglas2000locally} for original motivation foundations on models for the flow of fluids and the transport of {\it mass conservation} in porous geologic media.
	Based on this, the flux through curves 
	$\sigma_{j\pm\frac{1}{2}}^n(t)$ must be zero. With 
	the integration of (\ref{noLi}) and the Divergence
	Theorem, and because the line integrals 
	over curves $\sigma_{j\pm\frac{1}{2}}^n(t)$ vanish, 
	\begin{equation}
	\int_{\bar{x}_{j-\frac{1}{2}}^{n+1}}^{\bar{x}_{j+\frac{1}{2}}^{n+1}}
	u(x,t^{n+1})dx = \int_{x_{j-\frac{1}{2}}^{n}}^{x_{j+\frac{1}{2}}^{n}}u(x,t^{n})dx.
	\label{conserv}
	\end{equation}
        Here 
	curves $\sigma_{j\pm\frac{1}{2}}^n(t)$ are not straight 
	lines in general but rather solutions to the
	set of local nonlinear differential equations 
	(as in \cite{EAJP19,AMPB21} for systems):
	$
	\frac{d }{dt}[\sigma_{j\pm\frac{1}{2}}^n(t)]= \frac{H(u)}{u}$, 
	for $t^n < t \leq t^{n+1},
	$
	with the initial condition $\sigma_{j\pm\frac{1}{2}}^n(t^n) = x_{j\pm\frac{1}{2}}^n$, 
	assuming $u\neq0$ (for the sake of presentation).  
	This construction follows naturally from the finite 
	volume formulation of the linear La\-gran\-gian--Eu\-le\-rian 
	scheme (see Figure \ref{figura01}) as a building block 
        to construct {\it local} approximations for $\frac{H(u)}{u}$, such as 
	$f_{j\pm\frac{1}{2}}^n = \frac{H(U_{j\pm\frac{1}{2}}^n)}{U_{j\pm\frac{1}{2}}^n}$, with the initial condition 
	$\sigma_{j\pm\frac{1}{2}}^n(t^n) = x_{j\pm\frac{1}{2}}^n$. \begin{remark}
	In fact, distinct and high-order approximations are 
	also acceptable for $\frac{d }{dt}[\sigma_{j\pm\frac{1}{2}}^n(t)]$ 
	and can be regarded as ingredients to improve the accuracy 
	of the new family of La\-gran\-gian--Eu\-le\-rian methods.\end{remark}
	Eq. (\ref{conserv}) defines the mass conservation 
	but in a different mesh cell-centered in points 
	$\bar{x}_{j+\frac{1}{2}}^n$ of width $h_j^{n+1}$ defined by
    $h^{n+1}_j=h^n_j + (f^n_{j+\frac{1}{2}}
	-f^n_{j-\frac{1}{2}})\Delta t$, which gives us $\frac{h^n_j}{h^{n+1}_{j}} = 1-\frac{{f^n_{j+\frac{1}{2}}-f^n_{j-\frac{1}{2}}}}{{h^{n+1}_{j}}}\Delta t.\label{hnj}
	$
	Along the linear approximations of $f_{j\pm\frac{1}{2}}^n$, 
	we have that $\bar{x}_{j\pm\frac{1}{2}} = 
	x_{j\pm\frac{1}{2}}+f^n_{j\pm\frac{1}{2}} \Delta t$. 
	\subsection{\textbf{Conservation relations}}\label{secConRel}
	Equation (\ref{conserv}) defines a local mass balance between space 
	intervals at times $t^n$ and $t^{n+1}$. 
	By defining
		\begin{equation*}
			\overline{U}_j^{n+1} := \frac{1}{h_j^{n+1}}
	\int_{\bar{x}_{j-\frac{1}{2}}^{n+1}}^{\bar{x}_{j+\frac{1}{2}}^{n+1}}
	u(x,t^{n+1})dx, \quad \text{and} \quad U_j^n := 
	\frac{1}{h^n_j}\int_{x_{j - 
			\frac{1}{2}}^{n}}^{x_{j + \frac{1}{2}}^{n}}u(x,t^{n})dx,
	\end{equation*}
	the discrete version of Eq. (\ref{conserv}) is given by
	\begin{equation}
	\overline{U}_j^{n+1}  =  \frac{1}{h_j^{n+1}}
	\int_{\bar{x}_{j-\frac{1}{2}}^{n+1}}^{\bar{x}_{j+\frac{1}{2}}^{n+1}}
	u(x,t^{n+1})dx  = \frac{1}{h_j^{n+1}}
	\int_{x_{j-\frac{1}{2}}^{n}}^{x_{j+\frac{1}{2}}^{n}}u(x,t^{n})dx
	= \frac{h^n_j}{h_j^{n+1}} U_j^n.
	\label{eq32}
	\end{equation}
	
	Solutions $\sigma_{j\pm\frac{1}{2}}^n(t)$ to the differential 
	system are also obtained using linear approximations 
	$L(x,t)$. A piecewise constant 
	numerical data is reconstructed into a piecewise linear 
	approximation (although high-order reconstructions are 
	acceptable) by means of MUSCL-type interpolants 
	$L_j(x,t) = u_j(t)+(x-x_j)\frac{1}{h^n_j} u'_j$. 
	
	\begin{remark}For the numerical derivative $\frac{1}{h^n_j}u'_j$, 
	there are many choices of slope limiters (see, e.g., 
	\cite{AMPB21,APS17a}). Even though selecting such slope limiters 
	a priori is quite hard, this choice is based on the underlying 
	model problem under study.\end{remark} 
	
	The approximation of $U^n_{j-\frac{1}{2}}$ is
	\begin{equation}
	\begin{array}{ll}
	U^n_{j-\frac{1}{2}} & = \frac{1}{h^n_j}\int_{x_{j-1}^n}^{x_{j}^n} L(x,t) dx 
	=\frac{1}{h^n_j}\left( \int_{x_{j-1}^n}^{x_{j-\frac{1}{2}}^n} L_{j-1}(x,t) dx
	+ \int_{x_{j-\frac{1}{2}}^n}^{x_{j}^n} L_j(x,t) dx \right) \\ \\
	& = \frac{1}{2} (U^n_{j-1} + U^n_j) + \frac{1}{8} ({U_j^n}' - {U_{j-1}^n}').
	\end{array}   
	\label{umeio}
	\end{equation}
	\subsection{\textbf{Projection step}}\label{projStep}
	Next, for a partition with constant spacing  
        ($h^n_j = h, \forall j$), we obtain the resulting projection formula as follows:
	\begin{equation}
	\label{projection}
	U_{j}^{n+1} = \frac{1}{h}
	\left( c_{-1,j}\overline{U}_{j-1}^{n+1}+ c_{0,j}
	\overline{U}_{j}^{n+1}+ c_{+1,j}\overline{U}_{j+1}^{n+1} 
	\right),
	\end{equation}
 \noindent where the projection coefficients are 
    	\begin{equation}
    	c_{-1,j} = f^+(U^n_{j-\frac{1}{2}}) \Delta t,\quad 
    	c_{+1,j} = f^-(U^n_{j+\frac{1}{2}}) \Delta t,\quad
    	c_{0,j} =h-c_{-1,j}-c_{+1,j},\label{c2}
    	\end{equation}    
\noindent and $f^+$, $f^-$ defined as
	\begin{equation}
	f^+=\max(f,0) \quad\text{ and }\quad 
	f^-=\max(-f,0).\label{fmm}
	\end{equation} 	
	Note that, for function $f$, we have 
$	f=f^+-f^-$ and $
	|f|=f^++f^-.$	
Here $\Delta t$ is obtained under a {\it weak CFL condition}, that does not depend on the derivative of the flux function $H(u)$, but only based in the so-called no-flow curves $f^n_{j-\frac{1}{2}}$, 
	\begin{equation}
	{ \max_j \left\{ |f^n_{j-\frac{1}{2}}  \, 
	\Delta t| \right\} \leq \frac{h}{2},\label{CFL}
	}\end{equation}
	which is taken by construction of the method. We note that, 
	in the linear case, when $a(x,t) = a > 0$ (or $a<0$), numerical
	scheme (\ref{eq32})--(\ref{projection}) is a generalization
	of the upwind scheme. However, our scheme can provide an approximate solution
	in both scenarios: $a>0$ and $a<0$. In this case, the CFL condition
	is $|a \, \Delta t | \leq h$, as in the upwind scheme.
For the sake of clarity and completeness, we include the numerical experiment with the sonic rarefaction for the invicid Burgers' problem. Here, as the rarefaction wave is crossed, there is a sign change in the characteristic speed $u$, in the range $-1 < u < 1$ and thus there is one point at which $u = 0$, the sonic point. Our approach does not use Riemann solvers neither Godunov type implementation. {We can summarize the ideas presented so far in the following algorithm.}

    \begin{algorithm}
        \renewcommand{\thealgorithm}{}
        \caption{Nonstaggered Lagrangian-Eulerian method}\label{alg:nsle}
        \begin{algorithmic}
        	{
            \Procedure{NSLEstep}{$h^n$, $U^n$, $t^{n}$, $t^{n+1}$}
            \State \textbf{1.} Approximate the no-flow boundaries slope with $f_{j+\frac{1}{2}}^n \gets \frac{H(U^n_{j+\frac12})}{U^n_{j+\frac12}}$.
            \State \textbf{2.} Find the new mesh $h_j^{n+1} \gets h_j^n +(f_{j+\frac12}-f_{j-\frac12})\Delta t$.
            \State \textbf{3.} Compute the conserved quantity $\overline{U}_j^{n+1}$ with (\ref{eq32});  
            \State \textbf{4.} Compute $U^n_{j-\frac12}$ using an appropriate slope limiter with (\ref{umeio}).
            \State \textbf{5.} Compute the projection coefficients $c_{-1,j}, c_{0,j}$ and $c_{+1,j}$ with (\ref{c2}).
            \State \textbf{return} the values $U^{n+1}_j$ projected from $\overline{U}_j^{n+1}$ into the original mesh with (\ref{projection}).
            \EndProcedure
        }
        \end{algorithmic}
    \end{algorithm}

	We now examine the theoretical properties of the 
	La\-gran\-gian--Eu\-le\-rian scheme via weak asymptotic solutions. 

	\subsection{\textbf{Slope limiters}}
	\label{slope}
	States $U^n_{j-\frac{1}{2}}$ are obtained from states 
	$U^n_{j-1}$ and $U^n_j$ and functions $U^\prime_{j-1}$ and 
	$U^\prime_{j}$ in timestep $n$. These functions (obtained using
	\textit{slope limiters}) compute variations 
	of function $U$ at the neighborhood of $U_{j-\frac{1}{2}}$. 
	We would like to stress that said functions ($U^\prime_{j-1}$ and 
	$U^\prime_{j}$) are not numerical derivatives because 
	they can compute variations even for noncontinuous 
	$U$. In Appendix $\ref{lip}$, we demonstrate that they 
	are Lipschitz continuous.
	The first option of slope limiter is
	\begin{equation}
	U'_j = \MM_2\left(\Delta u_{j+\frac{1}{2}}, 
	\Delta u_{j-\frac{1}{2}}\right),
	\label{mm2}
	\end{equation}
	where $\MM_2$ corresponds to the most common limiter 
        (see, e.g., \cite{APS17a,ALPS1719}), with
	$\Delta u_{j+\frac{1}{2}}=u_{j+1}-u_{j}$, and can be written as
	\begin{equation}
	\MM_2\left(\sigma,\tau\right) = 
	\frac{1}{2}\lambda(\sigma,\tau) \min
	\left(|\sigma|,|\tau|\right),
	\label{mm}
	\end{equation}
	where $\lambda(u,v) = \sgn(u)+\sgn(v)$.
	The following option, which allows steeper slopes near 
	discontinuities and retains accuracy in smooth 
	regions (obtained between three values), can also be 
        used and is given by\footnote{The range of parameter 
		$\alpha$ is typically guided by the CFL condition.}
	\begin{equation}
	U'_j = \MM_3\left(\alpha \Delta u_{j+\frac{1}{2}}, 
	\frac{1}{2}(u_{j+1}-u_{j-1}),\alpha 
	\Delta u_{j-\frac{1}{2}}\right),
	\label{mm3}
	\end{equation}
	where $\MM_3$ can be written as
	\begin{equation}
	\MM_3\left(\sigma,\tau,\gamma\right)=
	\frac{1}{8}\lambda(\sigma,\tau)\lambda(\sigma,\gamma)\lambda(\tau,\gamma)\min\left(|\sigma|,|\gamma|,|\tau|\right),\label{mmn}
	\end{equation}
	One can prove that
$	\MM_3(x,y,z)=\MM_2\left(\MM_2(x,y),z\right).
$	
	A third option is the high-order slope limiter UNO, which is given by 
	\begin{equation}
	U'_j = \MM_2\left( \Delta u_{j+\frac{1}{2}}-\delta^2_1,
	\Delta u_{j-\frac{1}{2}}+\delta^2_2 \right),
	\label{uno}
	\end{equation}
	where $\delta^2_1$ is a function of $u_{j+2}$, $u_{j+1}$, $u_{j}$ 
and $u_{j-1}$ defined by $\delta^2_1 = \frac{1}{2}\MM\left(\Delta^2 
	u_{j+1},\Delta^2 u_{j} \right)$; and $\delta^2_2$, of $u_{j+1}$, 
$u_{j}$, $u_{j-1}$ and $u_{j-2}$ defined by $\delta^2_2 = \frac{1}{2}\MM\left(\Delta^2 
	u_{j},\Delta^2 u_{j-1} \right)$ and $\Delta^2u_j = u_{j+1}-2u_j+u_{j-1}$. 	
	In Appendix $\ref{lip}$, we demonstrate that the 
	reconstructions of $U_{j+\frac{1}{2}}$ are 
	Lipschitz continuous. 

	\section{The convergence proof of the La\-gran\-gian--Eu\-le\-rian scheme}\label{convergencee}
	
	In the  classical treatment for proving convergence of a numerical scheme, first it is proved that the approximate solutions, generated
by the numerical scheme, are of finite total variation (or bounded variation), actually, for scalar equations are proved that the total variation are not increasing. This condition gives the compact embedding of $BV$ functions  in $L^1$,  one can use the Helly's selection theorem to show compactness for the approximations and establish
pointwise convergence. After, it is proved that the solution is a weak solution for the scalar equations and, finally, it is proved that the solution satisfies an entropy criterion, the most commom used is the Kruzhkov entropy condition, see \cite{GDL}. In this Section, we show the convergence of our scheme using the weak asymptotic theory, that is very suitable to handle with the Eulerian-Lagrangian approach. From this theory, we obtain explicit estimates for the stability condition, these estimates are implemented and give very good numerical solutions as one can see  in the Section \ref{NumExp}. 

To establish the weak asymptotic solution, we consider 
	the one-dimensional scalar equation in (\ref{noLi}). 
	In order to avoid boundary conditions in the bounded 
	domain due to numerical purposes, we consider that 
	$x\in \mathbb{S}^1=\mathbb{R}/\mathbb{Z}$; variable 
	$t\in \mathbb{R}^+$, thus $u=u(x,t):\mathbb{S}^1\times 
	\mathbb{R}^+\to\Omega\subset\mathbb{R}$, and
	the flux function $H=H(u):\Omega \to \mathbb{R}$  is assumed to be a locally 
	Lipschitz function in $u$, i.e., 
	
	\textit{Assumption 1 -}  for all $c>0$, 
	$\exists L>0$, such that
\begin{equation}
	 |u_1|\leq c, |u_2|\leq c \Longrightarrow
	|H(u_1)-H(u_2)|\leq L|u_1-u_2|.     \label{h1}
	\end{equation}

	The 
	weak asymptotic solution is a sequence of solutions 
	$(u_\epsilon)_\epsilon=(u(x,t,\epsilon))_\epsilon$ 
	of class $\mathcal{C}^1$ in $t$ and class
	$L^\infty$ and is piecewise continuous in $x$ such 
	that, for all $\psi=\psi(x) \in 
	\mathcal{C}^\infty_c(\mathbb{R})$ (smooth 
	with compact support), and, 
	for all $t$,
	\begin{equation}
	\lim_{\epsilon \to 0} \int_R  
	((u_\epsilon)_t\psi-H(u_\epsilon)\psi_x)dx=0 
	\quad \text{and} \quad
	u_\epsilon(x,0)=u_0(x).
	\label{weaksol}
	\end{equation}
	In the weak asymptotic method, a PDE with a special flux (using 
	parameter $\epsilon$) is first 
	proposed. Then, for each fixed $\epsilon$, 
	we obtain an Ordinary Differential Equation (ODE) for variable $t$. 
	Based on the theory of ODEs, we prove the existence and 
	stability of the solution. Finally, we demonstrate 
	that, when taking $\epsilon \to 0$, the limit 
	satisfies ($\ref{weaksol}$). The idea is for 
	the flux to represent the numerical method, which is why 
	the existence and stability of the PDE with 
	special class may be an extension of 
	the numerical method.
	In our method, we use an auxiliary function 
	($f(u)=H(u)/u$) and assume that $u\neq 0$ 
	to avoid technical details (we would like to stress that 
	the convergence of the method can be proven
	in this case). Moreover, we assume that there 
	exists a $a>0$ such that $u>a>0$. Note that, 
	in this case, $f(u)$ is   a locally Lipschitz function in $u$ 
	because for all $c>a$, 
	$\exists \tilde{K}>0$ such that
	\begin{equation*}
	 a<u_1\leq c, \;  a<u_2\leq c \Longrightarrow
	|f(u_1)-f(u_2)|\leq \tilde{K}|u_1-u_2|.     \label{f1}
	\end{equation*} 
	Notice that we can take $\tilde{K}=L/a+\tilde{M}_c/a^2$, where $\displaystyle{\tilde{M}_c=\max_{a<u\leq c}|H(u)|}$.
	
	For our method, we propose the following {ODE}:
	\begin{equation}
	\partial_t (u_{_{\epsilon}})=\frac{1}{\epsilon}
	\left[u_{_{\epsilon-1}}f^+(\hat{u}_{_{\epsilon-1/2}})-
	u_{_{\epsilon}}f^+(\hat{u}_{_{\epsilon+1/2}})-
	u_{_{\epsilon}}f^-(\hat{u}_{_{\epsilon-1/2}})+
	u_{_{\epsilon+1}}f^-(\hat{u}_{_{\epsilon+1/2}})
	\right],\label{ODE}
	\end{equation}  
	with initial condition $u_\epsilon(x,0)=u_0(x,0)$. 
	Here $f^+$ and $f^-$ are defined in $(\ref{fmm})$  
	and $u_{_{\epsilon-i}}$, $u_{_{\epsilon+i}}$, 
	and $u_{_{\epsilon}}$ are defined as
	$\displaystyle{u_{_{\epsilon-i}}=u(x-i\epsilon,t,\epsilon)}$, $\displaystyle{ 
	u_{_{\epsilon+i}}=u(x+i\epsilon,t,\epsilon)}$  and  
	$\displaystyle{u_{_{\epsilon}}=u(x,t,\epsilon)}$.
	Functions ${f^+\left(\hat{u}_{_{\epsilon-1/2}}\right)}$ and 
	$f^-\left(\hat{u}_{_{\epsilon-1/2}}\right)$ are evaluated in 
	the middle of each cell, i.e, they are given by following identities ${f^+\left(\hat{u}_{_{\epsilon-1/2}}\right)
	=f^+\left(\hat{u}_{_{\epsilon-1/2}},x-\frac{\epsilon}{2},t\right)}$,\\ 
	${f^+\left(\hat{u}_{_{\epsilon+1/2}}\right)=f^+\left(\hat{u}_{_{\epsilon+1/2}},x+\frac{\epsilon}{2},t\right)}$,  
	${f^-\left(\hat{u}_{_{\epsilon-1/2}}\right)=f^-\left(\hat{u}_{_{\epsilon-1/2}},x-\frac{\epsilon}{2},t\right)},$ 
	and ${f^-\left(\hat{u}_{_{\epsilon+1/2}}\right)=f^-\left(\hat{u}_{_{\epsilon+1/2}},x+\frac{\epsilon}{2},t\right)}$. 

 States $\hat{u}_{_{\epsilon-1/2}}$ and 
	$\hat{u}_{_{\epsilon+1/2}}$ are obtained from a 
	combination of known states from equations such 
	as $(\ref{umeio})$ and are also determined in 
	the middle of each cell. 
	For the general case, we assume that there is a 
	function $L(x_{-p},\dots,x_p)$ so that we 
	can define, for instance, $\hat{u}_{{\epsilon+1/2}}$ as:
	\begin{equation}
	\hat{u}_{{\epsilon+1/2}}=L(u(x-p\epsilon,t,\epsilon),
	u(x-(p-1)\epsilon,t,\epsilon),\dots,u(x+(p-1)
	\epsilon,t,\epsilon),u(x+p\epsilon,t,\epsilon)).\label{ffdse}
	\end{equation}  
	Function $L$ is generated using the slope limiters 
	(defined in Section $\ref{slope}$). In addition, we 
	assume the following compatibility condition: for any $x$,
$	x=L(x,x,\dots,x).$	In Appendix $\ref{lip}$, we prove that the slope limiters, 	as well as the reconstructions of type $(\ref{umeio})$,
	are Lipschitz continuous. In this case, if 
	$u(x+i\epsilon,t,\epsilon)$ are continuous 
	functions
\begin{equation}
	\lim_{\epsilon\longrightarrow 0}|u_{\epsilon}-\hat{u}_{\epsilon-1/2}|=\lim_{\epsilon\longrightarrow 0}|u_{\epsilon}-\hat{u}_{\epsilon+1/2}|=0.\label{fds}
	\end{equation}
	Note that since $f(u)$ and $L$ are Lipschitz continuous,  then $f$ 
	applied to $\hat{u}_{_{\epsilon-1/2}}$ is also a Lipschitz function.
	
	{In the next result, we state the existence and stability result of the solution to $(\ref{ODE})$. As strategy of proof, we apply the  Taylor expansion with remainder term to substitute the derivative by a difference defined between $t$ and $t+dt$, where $dt$ represents a time step. We state our result  as follows}:
	\begin{proposition}
		\label{prop1}
		We construct, as a solution to $(\ref{ODE})$, a family 
		of functions $(x,t)\to 
		u(x,t,\epsilon):\mathbb{S}^1\times\mathbb{R}\to 
		\mathbb{R}$, for small enough $\epsilon$, that 
		are of class $\mathbb{C}^1$ { in $t$} for a fixed $\epsilon$
		and of class $\mathbb{L}^\infty$ for $x\in \mathbb{S}^1$ 
		and satisfy $(\ref{weaksol})$. If
		\begin{equation}\frac{dt}{\epsilon}
		\left(f^+\left(\hat{u}_{_{\epsilon+1/2}}\right)+
		f^-\left(\hat{u}_{_{\epsilon-1/2}}\right)\right)\leq 1,\quad \text{ for all } \quad
		\hat{u}_{_{\epsilon-1/2}} \text{ and } \hat{u}_{_{\epsilon+1/2}}\label{Cflcon}
		\end{equation}
		is satisfied, then the family 
		$\{u(\cdot,t,\epsilon)\}_\epsilon$ is bounded 
		in $\mathbb{L}^1(\mathbb{S}^1)$ uniformly in 
		$\epsilon$. In fact, we have that $||u(t,\epsilon)
		||_{\mathbb{L}^1(\mathbb{S}^1)}\leq 
		||u_0||_{\mathbb{L}^1(\mathbb{S}^1)}$ for all 
		$t$. Moreover, if initial condition $u_0(x)$ 
		and $H(u)$ are continuous, then $u(x,t,\epsilon)$ 
		is also continuous in $x$.  
	\end{proposition}
	Notice that $(\ref{Cflcon})$ is always valid if the CFL 
condition $(\ref{CFL})$ is satisfied. Indeed, 
the proof of Proposition $\ref{prop1}$ follows ideas from works 
\cite{ACP16,ACP17,MC12}, however, here is the first time that we 
prove the result for reconstructions of variable $u(x,t)$. Moreover, 
with the adaptation of our proof, we are able to obtain new 
conditions to guarantee the stability of the numerical method. 
These conditions are not obtained in previous works. 

	\textbf{Proof}. First, we fix 
	$\epsilon$ and obtain an ODE from (\ref{ODE}) 
	for variable $t$ as follows:
\begin{equation}
	{u}^\prime(x,t,\epsilon)=F_\epsilon({u}(x,t,\epsilon)), 
	\quad {u}(0,x)=u_0(x),\label{xxt}
	\end{equation}  
	Notice that $x$ is also a parameter in Eq. 
	$(\ref{xxt})$. 
	We define $F_\epsilon:L^\infty(\mathbb{S}^1)\longrightarrow L^\infty(\mathbb{S}^1)$ 
	as
	\begin{align}
	&F_\epsilon(\am{u}_\epsilon(x,t),x,t) = \frac{1}{\epsilon}
	\Big[\am{u}(x-\epsilon,t,\epsilon)f^+\left(\hat{\am{u}}\left(x-\frac{\epsilon}{2},t,\epsilon\right),x-\frac{\epsilon}{2},t\right)\notag\\
	&-u(x,t,\epsilon)f^+\left(\hat{u}\left(x+\frac{\epsilon}{2},t,\epsilon\right),x+\frac{\epsilon}{2},t\right)
	-u(x,t,\epsilon)f^-\left(\hat{u}\left(x-\frac{\epsilon}{2},t,\epsilon\right),x-\frac{\epsilon}{2},t\right)\notag\\
	&+u\left(x+{\epsilon},t,\epsilon\right)f^-\left(\hat{u}\left(x+\frac{\epsilon}{2},t,\epsilon\right),x+\frac{\epsilon}{2},t\right)
	\Big]\label{fluxeps}
	\end{align}
	Here $\hat{u}\left(x+{\epsilon/2},t,\epsilon\right)$ is obtained 
	from reconstruction $L(x_{-p},\dots,x_{p})$ 
	as in Eq. $(\ref{ffdse})$: 
	\begin{equation}
	\begin{split}
	\hat{u}(x+{\epsilon/2},t,\epsilon)=
	L(u(x-p\epsilon,t,\epsilon), \dots u(x+(p-1)
	\epsilon,t,\epsilon),u(x+p\epsilon,t,\epsilon)).
	\label{ffdse2}
	\end{split}
	\end{equation}   
	Since $f(\cdot)$ and reconstruction 
	$L(x_{-p},\dots,x_{p})$ are Lipschitz continuous,  so are $f(L(x_{-p},\dots,x_{p}))$. Since 
	flux $F_\epsilon$, defined in Eq. $(\ref{fluxeps})$, 
	is a combination of Lipschitz continuous (in a 
	remarkably simple way), we  get that $F_\epsilon$ 
	is also a Lipschitz function. Thus, 
	based on the classical theory of ODEs in Banach spaces, 
	in the Lipschitz continuous case, there is a local 
	solution to $t\in [0,\delta(\epsilon)]$ for 
	some $\delta(\epsilon)$ that depends on $\epsilon$. 
	For the global solution, since $f$ is bounded (because 
	$H$ is assumed to be Lipschitz continuous), 
	we can extend the solution to 
	$\delta(\epsilon)\to \infty$. From assumption 1, Eq. 
	$(\ref{h1})$, the Lipschitz constants of each 
	$F_\epsilon$ can be chosen uniformly on bounded 
	sets $L^\infty({\mathbb{S}^1})\times [0,\infty)$. To demonstrate
	the existence of the solution $(\ref{xxt})$ globally (in time), it suffices to prove that, 
	for fixed $\epsilon$, there exists a $c_\epsilon(t)<\infty$ 
	such that $||u(\cdot,t,\epsilon)||_\infty\leq c_\epsilon(t)<\infty$.  
	Here $c_\epsilon$ is a continuous function on $[0,\infty)$, which is not uniformly continuous in $\epsilon$.
	We also have that $H(u)$ is 
	a bounded function; hence, if $u\geq a>0$, then $f$ 
	is also bounded. 
	Let $M$ be defined as
	\begin{equation*}
	M=\sup_{\begin{array}{c}(u,x,t)\in\Omega\times[0,T]\times{\mathbb{S}^1}\end{array}}|f(u,x,t)|<\infty.
	\end{equation*}
	We have that $
	\displaystyle{|\partial_t u(x,t,\epsilon)|\leq \frac{4}{\epsilon}||
	u(\cdot,t,\epsilon)||_\infty M}$ \text{ where } $\displaystyle{
	||u(\cdot,t,\epsilon)||_\infty=\text{ess} \sup_{x\in {\mathbb{S}^1}}|u(x,t,\epsilon)|}$.
	Solving, we obtain:
	\begin{equation*}
	||u(\cdot,t,\epsilon)||_\infty\leq ||u_{0}(\cdot)
	||_\infty +\frac{4M}{\epsilon}\int_0^t||u(\cdot,\tau,\epsilon)||_\infty d\tau.
	\end{equation*}
	From Gr\"{o}nwall formula, we obtain that:
	\begin{equation}
	||u(\cdot,t,\epsilon)||_\infty\leq c_\epsilon (t), \quad \text{ where } \quad c_\epsilon(t)= ||u_{0}(\cdot)||_\infty\exp\left(\frac{4Mt}{\epsilon}\right).\label{ggr}
	\end{equation}
	Bound $(\ref{ggr})$ indicates the existence of 
	a global solution to the ODE $(\ref{ODE})$ for each 
	fixed $\epsilon$. However, note that the solutions 
	to the system of ODEs are not uniformly continuous 
	in $\epsilon$. To demonstrate that the 
	solutions to ODEs provide a weak asymptotic 
	solution for $(\ref{noLi})$, we will prove that 
	the solution is  $L^1$ 
	bounded uniformly with respect to $\epsilon$. 
	To do so, let $T>0$ for $t+dt\leq T$ and $dt>0$. 
	From the Taylor expansion with remainder term (for fixed 
	$\epsilon)$, it follows that $(\ref{ODE})$ can be written as
	\begin{align*}
	u(x,t+dt,\epsilon)=u_\epsilon&+\frac{dt}{\epsilon}
	\left[u_{_{\epsilon-1}}f^+\left(\hat{u}_{_{\epsilon-1/2}}\right)
	-u_{_{\epsilon}}f^+\left(\hat{u}_{_{\epsilon+1/2}}\right)
	-u_{_{\epsilon}}f^-\left(\hat{u}_{_{\epsilon-1/2}}\right)
	\right.\notag \\
	&\left.+u_{_{\epsilon+1}}f^-\left(\hat{u}_{_{\epsilon+1/2}}\right)\right]+dt \; r(x,t,\epsilon,dt),
	\end{align*} 
	where $||r(\cdot,t,\epsilon,dt)||_1\to 0$ (or $||r(\cdot,t,\epsilon,dt)||_\infty\to 0$), uniformly 
	in $t \in [0,T]$ and fixed $\epsilon$ (and not uniformly
	continuous in $\epsilon$), when $dt\to 0$. This 
	behavior results from the continuous differentiability 
	of map $t\longrightarrow u(\cdot, t,\epsilon)$, 
	$[0,\infty)\longrightarrow L^\infty ({\mathbb{S}^1})$ 
	for fixed $\epsilon$.  
	Since we are interested in obtaining the $L^1$ 
	bound, we take the absolute value,
	\begin{align}
	|u(x,t+dt,\epsilon)|\leq&|u_\epsilon|\left(1
	-\frac{dt}{\epsilon}(f^+\left(\hat{u}_{_{\epsilon+1/2}}\right)
	+f^-\left(\hat{u}_{_{\epsilon-1/2}}\right)\right)+
	\frac{dt}{\epsilon}\left[|u_{_{\epsilon-1}}|
	f^+\left(\hat{u}_{_{\epsilon-1/2}}\right)\right]\notag \\
	&+\frac{dt}{\epsilon}\left[|u_{_{\epsilon+1}}|
	f^-\left(\hat{u}_{_{\epsilon+1/2}}\right)\right]+dt|r(x,t,dt)|.\label{ODE2n}
	\end{align}
	 From the 
	definition of $f^+$ and $f^-$ in Eq. 
	(\ref{c2}), we get that if $(\ref{Cflcon})$ is satisfied, 
	then $(\ref{ODE2n})$ is true. Notice that $(\ref{Cflcon})$ is always valid if the CFL condition $(\ref{CFL})$ is satisfied. This proves that 
	the  condition $(\ref{CFL})$ provides the 
	method with stability because by integrating $(\ref{ODE2n})$ 
	and the appropriate translations of $\pm \epsilon$, we obtain
	\begin{equation}
	||u(\cdot,t+dt,\epsilon||_1=\int_{\mathbb{S}^1}|u(x,t+dt,\epsilon)|dx\leq 
	\int_{\mathbb{S}^1}|u(x,t,\epsilon)|dx+dt \; r_1(t,\epsilon,dt).\label{fdsee}
	\end{equation}
	Here the remainder value $r_1(t,\epsilon,dt)=\int_{\mathbb{S}} |r(\cdot, t,\epsilon,dt)|dx$ is 
	bounded, and $r_1(t,\epsilon,dt)\to 0$ when $dt\to 0$, 
	uniformly in $t\in [0,T]$ for each fixed $\epsilon$. 
	Notice that, for each $T>0$ given, we can divide 
	interval $[0,T]$ into $n$ subintervals 
	$[jdt_n,(j+1)dt_n]$, where 
	$dt_n=\frac{T}{n}$ and $0\leq j \leq n-1$. 
	Applying this in $(\ref{fdsee})$, we get
	\begin{equation*}
	\int_{\mathbb{S}^1}|u(x,T,\epsilon)|dx\leq 
	\int_{\mathbb{S}^1}|u(x,T-dt_n,\epsilon)|dx+dt_n \; r_1(T-dt,\epsilon,dt).
	\end{equation*} 
	Applying recursively for all intervals, we obtain
	\begin{equation*}
	\int_{\mathbb{S}^1}|u(x,T,\epsilon)|dx\leq 
	\int_{\mathbb{S}^1}|u_0(x)|dx+dt_n \; \sum_{i=1}^n r_1(T-idt,\epsilon,dt).
	\end{equation*} 
	Notice that
	\begin{equation*}
	dt_n \; \left|\sum_{i=1}^n r_1(T-idt,\epsilon,dt)\right|\leq 
	\frac{T}{n}n\max_{i}|r_1(T-idt,\epsilon,dt)|=T\max_{i}|r_1(T-idt,\epsilon,dt)|.
	\end{equation*}
		Thus, taking the limit $dt\longrightarrow 0$ and using  $r_1(t,\epsilon,dt)\longrightarrow 0$, when $dt\longrightarrow 0$ we obtain:
	\begin{equation}
	||u(\cdot,T,\epsilon)||_1=\int_{\mathbb{S}^1}|u(x,T,\epsilon)|dx\leq 
	\int_{\mathbb{S}^1}|u_0(x)|dx=||u_0(\cdot)||_1,\label{fdseeb}
	\end{equation}
	which gives us the $L^1$ uniform bounds in $\epsilon$. 
	
	To complete the proof of the proposition, we 
	define integral $I$ as
	\begin{equation*}
	I=\int_{{\mathbb{S}^1}}
	\left(\frac{1}{\epsilon}
	\left[u_{_{\epsilon-1}}f^+\left(\hat{u}_{_{\epsilon-1/2}}\right)
	-u_{_{\epsilon}}f^+\left(\hat{u}_{_{\epsilon+1/2}}\right)
	-u_{_{\epsilon}}f^-\left(\hat{u}_{_{\epsilon-1/2}}\right)+u_{_{\epsilon+1}}f^-\left(\hat{u}_{_{\epsilon+1/2}}\right)\right]\psi(x)-H(u_\epsilon)\psi_x(x)\right)dx
	\end{equation*}
	for a test function $\psi(x)$ in the sense of $(\ref{weaksol})$.
	Changing the order in the integration variable, we obtain
\begin{align*}
	I=\int_{{\mathbb{S}^1}}
	\left(u_{_{\epsilon}}f^+\left(\hat{u}_{_{\epsilon+1/2}}\right)\frac{\psi(x+\epsilon)-\psi(x)}{\epsilon}-
	u_{_{\epsilon}}f^-\left(\hat{u}_{_{\epsilon-1/2}}\right)\frac{\psi(x)-\psi(x-\epsilon)}{\epsilon} -H(u_\epsilon)\psi_x(x)\right)dx.
	\end{align*}

Using that 
	$({\psi(x+\epsilon)-\psi(x)})/{\epsilon}=\psi_x(x)+\mathbb{O}(\epsilon)$, $({\psi(x)-\psi(x-\epsilon)})/{\epsilon}=\psi_x(x)+\mathbb{O}(\epsilon)$,    $f=f^+-f^-$,  $uf=H$ and that $u(x,t,\epsilon)$ is bounded,  the integral $I$ satisfies:
	\begin{align*}
	I=\int_{{\mathbb{S}^1}}
	\left(u_{_{\epsilon}}\left(f^+\left(\hat{u}_{_{\epsilon+1/2}}\right)-f^+(u_\epsilon)\right)\psi_x(x)-
	u_{_{\epsilon}}\left(f^-\left(\hat{u}_{_{\epsilon-1/2}}\right)-f^-(u_\epsilon)\right)\psi_x(x)\right)dx+\mathbb{O}(\epsilon).
	\end{align*}	
	The function $u(x,t,\epsilon)$ is not necessarily 
	continuous, but their discontinuities are in a  set 
	of null measure, thus using that $f$, $f^+$ and $f^-$ are Lipschitiz functions, we have that (except in a set of null measure)
	\begin{equation}
	\left|f^+\left(\hat{u}_{_{\epsilon+1/2}}\right)-f^+(u_\epsilon)\right|\leq \bar{K}\left|\hat{u}_{_{\epsilon+1/2}}-u_\epsilon\right| \text{ and } \left|f^-\left(\hat{u}_{_{\epsilon-1/2}}\right)-f^-(u_\epsilon)\right|\leq \bar{K}\left|\hat{u}_{_{\epsilon-1/2}}-u_\epsilon\right|.\label{catr}
	\end{equation}
	For $\bar{K}$ the Lipschitz constant for function $f$. Using $(\ref{catr})$, we have that (except in a set of null measure)
	\begin{align*}
	|I|\leq \bar{K}\int_{{\mathbb{S}^1}}
	|u_{_{\epsilon}}|\left(\left|\hat{u}_{_{\epsilon+1/2}}-u_\epsilon\right|+\left|\hat{u}_{_{\epsilon-1/2}}-u_\epsilon\right|\right)\psi_x(x)dx+\mathbb{O}(\epsilon).
	\end{align*}	
	Taking the limit of $\epsilon\longrightarrow 0$,  using $(\ref{fds})$ and that $u(x,t,\epsilon)$ is bounded, we have that $|I|\longrightarrow 0$, that implies that $I\longrightarrow 0$ and thus $u(x,t,\epsilon)$ satisfies Eq. $(\ref{weaksol})$ and the proof is concluded. $\quad \square$
	
	\;
	
	\begin{remark}
	\label{remarkt}
		Whenever necessary, we can replace functions 
		$f^\pm$ with $k(t)+f^\pm$, where $k(t)$ is a positive 
		function that is large enough and bounded on any interval 
		$[0,T]$, so that functions $u\to 
		u(k(t)+f^\pm(\hat{u}, x, t))$ are strictly
		increasing on $\mathbb{R}$. Here $\hat{u}$ is a reconstruction of speed $u$.
		 The fact that 
		$u(k(t)+f^\pm({u}, x, t))$ is increasing  
		is proven in \cite{ACP16} (and references cited therein). 
		Notice that this substitution does not change the 
		proof of Proposition $\ref{prop1}$. As noticed by the authors 
		in \cite{ACP16}, function $k(t)$ only produces a vanishing viscosity solution. 
		
		In the case that there exists a reconstruction of variable $u$, the proof can be omitted because it is similar to that presented in \cite{ACP16} and because $f$ and the reconstruction of $u$ are Lipschitz continuous. In addition, one can demonstrate that, for a large enough $k(t)$, function $u(k(t)+f(\cdot))$ is independently increasing in the argument of $f$. This is important to prove the maximum principle.  
	\end{remark}
	\label{remark1}
	As the final step to establish the convergence of the numerical method, 
	we demonstrate that $(\ref{projection})$ can be written 
	as a particular case of $(\ref{ODE})$, taking 
	$\epsilon= h$.
Since our numerical method is a discrete numerical method, and we are using the theory for ODEs to prove the convergence of the solution, we need to show that if we take the limit of $\Delta t$ in $(\ref{projection})$ we obtain $(\ref{ODE})$.  
	\begin{proposition}
		\label{prop2}{
	 Consider the numerical method $(\ref{projection})$, taking the limit of $\Delta t \longrightarrow 0$, we obtain the ODE:}
		\begin{equation}{
		U_t=\frac{1}{h}\left(U_{j-1}f^+_{j-\frac{1}{2}}-U_{j}
		f^-_{j-\frac{1}{2}}-U_{j}f^+_{j+\frac{1}{2}}
		+U_{j+1}f^-_{j+\frac{1}{2}}\right).\label{ODEme}}
		\end{equation}
		{In this case, we say that the scheme $(\ref{projection})$ is compatible with the ODE $(\ref{ODEme})$, where $U_{j} = U_j(t)$ and $U_t = \lim_{\Delta t \to 0} \frac{U_j^{n+1}-U_j^{n}}{\Delta t}$}
	\end{proposition}
	\textbf{Proof}. The numerical scheme is given by 
	$(\ref{projection})$.  If we substitute Eq.  
	$(\ref{c2})$ in 
	Eq. $(\ref{projection})$, we obtain
	\begin{equation}
	U_{j}^{n+1} = \overline{U}_{j}^{n}+ 
	\frac{\Delta t}{h}\left( f^+(U_{j-\frac{1}{2}}) 
	\overline{U}_{j-1}^{n}- f^+(U_{j-\frac{1}{2}})
	\overline{U}_{j}^{n}- f^-(U_{j+\frac{1}{2}})
	\overline{U}_{j}^{n}+ f^-(U_{j+\frac{1}{2}})
	\overline{U}_{j+1}^{n} \right).
	\label{lere}
	\end{equation}
	Replacing $\overline{U}_{j}^{n}$, given by 
	Eq. $(\ref{eq32})$, in $(\ref{lere})$ reads   
	\begin{equation}
	U_{j}^{n+1} = 
	U_{j}^{n}\frac{h}{h^{n+1}_{j}}+ {\Delta t}\left( 
	f^+_{j-\frac{1}{2}}\frac{U_{j-1}^{n}}{h^{n+1}_{j-1}}
	-(f^+_{j-\frac{1}{2}}+ f^-_{j+\frac{1}{2}})
	\frac{U_{j}^{n}}{h^{n+1}_{j}}+ 
	f^-_{j+\frac{1}{2}}\frac{U_{j+1}^{n}}{h^{n+1}_{j+1}} \right).
	\label{eww}
	\end{equation}
	Using  $h^{n+1}_j$,  
	and 
	substituting this result in Eq. $(\ref{eww})$, we have
	\begin{equation*}\begin{split}
	U_{j}^{n+1} = &U_{j}^{n}
	\left(1-\Delta t\frac{f_{j+\frac{1}{2}}-
		f_{j-\frac{1}{2}}}{h^{n+1}_{j}}\right) + \\
	&\frac{\Delta t}{h}\left(f^+_{j-
		\frac{1}{2}}U_{j-1}^{n}\left(1+\Delta t
	\frac{f_{j-\frac{1}{2}}-f_{j-\frac{3}{2}}}{h^{n+1}_{j-1}}\right)- 
	f^+_{j-\frac{1}{2}}U_{j}^{n}\left(1+
	\Delta t\frac{f_{j+\frac{1}{2}}-
		f_{j-\frac{1}{2}}}{h^{n+1}_{j}}\right)-\right.\\
	&\left.-f^-_{j+\frac{1}{2}}U_{j}^{n}\left(1+
	\Delta t\frac{f_{j+\frac{1}{2}}-
		f_{j-\frac{1}{2}}}{h^{n+1}_{j}}\right)+ 
	f^-_{j+\frac{1}{2}}U_{j+1}^{n}\left(1+\Delta 
	t\frac{f_{j+\frac{3}{2}}-
		f_{j+\frac{1}{2}}}{h^{n+1}_{j+1}}\right) \right).
	\end{split}\end{equation*}  
	It can be written as
	\begin{equation}\begin{split}
	U_{j}^{n+1} - U^n_j &= -U_{j}^{n}\Delta t
	\left(\frac{f_{j+\frac{1}{2}}-f_{j-\frac{1}{2}}}{h^{n+1}_{j}}\right) 
	+\frac{\Delta t}{h}\left(f^+_{j-\frac{1}{2}}U_{j-1}^{n}- 
	(f^+_{j-\frac{1}{2}}+f^-_{j+\frac{1}{2}})U_{j}^{n}
	+f^-_{j+\frac{1}{2}}U_{j+1}^{n} \right) + o(\Delta t^2).
	\label{lere5}
	\end{split}\end{equation}
	Dividing both sides of $(\ref{lere5})$ by $\Delta t$, 
	using $f=f^+-f^-$, taking the limit of 
	$\Delta t \to 0$ in Eq. $(\ref{lere5})$, we obtain that the left side converges to $U_t$, moreover, $h_j^{n+1}$ converges to $h$ and Eq. $(\ref{lere5})$ leads to $(\ref{ODEme})$, i.e.,  the numerical method  $(\ref{projection})$  is compatible with $(\ref{ODEme})$. $\quad \square$
	
	Notice that Proposition $\ref{prop2}$ shows that 
	the numerical method is compatible with ODE 
	$(\ref{ODE})$ constructed in Proposition 
	$\ref{prop1}$, which is powerful for numerics.
	\subsection{Conditions for Total Variation Nonincreasing 
$TV_\epsilon(u(\cdot, t+dt,\epsilon))\leq TV_\epsilon(u(\cdot, t,\epsilon)$}
	\label{tvdsection}
	
	We now prove some further results regarding the scheme 
	described by ODEs  $(\ref{ODE})$. { More precisely we have that the total variation with respect to $x$ does not increase with time by using the weak asymptotic solution we have considered through the fully discrete La\-gran\-gian--Eu\-le\-rian scheme (\ref{eq32})--(\ref{projection}) for solving the one-dimensional scalar equation in (\ref{noLi}) after discretization in time of the ODEs  $(\ref{ODE})$.} 
	The first 
	property is that such scheme has a nonincreasing total variation that depends on $\epsilon$. This enables us to define a kind of total variation useful for this work.
	
	We say that a numerical scheme is $\epsilon$ Total Variation Nonincreasing, denoted as TVNI\textsubscript{$\epsilon$}, if
	$TV_\epsilon(u(\cdot, t+dt,\epsilon))\leq TV_\epsilon(u(\cdot, t,\epsilon),$
	where 
	\begin{equation}
	TV_\epsilon(u(\cdot, t,\epsilon))=
	\int_{{\mathbb{S}^1}}|u(x+\epsilon, t,\epsilon)-
	u(x, t,\epsilon)|dx.\label{tv}
	\end{equation}
	Note that the total variation for fixed $\epsilon$ can be obtained by $TV_\epsilon(u(\cdot, t+dt,\epsilon))/\epsilon$.
	
	Using a similar idea from Harten (see \cite{HART83}), we can prove the following result:

	\begin{lemma} 
		\label{lemar}
		If the numerical scheme can be written in the following semidiscrete form:
		\begin{equation}
		(u(x,t,\epsilon))_t=\frac{C_{x+\frac{\epsilon}{2}}\Delta_{\frac{\epsilon}{2}}u(x,t,\epsilon)-D_{x-\frac{\epsilon}{2}}\Delta_{-\frac{\epsilon}{2}}u(x,t,\epsilon)}{\epsilon},\label{uscheme}
		\end{equation}	
		with  $C_{x+\frac{\epsilon}{2}}$ and $D_{x-\frac{\epsilon}{2}}$ as the arbitrary values satisfying 
		\begin{align}
		&C_{x+\frac{\epsilon}{2}}\geq 0, \quad D_{x-\frac{\epsilon}{2}}\geq 0\text{ and } \quad  \frac{dt}{\epsilon}\left(C_{x+\frac{\epsilon}{2}}+ D_{x+\frac{\epsilon}{2}}\right)\leq 1,\label{cpos1}
		\end{align} 
		then the system is TVNI\textsubscript{$\epsilon$} and satisfies
		\begin{equation}
		TV_\epsilon(u(\cdot,t,\epsilon))\leq TV_\epsilon(u(\cdot,0), \quad \forall t\in[0,T].\label{lemv}
		\end{equation}	
	\end{lemma}
	In Eq. $(\ref{uscheme})$, we define
	\begin{equation}
	\Delta_{i\frac{\epsilon}{2}}u(x,t,\epsilon)=
	u\left(x+i\frac{\epsilon}{2},t,\epsilon\right)-u\left(x-i\frac{\epsilon}{2},t,\epsilon\right)\quad \text{for } i\in \mathbb{Z}.\label{deltad}
	\end{equation} 
	Notice that by using $(\ref{deltad})$, we can define $TV_\epsilon(u(\cdot, t,\epsilon))=
	\int_{{\mathbb{S}^1}}|\Delta_{\frac{\epsilon}{2}}u(x,t,\epsilon)|dx.\label{tv2}
	$
	
	\textbf{Proof of Lemma $\ref{lemar}$.} 
	From the Taylor expansion with remainder term (for fixed $\epsilon$), we can write $(\ref{uscheme})$ as
	\begin{equation}
	u(x,t+
	dt,\epsilon)=u_\epsilon +\frac{dt}{\epsilon}
	\left(C_{x+\frac{\epsilon}{2}}\Delta_{\frac{\epsilon}{2}}u(x,t,\epsilon)-D_{x-\frac{\epsilon}{2}}\Delta_{-\frac{\epsilon}{2}}u(x,t,\epsilon)\right)+dt r(x,t,\epsilon,dt),\label{uscheme2}
	\end{equation}	
	for $||r(\cdot,t,\epsilon,dt)||_1\longrightarrow 0$ ({or $||r(\cdot,t,\epsilon,dt)||_\infty\longrightarrow 0$}), uniformly in $t\in[0,T]$ and fixed $\epsilon$, when $dt\longrightarrow 0$.
	
	Subtracting $u(x,t+
	dt,\epsilon)$ from $u(x+\epsilon,t+
	dt,\epsilon)$, both given by $(\ref{uscheme2})$, we obtain
	\begin{equation*}
	\begin{split}
	\Delta_{\frac{\epsilon}{2}} u(x,t+
	dt,\epsilon)=&\Delta_{\frac{\epsilon}{2}}u(x,t,\epsilon)\left(1-\frac{dt}{\epsilon}\left(D_{x+\frac{\epsilon}{2}}+C_{x+\frac{\epsilon}{2}}\right)\right)+\\
	&+\frac{dt}{\epsilon}
	D_{x-\frac{\epsilon}{2}}\Delta_{\epsilon-\frac{1}{2}}u(x,t,\epsilon) +\frac{dt}{\epsilon}C_{x+\frac{3\epsilon}{2}}\Delta_{\epsilon+\frac{3}{2}}u(x,t,\epsilon) +dt \Delta_{\frac{\epsilon}{2}} r(x,t,\epsilon,dt),
	\end{split}
	\end{equation*}
	where $\Delta_{\frac{\epsilon}{2}} r(x,t,\epsilon,dt)=r(x+\epsilon,t,\epsilon,dt)-r(x,t,\epsilon,dt)$.
	Due to $(\ref{cpos1})$, all coefficients are nonnegative; therefore, we have that
	\begin{align}
	&|\Delta_{\frac{\epsilon}{2}} u(x,t+
	dt,\epsilon)|\leq |\Delta_{\frac{\epsilon}{2}}u(x,t,\epsilon)|\left(1-\frac{dt}{\epsilon}\left(D_{x+\frac{\epsilon}{2}}+C_{x+\frac{\epsilon}{2}}\right)\right) \notag\\
	&+\frac{dt}{\epsilon}
	D_{x-\frac{\epsilon}{2}}|\Delta_{\epsilon-\frac{1}{2}}u(x,t,\epsilon)|+\frac{dt}{\epsilon}C_{x+\frac{3\epsilon}{2}}|\Delta_{\frac{3\epsilon}{2}}u(x,t,\epsilon)| +dt |\Delta_{\frac{\epsilon}{2}} r(x,t,\epsilon,dt)|.\label{uscheme4}
	\end{align}
	Integrating $(\ref{uscheme4})$ in $x\in{\mathbb{S}^1}$, we notice that, due to translations of $\pm \epsilon$, there are two-by-two simplifications of the terms of $(\ref{uscheme4})$, as in Eq. $(\ref{fdsee})$. Thus, we get 
	\begin{align*}
	&TV_\epsilon(u(x,t+dt,\epsilon)=
	\int_{{\mathbb{S}^1}}|\Delta_{\frac{\epsilon}{2}} u(x,t+
	dt,\epsilon)|dx+dt\int_{{\mathbb{S}^1}}|\Delta_{\frac{\epsilon}{2}} r(x,t,\epsilon,dt)|dx\notag\\
	&\leq \int_{{\mathbb{S}^1}}|\Delta_{\frac{\epsilon}{2}} u(x,t,\epsilon)|dx+dt\int_{{\mathbb{S}^1}}|\Delta_{\frac{\epsilon}{2}} r(x,t,\epsilon,dt)|dx  =TV_\epsilon(u(x,t,\epsilon) +dt\int_{{\mathbb{S}^1}}|\Delta_{\frac{\epsilon}{2}} r(x,t,\epsilon,dt)|dx.
	\end{align*} 
	Since $\int_{{\mathbb{S}^1}}|\Delta_{\frac{\epsilon}{2}} r(x,t,\epsilon,dt)|dx\leq 2 \int_{{\mathbb{S}^1}}| r(x,t,\epsilon,dt)|dx=||r(x,t,\epsilon,dt)||_1\rightarrow 0$ when $dt\rightarrow 0$ and using an argument similar used to prove $(\ref{fdseeb})$, we obtain $(\ref{lemv})$. \quad $\square$
	
	To demonstrate that our scheme satisfies the TVNI\textsubscript{$\epsilon$} property, we must prove that our method satisfies the hypothesis of Lemma $\ref{lemar}$. For this purpose, we use the mean value theorem. However, since $f^+$ and $f^-$ are only continuous but not differentiable in some points, we need to extend the mean value theorem to this more general case. In Appendix $\ref{apendicef}$, we provide a smooth extension for the derivatives of 
	$f^+$ and $f^-$,  denoted as $\hat{f}^+$ and $\hat{f}^-$.
	Using these two functions, we can prove the following result:
	\begin{proposition}
		\label{Prop5}
		Let us assume that the reconstruction of $\hat{u}_{_{\epsilon\pm 1/2}}$ satisfies
		\begin{equation*}
		\hat{u}_{_{\epsilon+1/2}}-\hat{u}_{_{\epsilon-1/2}}=L_{1,\epsilon}(u_{_{\epsilon}}-u_{_{\epsilon-1}})+L_{2,\epsilon}(u_{_{\epsilon+1}}-u_{_{\epsilon}})
		\end{equation*}
		for some functions $L_{1,\epsilon}$ and $L_{2,\epsilon}$. If we define $f^+=\max(f,0)+k$ and $f^-=\max(-f,0)+k$, with $k$ satisfying \begin{equation}
		k=\max_{u\in\Omega}(|uf^\prime(u)|)|\max_{\epsilon}(L_{1,\epsilon},L_{2,\epsilon})|,\label{ksatis}
		\end{equation}
		then scheme $(\ref{ODE})$ is TVNI\textsubscript{$\epsilon$} if we take $dt/\epsilon$ satisfying 
		\begin{equation}
		\frac{d t}{\epsilon}\left(2k+\max_{u\in\Omega}(|(uf(u))^\prime|(L_1+L_2)\right)\leq  1, \label{etd4}
		\end{equation}
		where
		$
		L_1=\max_{\epsilon}L_{1,\epsilon} \text{ and } L_2=\max_{\epsilon}L_{2,\epsilon}.
		$

	\end{proposition}
	
	\textbf{Proof.}
	We write the Right-Hand Side (RHS) of Eq. $(\ref{ODE})$ (disregarding $1/\epsilon$) as
	\begin{align*}
	RHS=&u_{_{\epsilon-1}}f^+\left(\hat{u}_{_{\epsilon-1/2}}\right)-u_{_{\epsilon}}f^+\left(\hat{u}_{{\epsilon-1/2}}\right)+u_{_{\epsilon}}f^+\left(\hat{u}_{_{\epsilon-1/2}}\right)
	-u_{_{\epsilon}}f^+\left(\hat{u}_{_{\epsilon+1/2}}\right)\notag \\
	-&u_{_{\epsilon}}f^-\left(\hat{u}_{_{\epsilon-1/2}}\right)-u_{_{\epsilon}}f^-\left(\hat{u}_{_{\epsilon+1/2}}\right)+u_{_{\epsilon}}f^-\left(\hat{u}_{_{\epsilon+1/2}}\right)
	+u_{_{\epsilon+1}}f^-\left(\hat{u}_{_{\epsilon+1/2}}\right),
	\end{align*}  
	which can be written as
	\begin{align}
	RHS=&f^+\left(\hat{u}_{_{\epsilon-1/2}}\right)(u_{_{\epsilon-1}}-u_{_{\epsilon}}) +u_{_{\epsilon}}\left(f^+\left(\hat{u}_{_{\epsilon-1/2}}\right)
	-f^+\left(\hat{u}_{_{\epsilon+1/2}}\right)\right)\notag \\
	+&
	f^-\left(\hat{u}_{_{\epsilon+1/2}}\right)\left(u_{_{\epsilon+1}}-u_{_{\epsilon}}\right)
	+u_{_{\epsilon}}\left(f^-\left(\hat{u}_{_{\epsilon+1/2}}\right)-f^-\left(\hat{u}_{_{\epsilon-1/2}}\right)\right).\label{fst}
	\end{align}
	Using $(\ref{deltad})$; the mean value theorem; and the extensions for the derivative of ${f}^+$, denoted as $(\hat{f}^+)^\prime$, and for that of $f^-$, denoted as $(\hat{f}^+)^\prime-f^\prime$ (because $f^-=f^+-f$), we can write $(\ref{fst})$ as
	\begin{align}
	&RHS=-\left[ \Delta_{-\frac{\epsilon}{2}}u f^+\left(\hat{u}_{_{\epsilon-1/2}}\right)+u_{_{\epsilon}}(\hat{f}^+)^\prime(\xi_{\epsilon})\left(
	L_{1,\epsilon}\Delta_{-\frac{\epsilon}{2}}u+L_{2,\epsilon}\Delta_{\frac{\epsilon}{2}}u\right)\right] \notag \\
	&+\left[ \Delta_{\frac{\epsilon}{2}}u f^-\left(\hat{u}_{_{\epsilon+1/2}}\right)+u_{_{\epsilon}}((\hat{f}^+)^\prime(\xi_{\epsilon})-f^\prime(\eta_{\epsilon}))\left(
	L_{1,\epsilon}\Delta_{-\frac{\epsilon}{2}}u+L_{2,\epsilon}\Delta_{\frac{\epsilon}{2}}u\right)\right] \label{fst2}
	\end{align}
	where $\xi_\epsilon$ and $\eta_\epsilon$ are values between $\hat{u}_{_{\epsilon-1/2}}$ and $\hat{u}_{_{\epsilon+1/2}}$, and functions 
	$L_{1,\epsilon}$ and $L_{2,\epsilon}$ depend on the reconstruction of $\hat{u}_{_{\epsilon+1/2}}$. Here, we assume that we can estimate such reconstruction between a state $u_{_{\epsilon}}$ using $u_{_{\epsilon}}$, $u_{_{\epsilon-1}}$, and $u_{_{\epsilon+1}}$, even if the original reconstruction depends on more points. 
	Rearranging $(\ref{fst2})$, we get
	\begin{align*}
	RHS=& -\Delta_{-\frac{\epsilon}{2}}u \left[ f^+\left(\hat{u}_{_{\epsilon-1/2}}\right)+u_{_{\epsilon}}f^\prime(\eta_{\epsilon})
	L_{1,\epsilon}\right]+\Delta_{\frac{\epsilon}{2}}u \left[ f^-\left(\hat{u}_{_{\epsilon+1/2}}\right)-u_{_{\epsilon}}f^\prime(\eta_{\epsilon})
	L_{2,\epsilon}\right].
	\end{align*}
	Notice that the estimate does not depend on the extension of the derivative of $f^+$, which is canceled, thus obtaining a result that only depends on the derivative of $f$.
	
	Using Lemma $\ref{lemar}$, we observe that scheme $(\ref{ODE})$ is TVNI\textsubscript{$\epsilon$} if 
	\begin{equation}
	f^+\left(\hat{u}_{_{\epsilon-1/2}}\right)+u_{_{\epsilon}}f^\prime(\eta_{\epsilon})
	L_{1,\epsilon}\geq 0,\quad \text{ and }\quad 
	f^-\left(\hat{u}_{_{\epsilon+1/2}}\right)-u_{_{\epsilon}}f^\prime(\eta_{\epsilon})
	L_{2,\epsilon}\geq 0 \label{etd}
	\end{equation}
	and 
	\begin{equation}
	\frac{d t}{\epsilon}\left(f^+\left(\hat{u}_{_{\epsilon+1/2}}\right)+u_{_{\epsilon+1}}f^\prime(\eta_{\epsilon+1})
	L_{1,\epsilon+1}+ 
	f^-\left(\hat{u}_{_{\epsilon+1/2}}\right)-u_{_{\epsilon}}f^\prime(\eta_{\epsilon})
	L_{2,\epsilon}\right)\leq  1. \label{etd2}
	\end{equation}
	Condition $(\ref{etd})$ can be satisfied because it is possible to take $f^+=\max(0,f)+k$ for a positive constant $k$; since $f=f^+-f^-$, this choice does not change function $f$. For example, given that the minimum value for $f^+$ is 0, we can choose $k$ satisfying $(\ref{ksatis})$.	Condition $(\ref{etd2})$ can be rewritten as
	\begin{equation}
	\frac{d t}{\epsilon}\left(f\left(\hat{u}_{_{\epsilon+1/2}}\right)+2k+ u_{_{\epsilon+1}}f^\prime(\eta_{\epsilon+1})
	L_{1,\epsilon+1}-u_{_{\epsilon}}f^\prime(\eta_{\epsilon})
	L_{2,\epsilon}\right)\leq  1. \label{etd3}
	\end{equation}
	Note that ${u_{_{\epsilon+1}}f^\prime(\eta_{\epsilon+1})+f\left(\hat{u}_{_{\epsilon+1/2}}\right)}$ is close to the derivative of $uf(u)$; therefore, we can estimate $(\ref{etd3})$ satisfying $(\ref{etd4})$, and scheme $(\ref{ODE})$ is $TVNI_{\epsilon}$. $\quad \quad \square$
	\begin{example}
		For the reconstruction given by Eq. $(\ref{umeio})$ and using the MinMod slope limiter, given by Eqs. $(\ref{mm2})$--$(\ref{mm})$, for the derivative of $U^\prime_j$, we have that, for ${\Delta\hat{u}_{_{\epsilon}}=\hat{u}_{_{\epsilon+1/2}}-\hat{u}_{_{\epsilon-1/2}}}$ (substituting $j$ by $\epsilon$ in Eq. $(\ref{umeio})$),
		\begin{equation}
		\Delta\hat{u}_{_{\epsilon}}=\hat{u}_{_{\epsilon+1/2}}-\hat{u}_{_{\epsilon-1/2}}=\frac{u_{_{\epsilon+1}}+u_{_{\epsilon}}}{2}-\frac{u_{_{\epsilon}}+u_{_{\epsilon-1}}}{2}+\frac{u_{_{\epsilon+1}}^\prime-u_{_{\epsilon}}^\prime}{8}-\frac{u_{_{\epsilon}}^\prime-u_{_{\epsilon-1}}^\prime}{8}.\label{rhs1}
		\end{equation}
		We can write the two first terms of the RHS of Eq. $(\ref{rhs1})$ as
		\begin{equation*}
		\frac{u_{_{\epsilon+1}}+u_{_{\epsilon}}}{2}-\frac{u_{_{\epsilon}}+u_{_{\epsilon-1}}}{2}=\frac{u_{_{\epsilon+1}}-u_{_{\epsilon}}}{2}+\frac{u_{_{\epsilon}}-u_{_{\epsilon-1}}}{2}=\frac{\Delta u_{\frac{\epsilon}{2}}+\Delta u_{-\frac{\epsilon}{2}}}{2}.
		\end{equation*}
		For the derivatives and using the MinMod limiter, we know that
		\begin{align*}
		u_{_{\epsilon}}^\prime=\frac{1}{2}\left(\sgn\left(\Delta u_{\frac{\epsilon}{2}}\right)+\sgn\left(\Delta u_{-2\frac{\epsilon}{2}}\right)\right)\min\left(\left|\Delta u_{\frac{\epsilon}{2}}\right|,\left|\Delta u_{-\frac{\epsilon}{2}}\right|\right)=\theta_{1,\epsilon}\Delta u_{\frac{\epsilon}{2}} \quad\text{ or }\quad
		\theta_{2,\epsilon}\Delta u_{-\frac{\epsilon}{2}}.
		\end{align*}
		Choosing $\theta_{1,\epsilon}$ or $\theta_{2,\epsilon}$
		depends on the method. However, notice that these values satisfy
		$0\leq \theta_{1,\epsilon},\;\theta_{2,\epsilon}\leq 1$. We write the two last terms of the RHS of Eq. $(\ref{rhs1})$ as
		\begin{equation*}
		\frac{\theta_{2,\epsilon+1}\Delta u_{\frac{\epsilon}{2}}-\theta_{1,\epsilon}\Delta u_{\frac{\epsilon}{2}}-\left(\theta_{2,\epsilon}\Delta u_{-\frac{\epsilon}{2}}-\theta_{1,\epsilon-1}\Delta u_{-\frac{\epsilon}{2}}\right)}{8}=\frac{(\theta_{2,\epsilon+1}-\theta_{1,\epsilon})\Delta u_{\frac{\epsilon}{2}}-(\theta_{2,\epsilon}-\theta_{1,\epsilon-1})\Delta u_{-\frac{\epsilon}{2}}}{8}.
		\end{equation*}
		Note that $-1\leq\theta_{2,\epsilon}-\theta_{1,\epsilon-1}\leq 1$.
		Then, $L_{1,\epsilon}$ and $L_{2,\epsilon}$ can be written as $
		\displaystyle{L_{1,\epsilon}=\frac{1}{2}+\frac{(\theta_{2,\epsilon+1}-\theta_{1,\epsilon})}{8}}$,\linebreak$\displaystyle{\quad L_{2,\epsilon}=\frac{1}{2}-\frac{(\theta_{2,\epsilon}-\theta_{1,\epsilon-1})}{8}}.$
		We  estimate $3/8\leq L_{1,\epsilon}, L_{2,\epsilon}\leq 5/8$. Here, $k$, given by Eq. $(\ref{ksatis})$, can be written as 
		$\displaystyle{k=\max_{u\in\Omega}(|uf^\prime(u)|)|\frac{5}{8}}.$
		And estimate $(\ref{etd4})$ satisfies
		\begin{equation*}
		\frac{d t}{\epsilon}\left(2k+\frac{5}{4}\max_{u\in\Omega}(|(uf(u))^\prime|\right)=\frac{d t}{\epsilon}\left(\frac{5}{2}\max_{u\in\Omega}(|(uf(u))^\prime|\right)\leq  1. 
		\end{equation*} 
	\end{example}
	
	Using a similar idea, one can obtain an estimate for $L_{1,\epsilon}$ and $L_{2,\epsilon}$ for the other slope limiters presented in Section $\ref{slope}$.

\subsection{The maximum principle and the entropy solution}
 \label{maximumsec}
 
 In this section, we demonstrate that scheme $(\ref{ODE})$ leads to a solution that satisfies the maximum principle and the Kruzhkov entropy condition. We first prove the maximum principle with ideas similar to those reported in \cite{ACP16}. In this case, $f^+=\max(f,0)+k$ and $f^-=\max(-f,0)+k$ must be used in such way that both $uf^+(\hat{u})$ and $uf^-(\hat{u})$ are increasing functions, for all $\hat{u}$, as stated in Remark $\ref{remarkt}$. Moreover, we can choose a large enough $k$ such that $uf^\pm(\cdot)$ is independently increasing in the argument of functions $f^\pm$. This fact is useful to prove the maximum principle. According to the numerical experiments, if this condition is not satisfied, the maximum principle is not satisfied either. In the next proposition we denote  $u_0(x,\epsilon)$ as a continuous approximation of $u_0(x)$, the initial data for $(\ref{noLi})$, and we state our result as:

 \begin{proposition}
 \label{propmax}
 Let $k$ be large enough such that $uf^+(\cdot)$ and $uf^-(\cdot)$ are increasing functions. Then, any local solution on $[0,T)$, for $T>0$, of $(\ref{noLi})$ using scheme $(\ref{ODE})$ takes its values between range $[\min_{x\in{\mathbb{S}^1}}u_0(x),\max_{x\in{\mathbb{S}^1}}u_0(x)]$. 
 \end{proposition}

{ The proof of Proposition $\ref{propmax}$ follows the same steps from the maximum principle lemma in \cite{ACP16} (page 15). However, in this work we adapted such a proof to the La\-gran\-gian--Eu\-le\-rian scheme with reconstruction $(\ref{ODE})$. }

 \textbf{Proof.} We first consider $x\in {\mathbb{S}^1}$. Also, we consider values $\epsilon$ so that $\{n\epsilon\}_{n\in\mathbb{Z}}$ forms a dense set in ${\mathbb{S}^1}$. By contradiction, we assume that there exists a $\epsilon_0>0$ satisfying, for $T>0$,
\begin{equation}
\sup_{x\in{\mathbb{S}^1}}u(x,t,\epsilon_0)>\sup_{x\in{\mathbb{S}^1}}u_0(x,\epsilon_0)
\quad \text{ for some } t\in[0,T].\label{sup2t}
\end{equation}
{Since $u_0(x,\epsilon)$ is continuous, we can choose a small enough $\epsilon_0$ e $T>0$ so that $\{u(x,t,\epsilon_0)\} \subset [\min_{x\in{\mathbb{S}^1}}u_0(x,\epsilon)-\eta,\max_{x\in{\mathbb{S}^1}}u_0(x,\epsilon)+\eta]$.} Given that $u_0(x,\epsilon_0)$ is smooth,  solution $u(x,t,\epsilon_0)$ from Eq. $(\ref{ODE})$ is also smooth because this space can be considered a Banach space using the $L^\infty$ norm. Thus, there exists $t_0$, $x_0$ such that $\sup_{x\in{\mathbb{S}^1}}u(x,t,\epsilon_0)=u(x_0,t_0,\epsilon_0)$. Since $(t_0,x_0)$ is a maximum, solution $u(x,t,\epsilon_0)$ satisfies
\begin{equation}
    \partial_tu(x_0,t_0,\epsilon_0)\geq 0.\label{dtposi}
\end{equation}
Moreover, if we use scheme $(\ref{ODE})$, we obtain
\begin{align}
	\partial_t u(x_0,t_0,\epsilon_0)&=\frac{1}{\epsilon_0}
	\Big\{u(x_0-\epsilon_0,t_0,\epsilon_0)f^+\left(\hat{u}\left(x_0-\frac{\epsilon_0}{2},t_0,\epsilon_0\right)\right)\notag\\
	&-	u(x_0,t_0,\epsilon_0)f^+\left(\hat{u}\left(x_0+\frac{\epsilon}{2},t_0,\epsilon_0\right)\right)-	u(x_0,t_0,\epsilon_0)f^-\left(\hat{u}\left(x_0-\frac{\epsilon_0}{2},t_0,\epsilon_0\right)\right) \notag \\
	&+u(x_0+\epsilon_0,t_0,\epsilon_0)f^-\left(\hat{u}\left(x_0+\frac{\epsilon_0}{2},t_0,\epsilon_0\right)\right)
	\Big\}.\label{ODE2t}
	\end{align}
	Since $uf^\pm(\cdot)$ are increasing, $u(x_0-\epsilon_0,t_0,\epsilon_0)\leq u(x_0 ,t_0,\epsilon_0)$, and $u(x_0+\epsilon_0,t_0,\epsilon_0)\leq u(x_0 ,t_0,\epsilon_0)$ we have that
	\begin{equation*}
u(x_0-\epsilon_0,t_0,\epsilon_0)f^+\left(\hat{u}\left(x_0-\frac{\epsilon_0}{2},t_0,\epsilon_0\right)\right)\leq u(x_0,t_0,\epsilon_0)f^+\left(\hat{u}\left(x_0+\frac{\epsilon_0}{2},t_0,\epsilon_0\right)\right) 
	\end{equation*}
	and
	\begin{equation*}
	u(x_0+\epsilon_0,t_0,\epsilon_0)f^-\left(\hat{u}\left(x_0+\frac{\epsilon_0}{2},t_0,\epsilon_0\right)\right)
\leq u(x_0,t_0,\epsilon_0)f^-\left(\hat{u}\left(x_0-\frac{\epsilon_0}{2},t_0,\epsilon_0\right)\right).
	\end{equation*}
	Therefore, from Eq. $(\ref{ODE2t})$, we have that
	\begin{equation}
	\partial_tu(x_0,t_0,\epsilon_0)\leq 0. \label{dtposi2}
	\end{equation}
	From inequalities $(\ref{dtposi})$ and $(\ref{dtposi2})$, we get that 
	$\partial_tu(x_0,t_0,\epsilon_0)=0$. Thus, the second member of $(\ref{ODE2t})$ is null. Since function $uf^\pm(\cdot)$ are increasing, it means that
	$u(x_0-\epsilon_0,t_0,\epsilon_0)= u(x_0+\epsilon_0,t_0,\epsilon_0)=u(x_0,t_0,\epsilon_0)$, which, by recursion, leads to $u(x_0+n\epsilon_0,t_0,\epsilon_0)= u(x_0,t_0,\epsilon_0)$ for all $n$. In other words, $u$ is constant because $u$ is (at least) continuous and  $\mathbb{N}\epsilon_0$ is dense in ${\mathbb{S}^1}$ module 1 (since $\epsilon_0$ is taken as irrational). From { ODE} $(\ref{ODE})$, $u$ is constant, and the solution is trivial, leading to a contraction by the assumption. The same argument can be used by substituting $\sup$ by $\inf$ in Eq. $(\ref{sup2t})$, and the proof is completed. $\quad \square$.
	
	The next step of our construction is to prove that the proposed scheme satisfies some kind of entropy solution. In this work, we use Kruzhkov entropy solution.
	We say that the solution $u(x,t)$ satisfies the Kruzhkov entropy if
	\begin{align*}
\int_{0}^T\int_{{\mathbb{S}^1}}\Big(\left|u(x,t)-A\right|g_t(x,t)+\sgn(u(x,t)-A)[&u(x,t)f(u(x,t))-Af(A)]g_x(x,t)\Big)dxdt+\notag \\
&+ \int_{{\mathbb{S}^1}}\left|u_0(x)-A\right|g(x,0)dx\geq 0.
\end{align*}
for all $g(x,t)\in\mathcal{C}^\infty_0(\mathbb{S}^1\times[0,T))$.
	For this proof, we assume that the sequence generated by scheme $(\ref{ODE})$ is pre-compact, which is demonstrated in Appendix $\ref{precompact}$.

\begin{remark}
\label{rem}
Since $f$ is a Lipschitz function, then, for a sufficiently large constant $k$ in $f^+=\max(f,0)+k$ and $f^-=\max(-f,0)+k$,  we have that
\begin{equation}
uf^+(\hat{u})\leq af^+(a) \quad \text{ and } 
\quad uf^-(\hat{u})\leq af^-(a) \quad
\text{ if } u\leq a.\label{dift}
\end{equation}
\end{remark}

\begin{proposition}
Let us assume that constant $k$ is sufficiently large so that Eq. $(\ref{dift})$ is satisfied on segment $[-M,M]$. Then $u(x,t,\epsilon)\longrightarrow u(x,t)$ when $\epsilon\longrightarrow 0$ in $L_{loc}^1({\mathbb{S}^1}\times[0,\infty))$, when $u(x,t)$ is the only entropy solution to $(\ref{noLi})$.  
\end{proposition}
\textbf{Proof.}
 We consider a constant $A\in[-M,M]$. For almost $(x,t)\in {\mathbb{S}^1}\times(0,\infty)$ and fixed $x$, we differentiate $|u(x,t)-A|$ and then using $(\ref{ODE})$, we obtain
\begin{align}
&\frac{d}{dt}|u(x,t,\epsilon)-A|=
\sgn(u_{_{\epsilon}}-A)\frac{d}{dt}u(x,t,\epsilon)\notag\\
& =\frac{1}{\epsilon}\sgn(u_{_{\epsilon}}-A)\left[u_{_{\epsilon-1}}f^+\left(\hat{u}_{_{\epsilon-1/2}}\right)-
u_{_{\epsilon}}f^+\left(\hat{u}{_{\epsilon+1/2}}\right)-
u_{_{\epsilon}}f^-\left(\hat{u}_{_{\epsilon-1/2}}\right)+
u_{_{\epsilon+1}}f^-\left(\hat{u}{_{\epsilon+1/2}}\right)
\right]\notag\\
& =\frac{1}{\epsilon}\sgn(u_{_{\epsilon}}-A)\left[(u_{_{\epsilon-1}}f^+\left(\hat{u}_{_{\epsilon-1/2}}\right)-Af^+(A))+(u_{_{\epsilon+1}}f^-\left(\hat{u}{_{\epsilon+1/2}}\right)-Af^-(A))\right]\notag\\
&-\frac{1}{\epsilon}\sgn(u_{_{\epsilon}}-A)\left[(u_{_{\epsilon}}\left(f^+\left(\hat{u}_{_{\epsilon+1/2}}\right)+f^-\left(\hat{u}_{_{\epsilon-1/2}}\right)\right)-
A(f^+(A)+
f^-(A))\right]\label{fcontast1}.
\end{align}
Since Remark \ref{rem} and $(\ref{dift})$ are valid, for $u$ and $A$ in $[-M,M]$, we find that
\begin{align}
&\sgn(u_{_{\epsilon}}-A)\left[(u_{_{\epsilon}}\left(f^+\left(\hat{u}_{_{\epsilon+1/2}}\right)+
f^-\left(\hat{u}_{_{\epsilon-1/2}}\right)\right)-
A(f^+(A)+
f^-(A))\right]=\notag\\
&\left|u_{_{\epsilon}}f^+\left(\hat{u}_{_{\epsilon+1/2}}\right)-Af^+(A)\right|+\left|u_{_{\epsilon}}f^-\left(\hat{u}_{_{\epsilon-1/2}}\right)-Af^-(A)\right|\label{fcontast}.
\end{align}
 Substituting $(\ref{fcontast})$ in $(\ref{fcontast1})$, we can estimate 
	\begin{align}
\frac{d}{dt}&|u(x,t,\epsilon)-A|\leq \frac{1}{\epsilon} \Big\{\left|u_{_{\epsilon-1}}f^+\left(\hat{u}_{_{\epsilon-1/2}}\right)-Af^+(A)\right|-\left|u_{_{\epsilon}}f^+\left(\hat{u}_{_{\epsilon+1/2}}\right)-Af^+(A)\right| +\notag\\
&+\left|u_{_{\epsilon+1}}f^-\left(\hat{u}_{_{\epsilon+1/2}}\right)-Af^-(A)\right|-\left|u_{_{\epsilon}}f^-\left(\hat{u}_{_{\epsilon-1/2}}\right)-Af^-(A)\right|\Big\}.\label{ineq1}
\end{align}
	
Multiplying inequality $(\ref{ineq1})$ by the nonnegative test function $g=g(x,t)\in C^\infty_0({\mathbb{S}^1}\times[0,T)),T>0$ and integrating by parts, we obtain
\begin{align}
-&\int_{{\mathbb{S}^1}}\left|u_0(x)-A\right|g(x,0)dx-\int_{0}^T\int_{{\mathbb{S}^1}}\left|u(x,t,\epsilon)-A\right|g_t(x,t)dx dt\leq \notag \\
&\int_0^T\int_{{\mathbb{S}^1}}\frac{1}{\epsilon} \Big\{\left|u_{_{\epsilon-1}}f^+\left(\hat{u}_{_{\epsilon-1/2}}\right)-Af^+(A)\right|-\left|u_{_{\epsilon}}f^+\left(\hat{u}_{_{\epsilon+1}}\right)-Af^+(A)\right| +\notag\\
&+\left|u_{_{\epsilon+1}}f^-\left(\hat{u}_{_{\epsilon+1/2}}\right)-Af^-(A)\right|-\left|u_{_{\epsilon}}f^-\left(\hat{u}_{_{\epsilon-1/2}}\right)-Af^-(A)\right|\Big\}g(x,t)dxdt.\label{notw}
\end{align}
{{
Note that if we perform the change of variable
$x=x+\epsilon$, we can rewrite
\begin{align}
\int_0^T\int_{{\mathbb{S}^1}}\Big\{\left|u_{_{\epsilon-1}}f^+\left(\hat{u}_{_{\epsilon-1/2}}\right)-Af^+(A)\right|\Big\}g(x,t)dxdt=\int_0^T\int_{{\hat{\mathbb{S}}^1}}\Big\{\left|u_{_{\epsilon}}f^+\left(\hat{u}_{_{\epsilon+1/2}}\right)-Af^+(A)\right|\Big\}g(x+\epsilon,t)dxdt\notag
\end{align}
where $\hat{\mathbb{S}}^1$ represents a shift of $\epsilon$ in $\mathbb{S}^1$.
Since $g(x,t)$ has compact support, we take the support in such way that 
\begin{align}
\int_0^T\int_{{\mathbb{S}^1}}\Big\{\left|u_{_{\epsilon-1}}f^+\left(\hat{u}_{_{\epsilon-1/2}}\right)-Af^+(A)\right|\Big\}g(x,t)dxdt&=\int_0^T\int_{{\hat{\mathbb{S}}^1}}\Big\{\left|u_{_{\epsilon}}f^+\left(\hat{u}_{_{\epsilon+1/2}}\right)-Af^+(A)\right|\Big\}g(x+\epsilon,t)dxdt\notag\\&=\int_0^T\int_{{{\mathbb{S}}^1}}\Big\{\left|u_{_{\epsilon}}f^+\left(\hat{u}_{_{\epsilon+1/2}}\right)-Af^+(A)\right|\Big\}g(x+\epsilon,t)dxdt \label{valt}
\end{align}
Performing the change of variables $x=x-\epsilon$ and using the similar argument used in $(\ref{valt})$, we prove that
\begin{align}
\int_0^T\int_{{\mathbb{S}^1}}\Big\{\left|u_{_{\epsilon+1}}f^-\left(\hat{u}_{_{\epsilon+1/2}}\right)-Af^-(A)\right|\Big\}g(x,t)dxdt=\int_0^T\int_{{{\mathbb{S}}^1}}\Big\{\left|u_{_{\epsilon}}f^-\left(\hat{u}_{_{\epsilon-1/2}}\right)-Af^-(A)\right|\Big\}g(x-\epsilon,t)dxdt \label{valt2}
\end{align}
Applying $(\ref{valt})$ and $(\ref{valt2})$ in the inequality $(\ref{notw})$, we obtain
}}
\begin{align}
-&\int_{{\mathbb{S}^1}}\left|u_0(x)-A\right|g(x,0)dx-\int_{0}^T\int_{{\mathbb{S}^1}}\left|u(x,t,\epsilon)-A\right|g_t(x,t)dx dt\leq \notag \\
&\int_0^T\int_{{\mathbb{S}^1}}\Big\{\left|u_{_{\epsilon}}f^+\left(\hat{u}_{_{\epsilon+1/2}}\right)-Af^+(A)\right|\frac{g(x+\epsilon,t)-g(x,t)}{\epsilon}+\notag\\
&+\left|u_{_{\epsilon}}f^-\left(\hat{u}_{_{\epsilon-1/2}}\right)-Af^-(A)\right|\frac{g(x,t)-g(x-\epsilon,t)}{\epsilon}dx dt= \notag \\
&= \int_0^T\int_{{\mathbb{S}^1}} \left(\left|u_{_{\epsilon}}f^+\left(\hat{u}_{_{\epsilon+1/2}}\right)-Af^+(A)\right|-\left|u_{_{\epsilon}}f^-\left(\hat{u}_{_{\epsilon-1/2}}\right)-Af^-(A)\right|\right)g_x dxdt+I(\epsilon).\label{notw2}
\end{align}
Since $g\in C^\infty_0({\mathbb{S}^1}\times[0,T))$, then $I(\epsilon)\longrightarrow 0$ when $\epsilon \longrightarrow 0$. Moreover, since Eq. $(\ref{dift})$ is satisfied, we have that
\begin{align}
&\left|u_{_{\epsilon}}f^+\left(\hat{u}_{_{\epsilon+1/2}}\right)-Af^+(A)\right|-\left|u_{_{\epsilon}}f^-\left(\hat{u}_{_{\epsilon-1/2}}\right)-Af^-(A)\right|=\notag\\
&\sgn(u_{_{\epsilon}}-A)\left[u_{_{\epsilon}}f^+\left(\hat{u}_{_{\epsilon+1/2}}\right)-Af^+(A)-(u_{_{\epsilon}}f^-\left(\hat{u}_{_{\epsilon-1/2}}\right)-Af^-(A))\right]=\notag \\
&\sgn(u_{_{\epsilon}}-A)\left[u_{_{\epsilon}}\left(f^+\left(\hat{u}_{_{\epsilon+1/2}}\right)-f^-\left(\hat{u}_{_{\epsilon-1/2}}\right)\right)-(A(f^+(A)-f^-(A))\right]=\notag\\
&\sgn(u_{_{\epsilon}}-A)\left[u_{_{\epsilon}}\left(f^+\left(\hat{u}_{_{\epsilon+1/2}}\right)-f^-\left(\hat{u}_{_{\epsilon-1/2}}\right)\right)-Af(A)\right].\label{intd}
\end{align}
In Eq. $(\ref{intd})$, we use $f^+-f^-=f$. By substituting the result of $(\ref{intd})$ into Eq. $(\ref{notw2})$
\begin{align}
&\int_{0}^T\int_{{\mathbb{S}^1}}\Big(\left|u(x,t,\epsilon)-A\right|g_t(x,t)+\sgn(u_{_{\epsilon}}-A)[u_{_{\epsilon}}\left(f^+\left(\hat{u}_{_{\epsilon+1/2}}\right)-f^-\left(\hat{u}_{_{\epsilon-1/2}}\right)\right)-Af(A))]\Big)dxdt+\notag \\
&+ \int_{{\mathbb{S}^1}}\left|u_0(x)-A\right|g(x,0)dx\geq -I(\epsilon).\label{utte}
\end{align}
Here, $u_{_{\epsilon}}=u(x,t,\epsilon)$. In Appendix \ref{precompact}, we show that family $u(x,t,\epsilon)$ for $\epsilon>0$ is a pre-compact sequence in $L^1({\mathbb{S}^1}\times[0,T])$. Let $u(x,t)$ be an accumulation point of family $u(x,t,\epsilon)$, thus for a subsequence $\epsilon_r$, we have that $u(x,t,\epsilon_r)\longrightarrow u(x,t)$ when $r\longrightarrow \infty$ in $L^1({\mathbb{S}^1}\times [0,T])$ and from Eq. $(\ref{fds})$, we have that  $u_{\epsilon_r+1/2}={\hat{u}\left(x+\frac{1}{2},t,\epsilon_r\right)\longrightarrow u(x,t)}$ and $u_{\epsilon_r-1/2}={\hat{u}\left(x-\frac{1}{2},t,\epsilon_r\right)\longrightarrow u(x,t)}$. Taking  $\epsilon=\epsilon_r\longrightarrow 0$ in $(\ref{utte})$ and using that $f^+-f^-=f$, we obtain the entropy relation, remembering that $I(\epsilon)\longrightarrow 0$
\begin{align}
\int_{0}^T\int_{{\mathbb{S}^1}}&\Big(\left|u(x,t)-A\right|g_t(x,t)+\sgn(u(x,t)-A)[u(x,t)f(u(x,t)-Af(A))]g_x(x,t)\Big)dxdt+\notag \\
&+ \int_{{\mathbb{S}^1}}\left|u_0(x)-A\right|g(x,0)dx\geq 0.\label{utte2}
\end{align}
In Eq. $(\ref{utte2})$, $ A \in [-M,M]$. However, for $|A|\geq M$, notice that the inequality that is Eq. $(\ref{utte2})$ reduces to the equality (weak solution)  
\begin{align*}
\int_{0}^T\int_{{\mathbb{S}^1}}\Big(u(x,t)g_t(x,t)+u(x,t)f(u(x,t)g_x(x,t)\Big)dxdt+ \int_{{\mathbb{S}^1}}u_0(x)g(x,0)dx= 0.
\end{align*}
From these results, we obtain that $(\ref{utte2})$ holds for all $A\in \mathbb{R}$. Since $T>0$ and $g=g(x,t)\in C^\infty_0({\mathbb{S}^1}\times[0,T))$ are arbitrary, inequality $(\ref{utte2})$ leads to solution $u(x,t)$, which is the entropy solution to $(\ref{noLi})$. This solution is unique; in particular, an accumulation point $u(x,t)$ of $u(x,t,\epsilon)$ using $(\ref{ODE})$ is what is unique about it. This implies that family $u(x,t,\epsilon)$ converges to $u(x,t)$ as $\epsilon\longrightarrow 0$ in $L^1_{loc}({\mathbb{S}^1}\times[0,\infty))$ because $T$ is arbitrary, which completes the proof. $\quad \square$.

	\section{Numerical experiments}\label{NumExp}
	
	In order to illustrate the robustness of the proposed 
	numerical scheme, we present numerical experiments 
	describing the explicit calculation of the weak 
	asymptotic approximations for concrete conservation 
	law equations. We also provide examples for systems of equations.
        All the calculations were performed in 
	the order of seconds with MATLAB on a standard 
	desktop computer.

	\subsection{Comparison between numerical studies and the W1 distance}
	
\begin{figure}[p]
		\centering
		\includegraphics[width=0.3\textwidth]{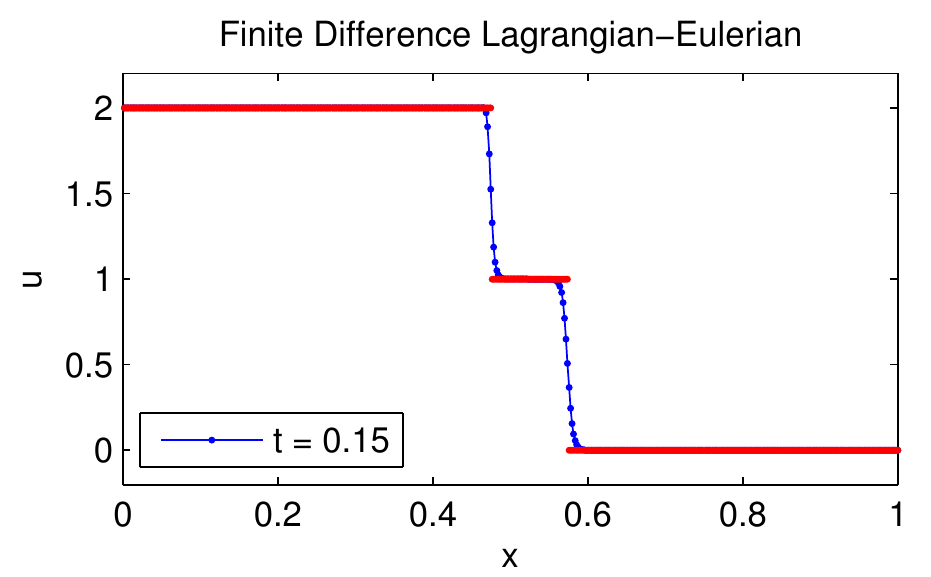}
		\includegraphics[width=0.3\textwidth]{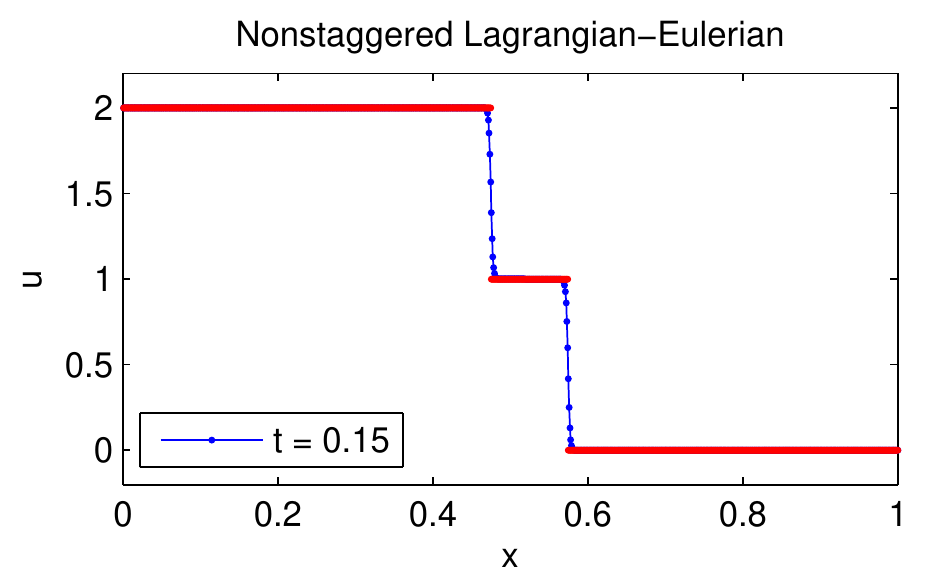}
		\includegraphics[width=0.3\textwidth]{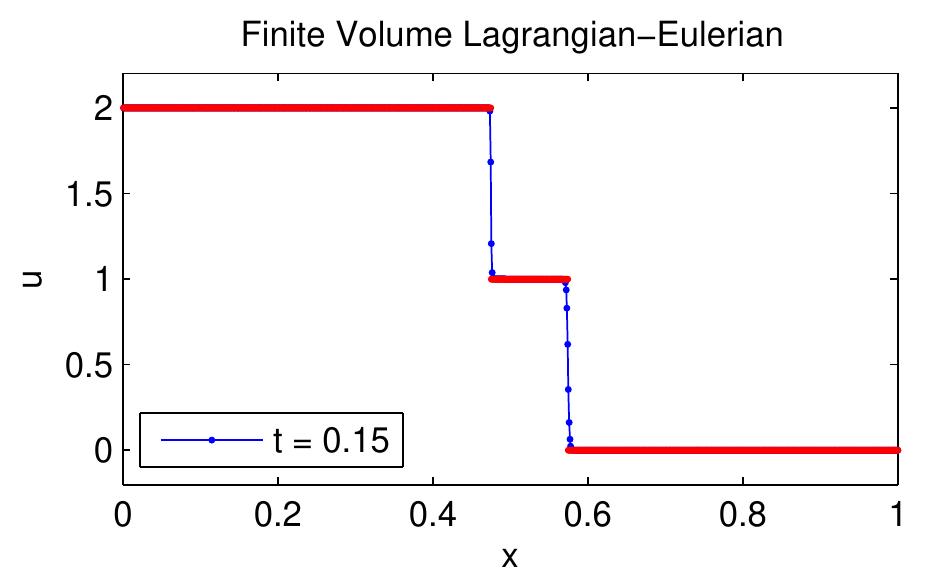}\\
		\includegraphics[width=0.3\textwidth]{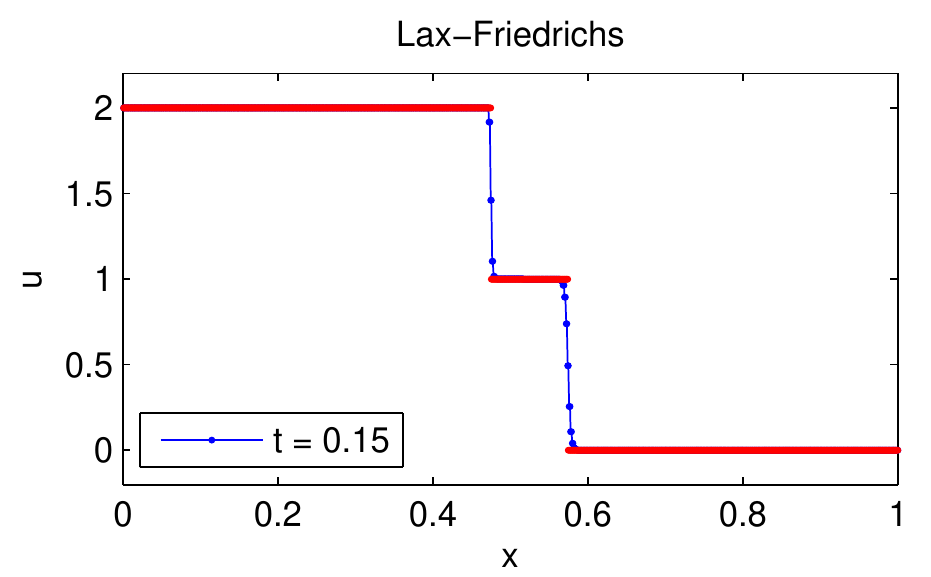}
		\includegraphics[width=0.3\textwidth]{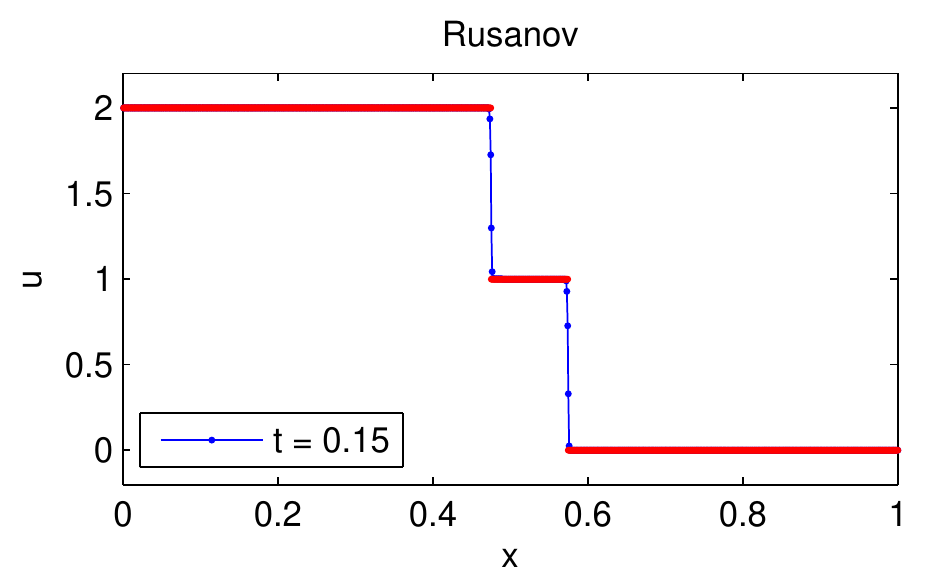}
		\includegraphics[width=0.3\textwidth]{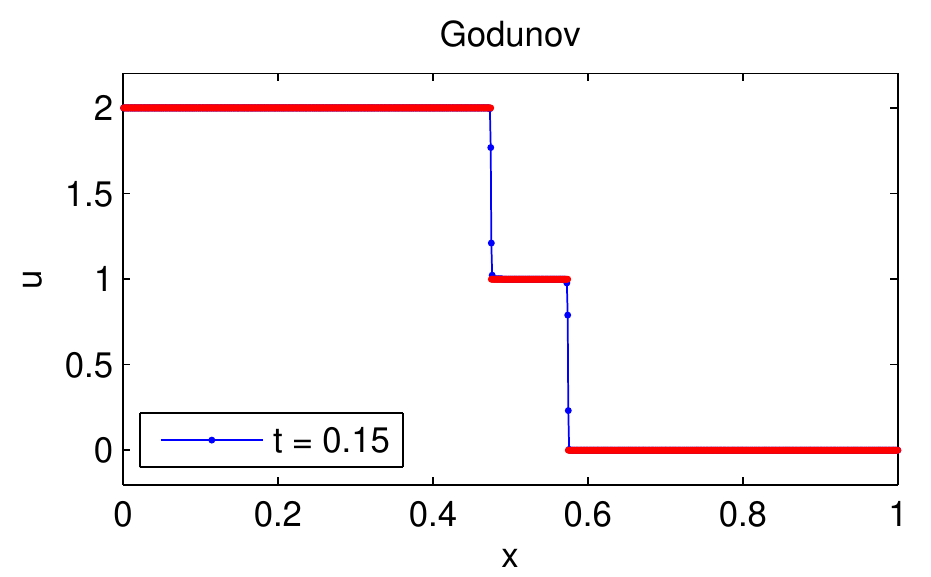}
		\caption{Numerical solutions to problem P1 at time $t=0.15$ 
			(before shock).}
		\label{figL1W1Burgers1}
	\end{figure}

\begin{figure}[p]
		\centering
		\includegraphics[width=0.3\textwidth]{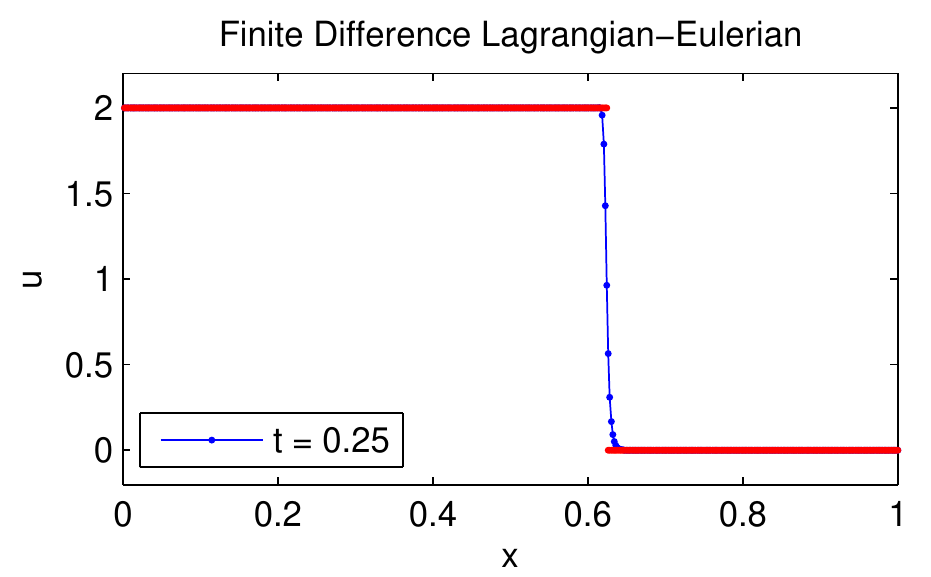}
		\includegraphics[width=0.3\textwidth]{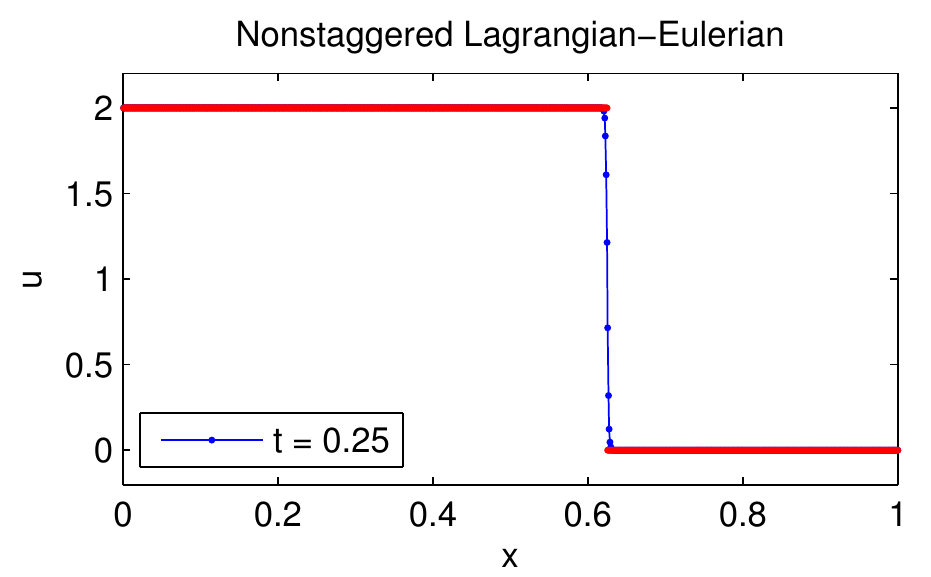}
		\includegraphics[width=0.3\textwidth]{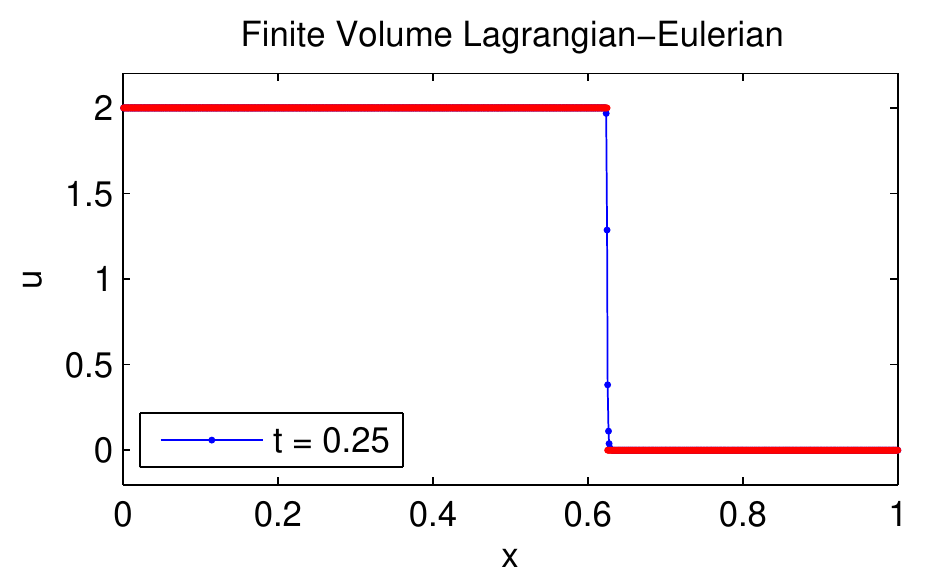}\\
		\includegraphics[width=0.3\textwidth]{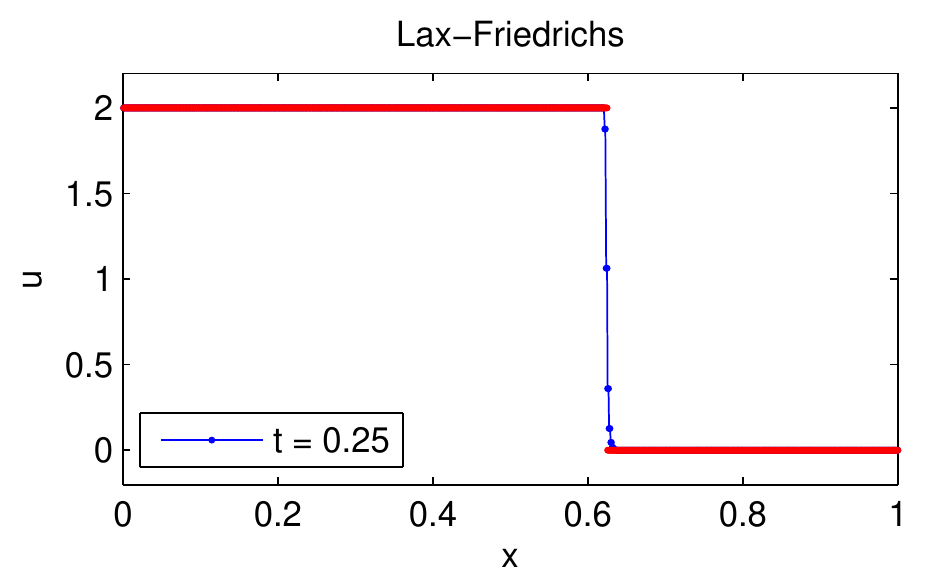}
		\includegraphics[width=0.3\textwidth]{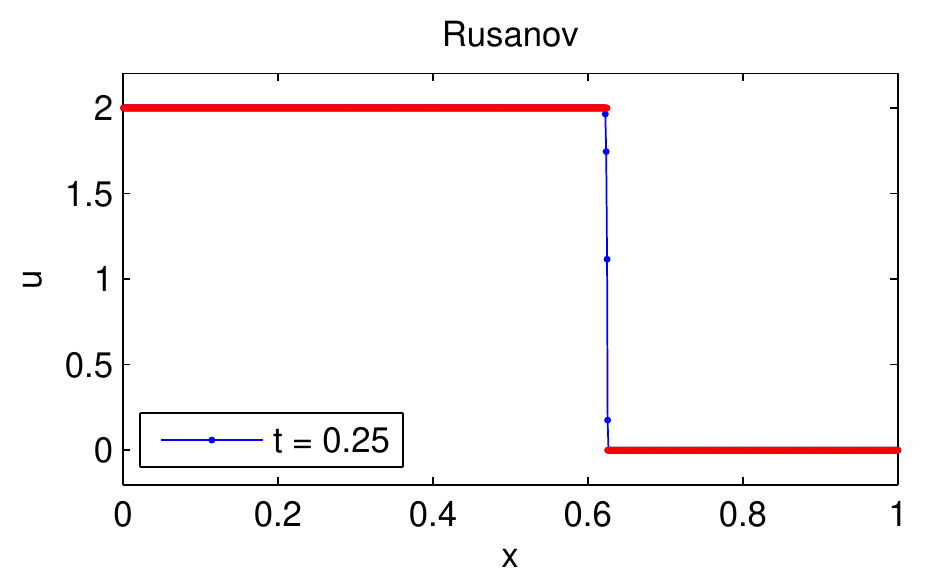}
		\includegraphics[width=0.3\textwidth]{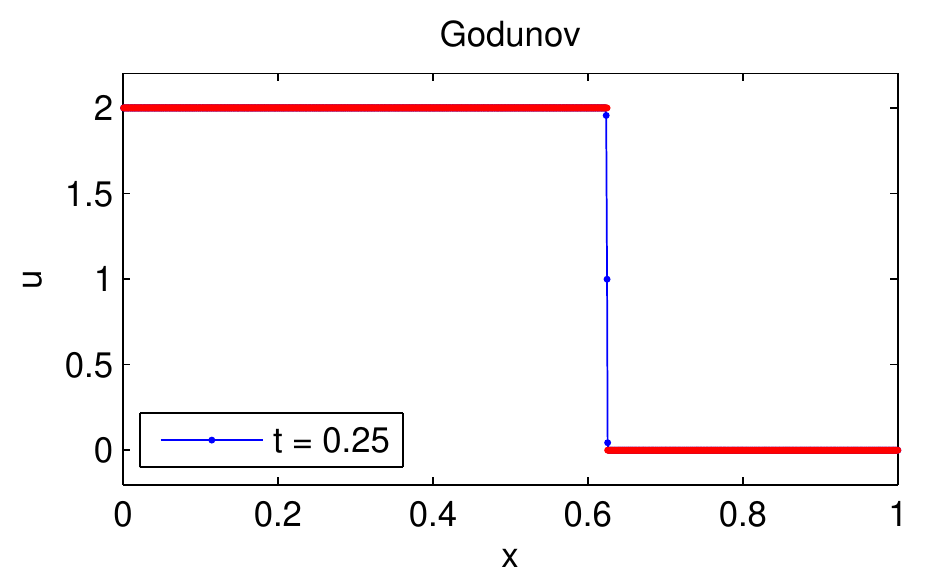}
		\caption{Numerical solutions to problem P1 at time $t=0.25$ 
			(after shock).}
		\label{figL1W1Burgers2}
	\end{figure}

\begin{figure}[p]
		\centering
		\includegraphics[width=0.3\textwidth]{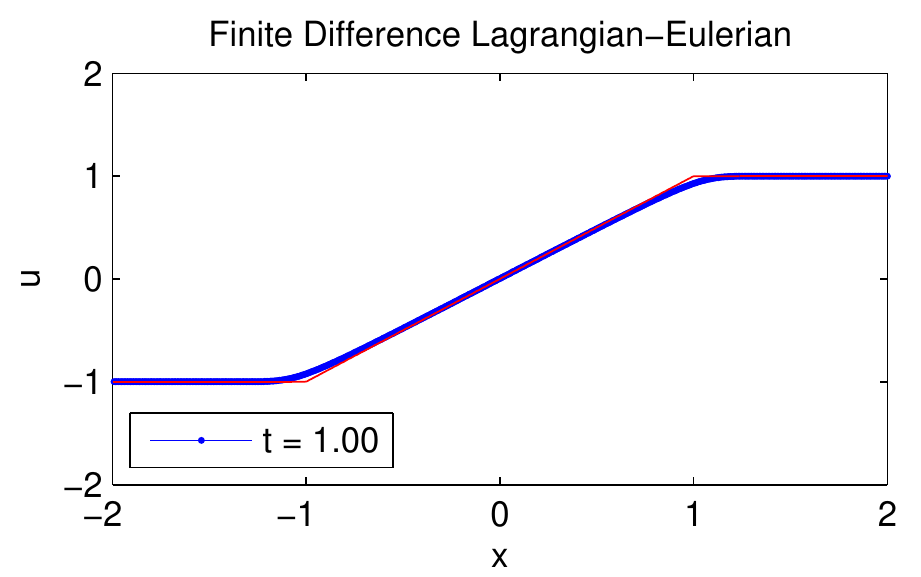}
		\includegraphics[width=0.3\textwidth]{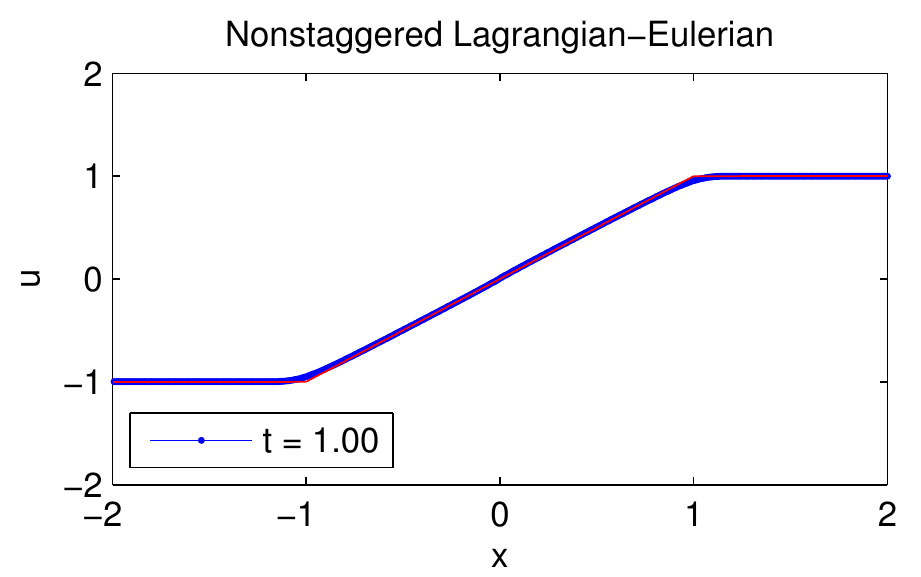}
		\includegraphics[width=0.3\textwidth]{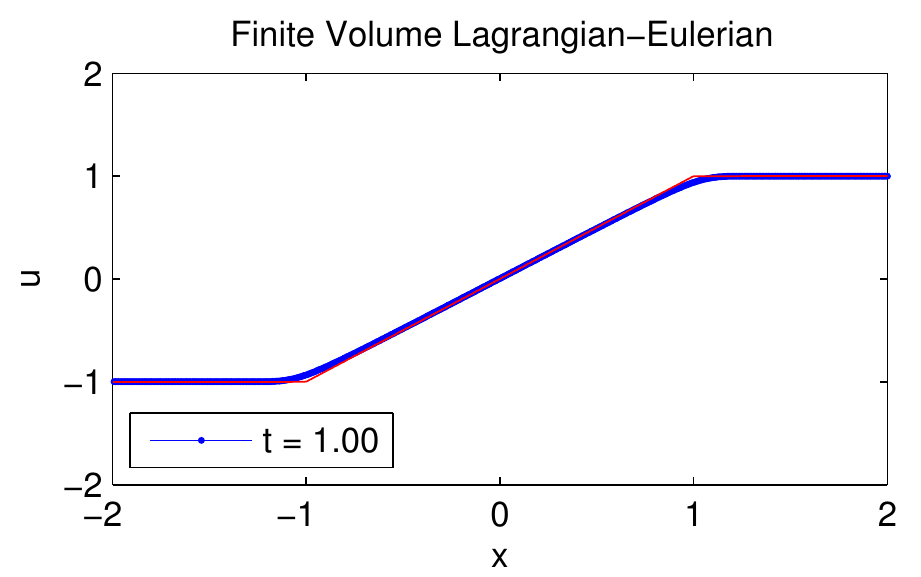}\\
		\includegraphics[width=0.3\textwidth]{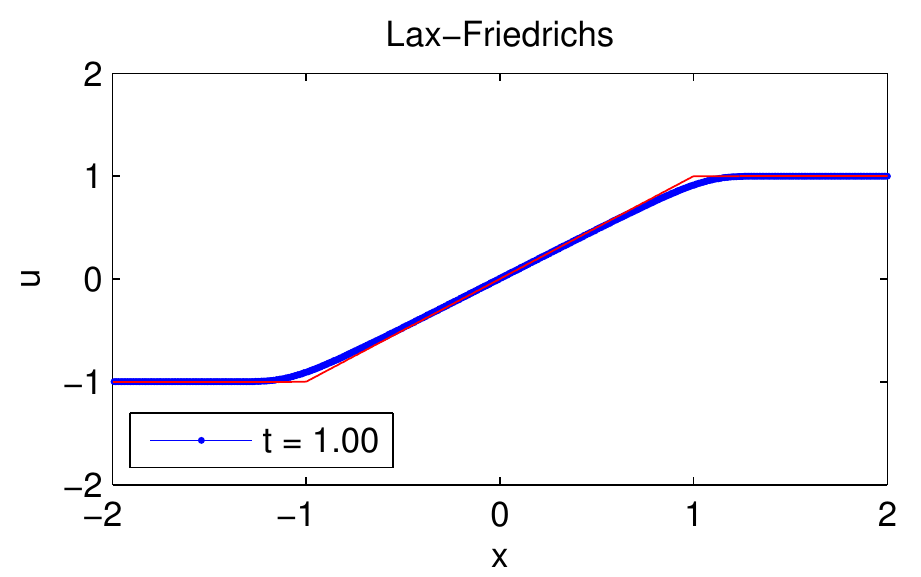}
		\includegraphics[width=0.3\textwidth]{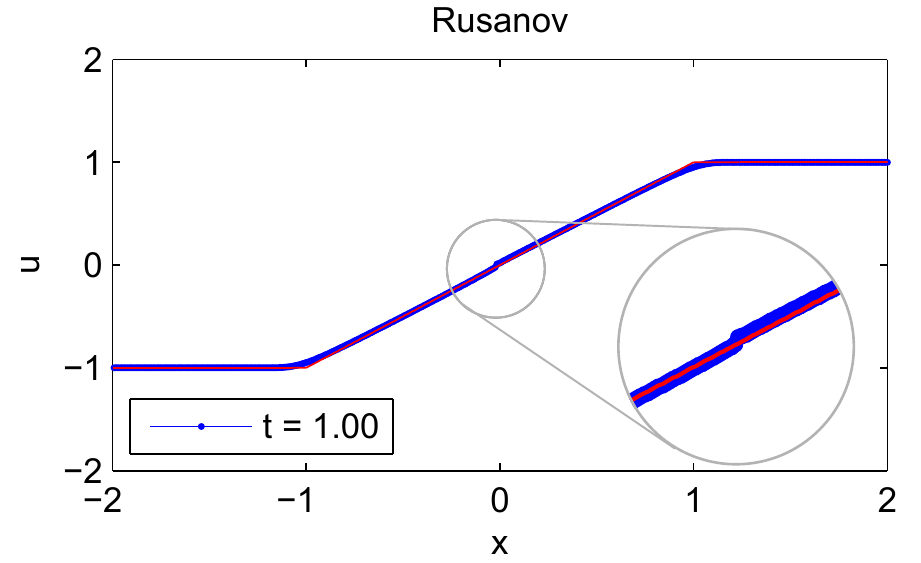}
		\includegraphics[width=0.3\textwidth]{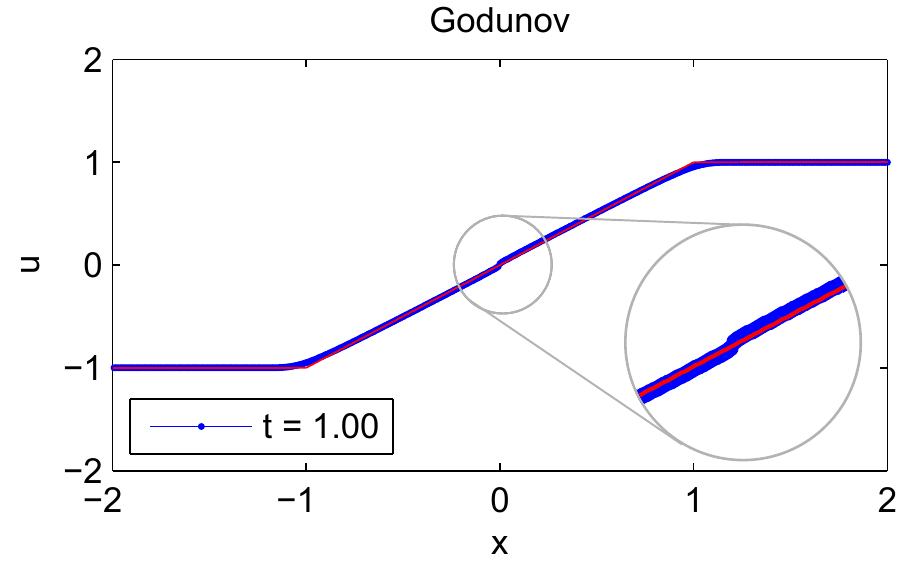}
		\caption{Numerical solutions to problem P2 at time $t=0.5$ 
			(rarefaction). Notice that the Lagrangian-Eulerian scheme does not produce the well-known spurious glitch effect in the sonic rarefaction present in Godunov and Rusanov's simulations.}
		\label{figL1W1Burgers3}
	\end{figure}

\begin{figure}[p]
		\centering
        \includegraphics[width=0.3\textwidth]{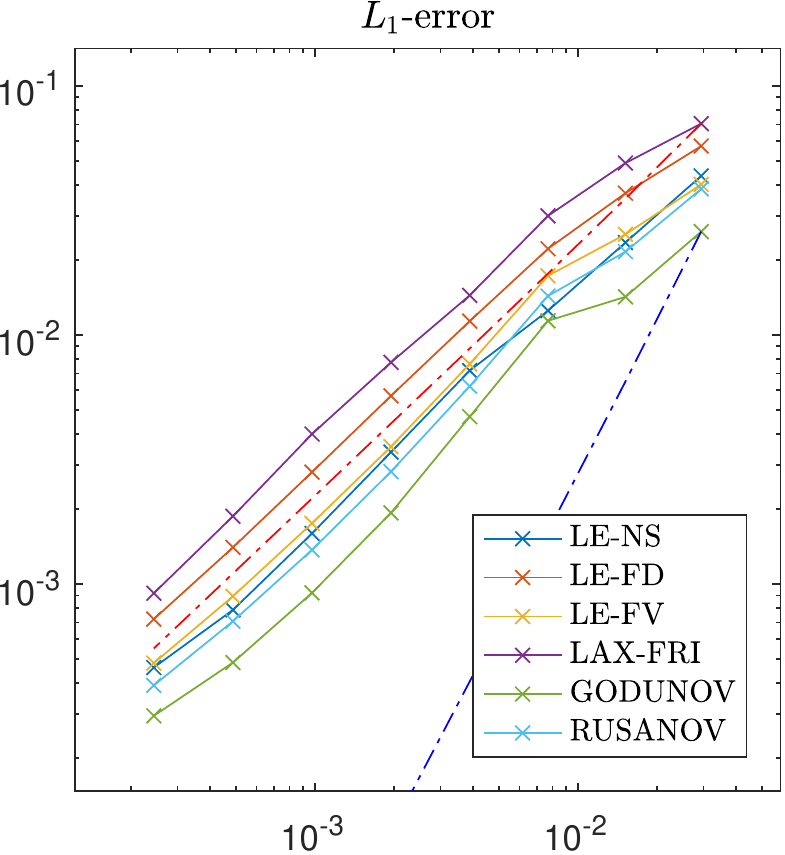}
        \includegraphics[width=0.3\textwidth]{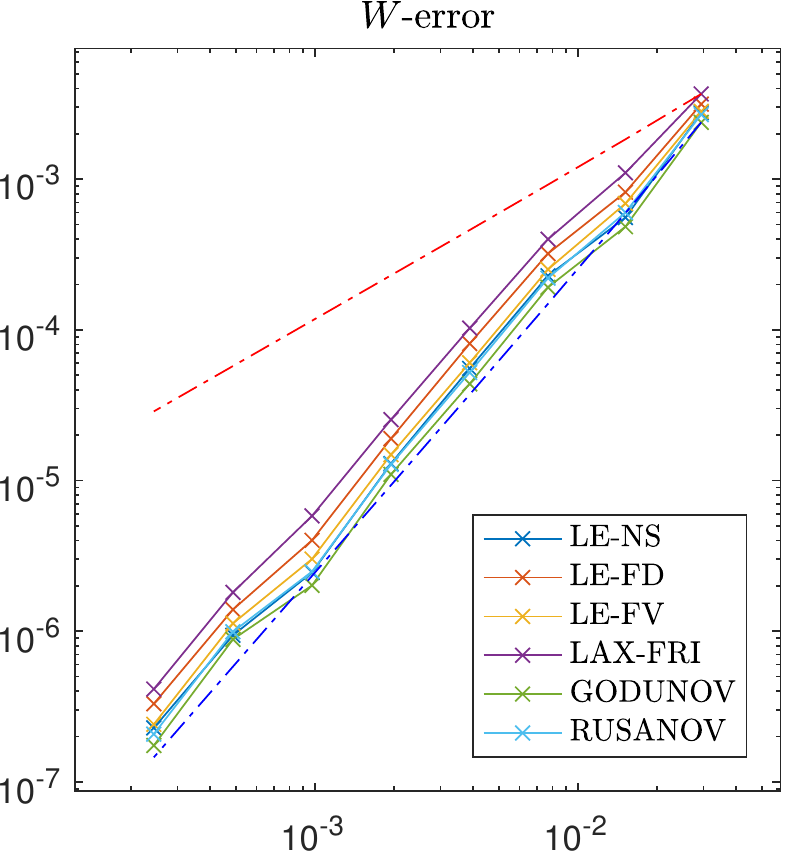}\\
		\includegraphics[width=0.3\textwidth]{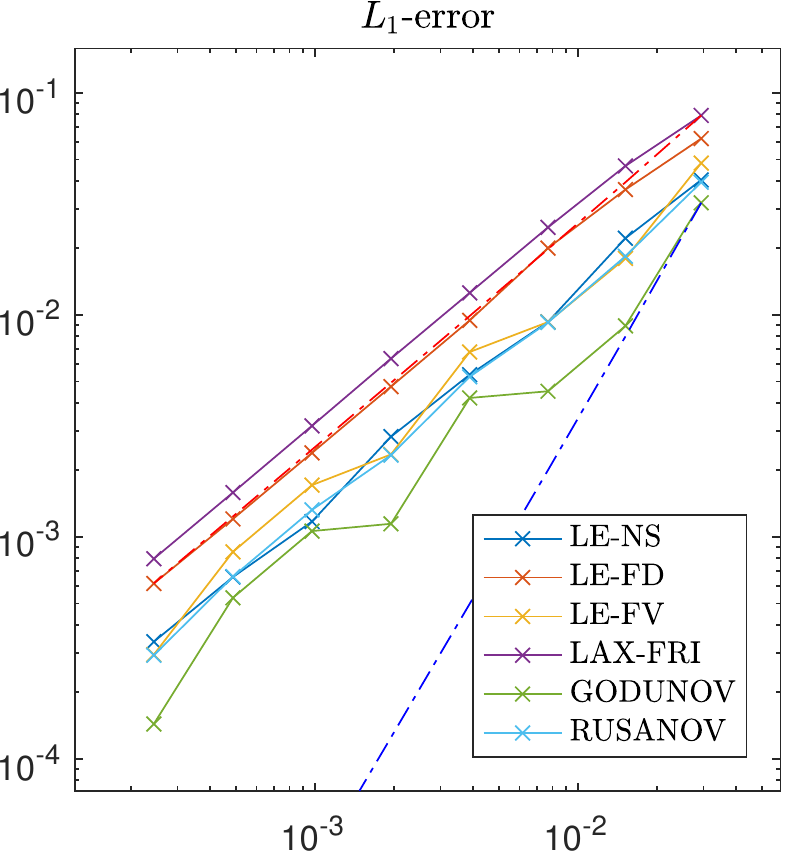}
		\includegraphics[width=0.3\textwidth]{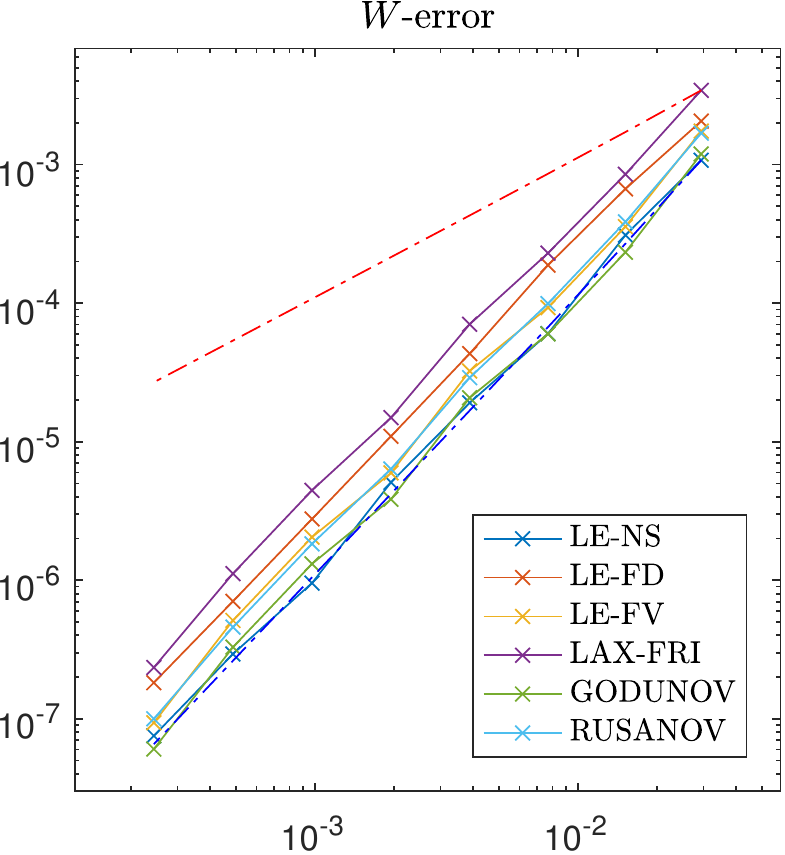}
		\caption{Log-log plots for norms $L^1$ 
			 and $W$ of the 
			error versus the cell sizes for problem P1 at time $t=0.15$ 
			(before shock; first two pictures) and at time $t=0.25$ (after shock; last two pictures). The top solid line  represents the convergence of the Lax--Friedrichs numerical scheme, while the bottom solid line marks the convergence of the Godunov method. The error obtained with the Nonstaggered Lagrangian--Eulerian scheme in these cases approaches that of Godunov and is sometimes lower than that of Rusanov.}
		\label{loglogerrorB1B2}
	\end{figure}

	The Wasserstein 
	distance between two
	probability measures $\mu$ and $\nu$ on $\mathbb{R}$ 
	can equivalently be defined as
	\begin{equation}
	W_1(\mu, \nu) := \sup_{||\varphi||_{\text{Lip}} \leq 1} \int_{\mathbb{R}}\varphi(x)d (\mu-\nu)(x).
	\end{equation}
	Here, the supremum is taken over all functions 
	$\varphi:\mathbb{R}\to\mathbb{R}$ with Lipschitz 
	semi-norm $	
	||\varphi||_{\text{Lip}} := 
	\sup_{x\neq y}|\frac{\varphi(x)-\varphi(y)}{x-y}|,
$ 
	at most 1. Given Borel measurable functions $u, 
	v : \mathbb{R}\to\mathbb{R}$ satisfying the analogous properties,
$	\int_{\mathbb{R}}(u-v)(x)dx = 0, \quad \int_{\mathbb{R}}|x||u-v|(x)dx < \infty.
$	
Following \cite{USFSS16}, given an exact and an approximate solution to (\ref{noLi}), the difference between them has zero mass when the numerical scheme is conservative, and decays sufficiently fast. The Wasserstein error $W_1$ {(i.e., computing the error in the Lip'-norm)} must be well-defined and finite by measuring the amount of work that goes into moving the surplus of mass to behind the shock, where there is a shortage of mass. 
{ 
In addition, we implemented the nonstaggered Lagrangian--Eulerian scheme (\ref{eq32})--(\ref{projection}) presented in Section \ref{sec:lefvm} and we reproduced numerical experiments introduced in \cite{USFSS16} for Burgers’ equation on interval $[0,1]$ (see Figures \ref{figL1W1Burgers1} to \ref{loglogerrorB1B2}). \\
}  
{
For {\bf Model Problem P1}:
}
$	u_t + \left(\frac{u^2}{2}\right)_x = 0,
$	with initial data containing two jumps 
$	u_0(x) =  \left\{\begin{array}{ll}
	 2, \quad x \leq \frac14,\\
	 1, \quad \frac14 \leq x < \frac12,\\
	 0, \quad \frac12 \leq x,
	\end{array}\right.
$	
{we found that our scheme to the underlying set up P1
is $O(\Delta x)$ in $L_1$ and $O(\Delta x^2)$ in $W_1$ 
(see Figure \ref{loglogerrorB1B2}) in the presence of shocks;
see Figure \ref{figL1W1Burgers1} and Figure \ref{figL1W1Burgers2}.}
The exact solution to P1 is, for $t < 0.25$, is
$	u(x,t) =  \left\{\begin{array}{ll}
	2, \quad x \leq \frac{1+3t}{4},\\
	1, \quad \frac{1+3t}{4} \leq x < \frac{1+t}{2},\\
	0, \quad \frac{1+t}{2} \leq x,
	\end{array}\right.
$	and, for $t > 0.25$,
$	u(x,t) =  \left\{\begin{matrix}
	2, & x \leq \frac38+t,\\
	0, & x \geq \frac38+t.
	\end{matrix}\right.$ 
We observe that the new reconstruction step to the nonstaggered Lagrangian--Eulerian scheme (\ref{eq32})--(\ref{projection}) does not tends to smooth the reconstruction variations by using slope limiters without introducing excessive numerical diffusion or spurious oscillations in the interaction of the discontinuities in the solution as times evolves.

	We also considered Burger's equation with 
	initial data ({\bf Model Problem P2}), given by
$	u(x,t) = \left\{ \begin{matrix}
	-1, & x \leq 0,\\
	1, & 0 \leq x,
	\end{matrix}\right.
$	whose exact solution is a rarefaction wave, namely,
$	u(x,t) = \left\{ \begin{matrix}
	-1, & x \leq -t,\\
	x, & -t \leq x < t,\\
	1, & t \leq x,
	\end{matrix}\right.$ { 
For this set up P2, we have a rarefaction solution of the inviscid Burgers' equation as times evolves. In particular, the proposed scheme is able to capture with good resolution the rarefaction wave in the vicinity where a sign change in the wave speeds is observed at point $u=0$. On the other hand, we might see from the Figure \ref{figL1W1Burgers3} that both the classical Godunov and the Rusanov schemes produce the spurious {\it sonic glitch or entropy glitch} effect, located at the point  $u=0$. Such phenomenon  arises in the presence of sonic rarefaction waves due to the change in signal of the wave speeds. For this set up P2, we also observed that our method is $O(\Delta x)$ in $L_1$-norm, and $O(\Delta x^2)$ in $W_1$-norm.}


	\subsection{A nonlocal traffic model}
	We present numerical approximations 
	of the classical Lighthill--Whitham--Richards (LWR) 
	model for vehicular traffic \cite{amorim2015}, which 
	consists of a continuity equation
	\begin{equation}
	\partial_{t} \rho + \partial_{x} \left( \rho V \right) = 0, \qquad \rho \in [0,1],
	\end{equation}
	where $\rho$ is the
	(average) vehicular density. Density $\rho$ is a 
	function of $t$ (time), and $x$ is a position 
	along a road with neither entries nor exits. In 
	this equation, a \textit{speed law} function 
	$V = V(\rho)$ is defined as follows:
	\begin{equation}
	V(\rho) = V_{\max} (1-\rho)(1-\rho \ast \eta).
	\label{eq:3.2}
	\end{equation}
	
	By setting $V_{\max} > 0$, this flux function 
	can be used as an LWR-type macroscopic model
	for vehicular traffic, where drivers adjust 
	their speed according to the local traffic density. 
	The convolution is realized with $\eta$ ($\alpha$ 
	is chosen so that $\int_{\mathbb{R}} \eta = 1$), 
	which is defined as	
	\begin{equation}
	\eta(x) = \begin{cases}
	\alpha \left((x_1-x)(x-x_2)\right)^{5/2}, \quad -x_1 \leq x \leq x_2 \\
	0, \quad \text{ otherwise.}
	\end{cases} 
	\label{eq:3.1}
	\end{equation}
	
\begin{figure}[p]
		\centering
		\includegraphics[width=0.24\textwidth]{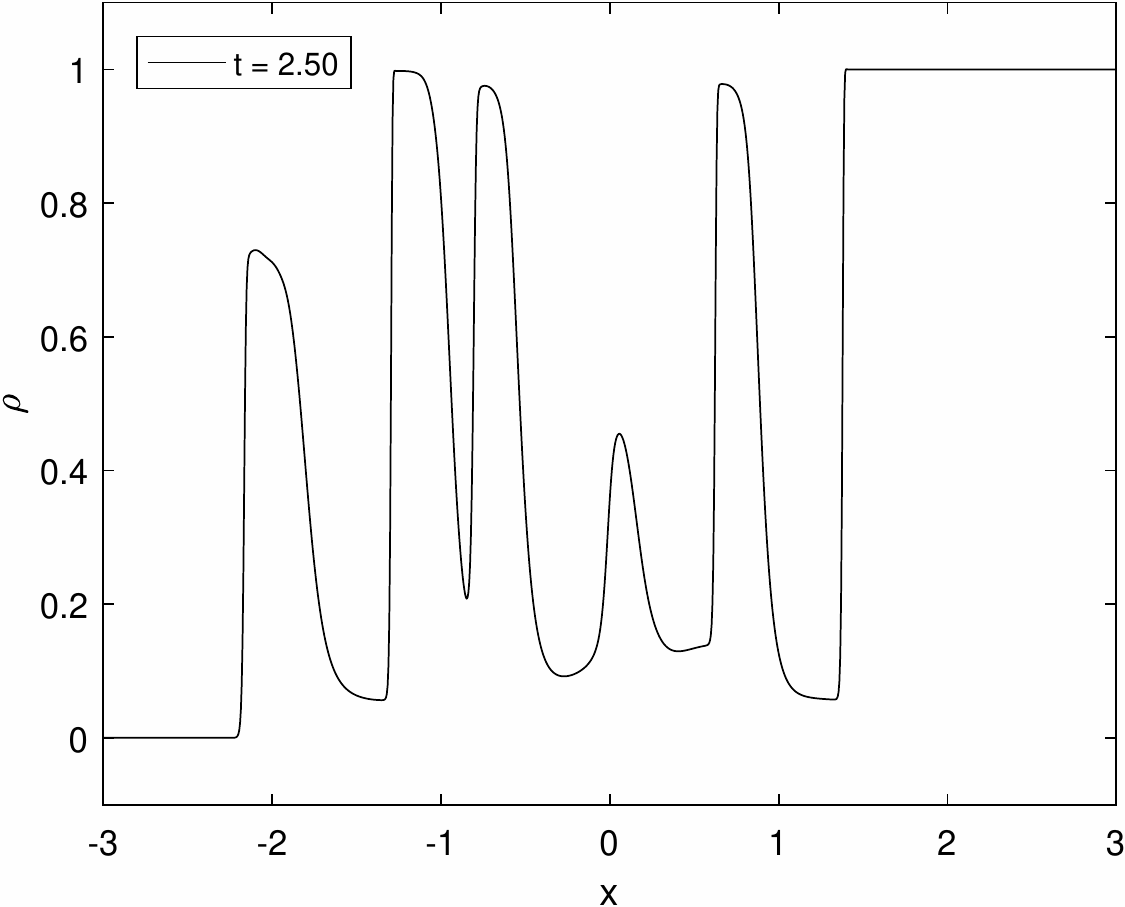}
		\includegraphics[width=0.24\textwidth]{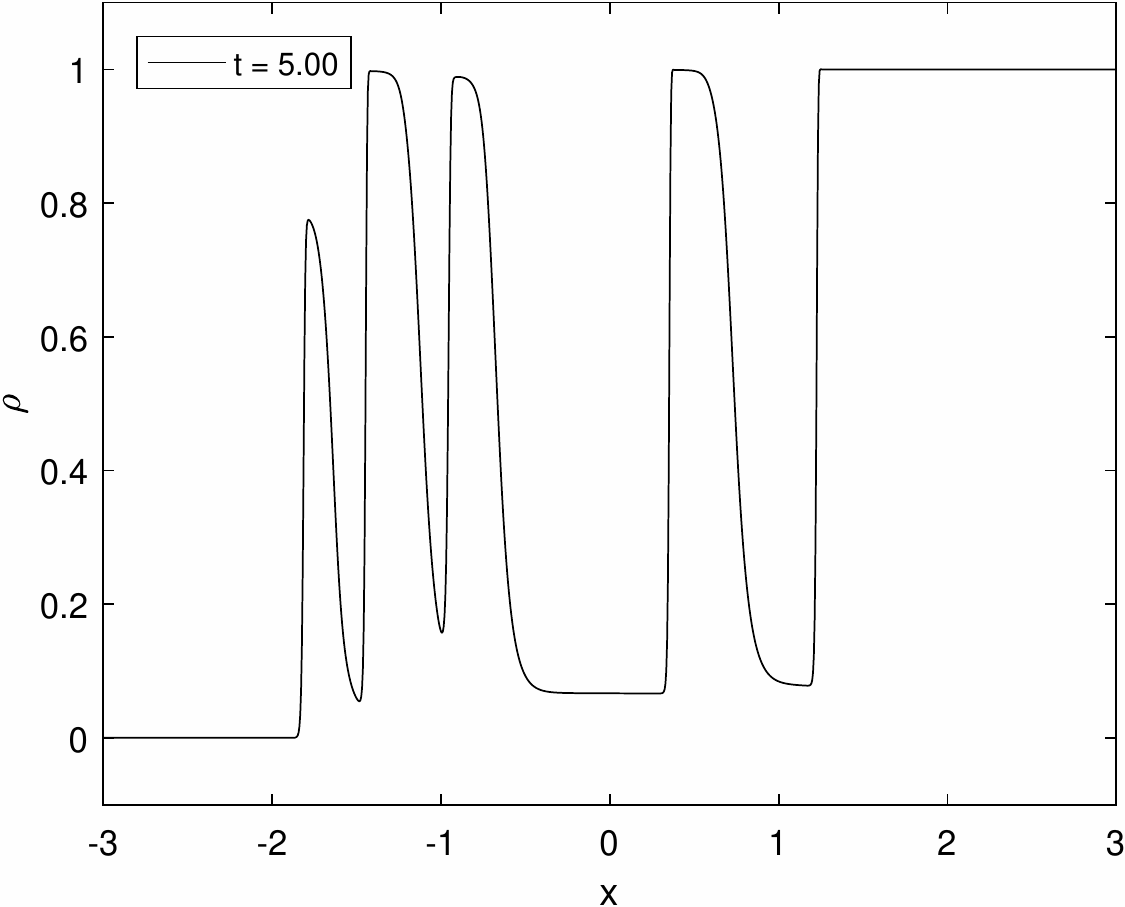}
		\includegraphics[width=0.24\textwidth]{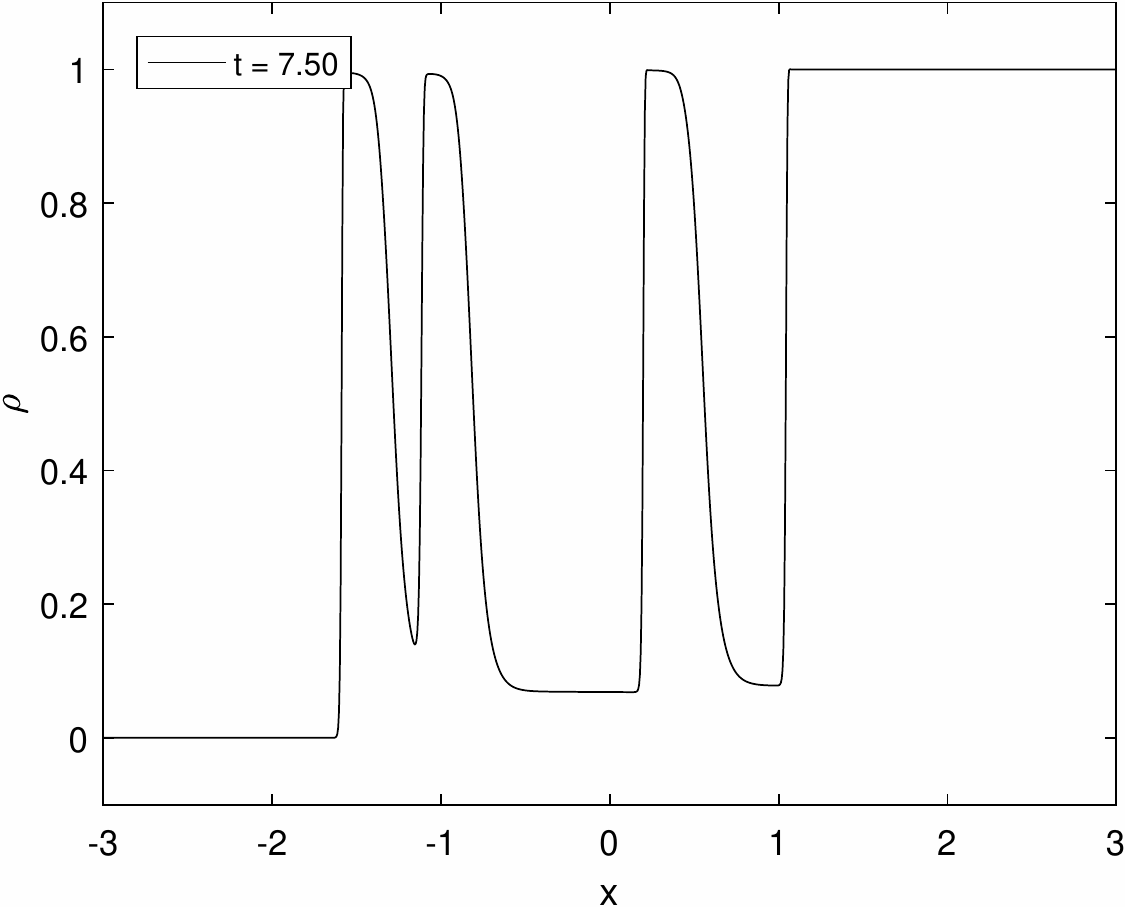}
		\includegraphics[width=0.24\textwidth]{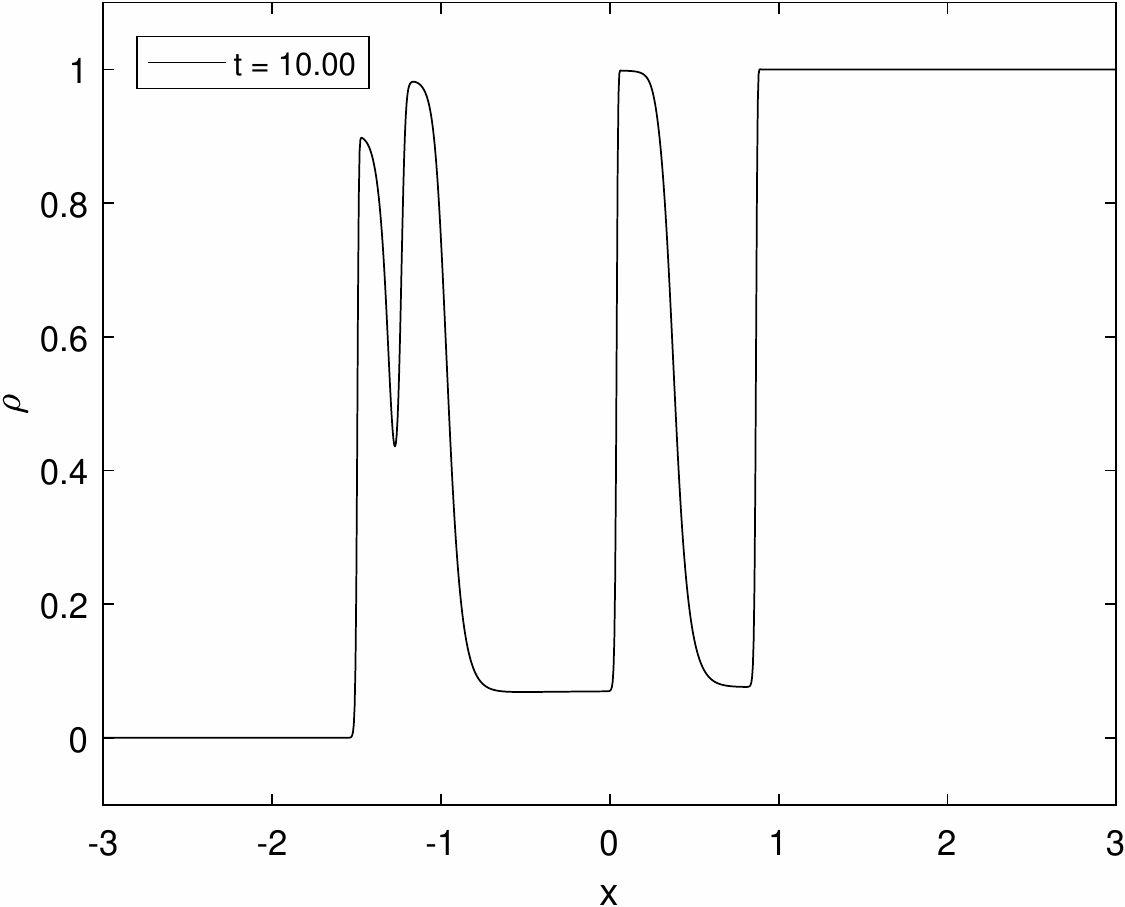}	
		\caption{Backward horizon case with 2048 mesh points at times t = 2.50, 5.01, 7.50, 10.00. The shock heights and velocities agree with the results provided in \cite{amorim2015}.
In Figure \ref{loglogerror}, we can see a first-order behavior of accuracy in the numerical solutions.}
		\label{figlwr2}
	\end{figure}

\begin{figure}[p]
		\centering
		\includegraphics[width=0.24\textwidth]{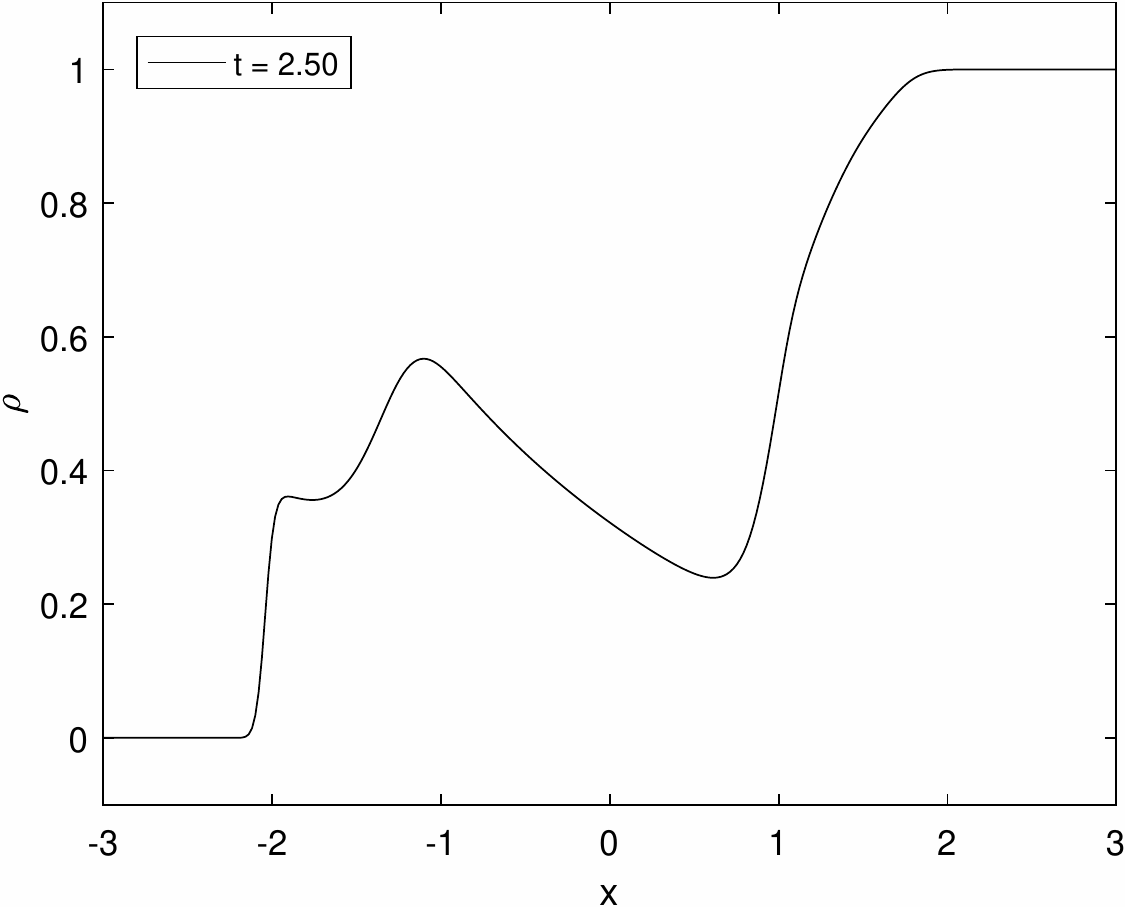}
		\includegraphics[width=0.24\textwidth]{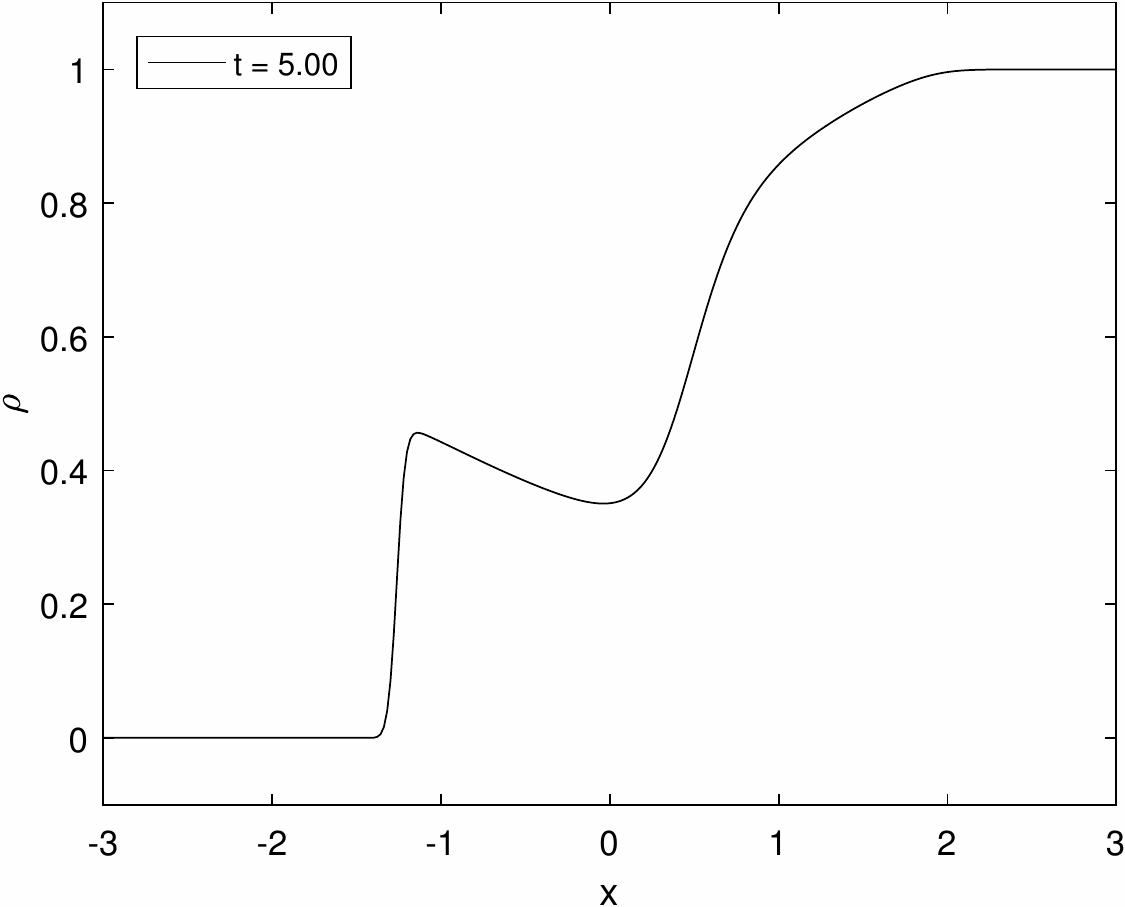}
		\includegraphics[width=0.24\textwidth]{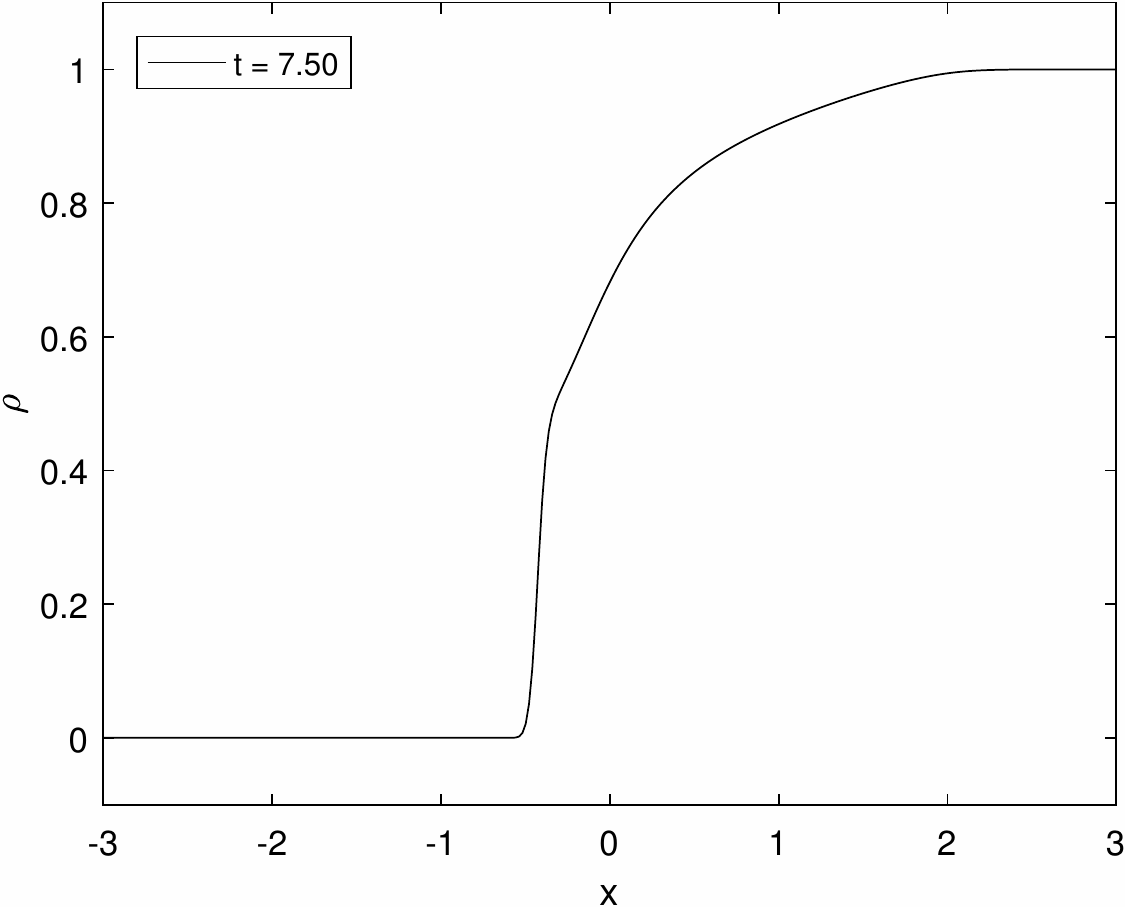}
		\includegraphics[width=0.24\textwidth]{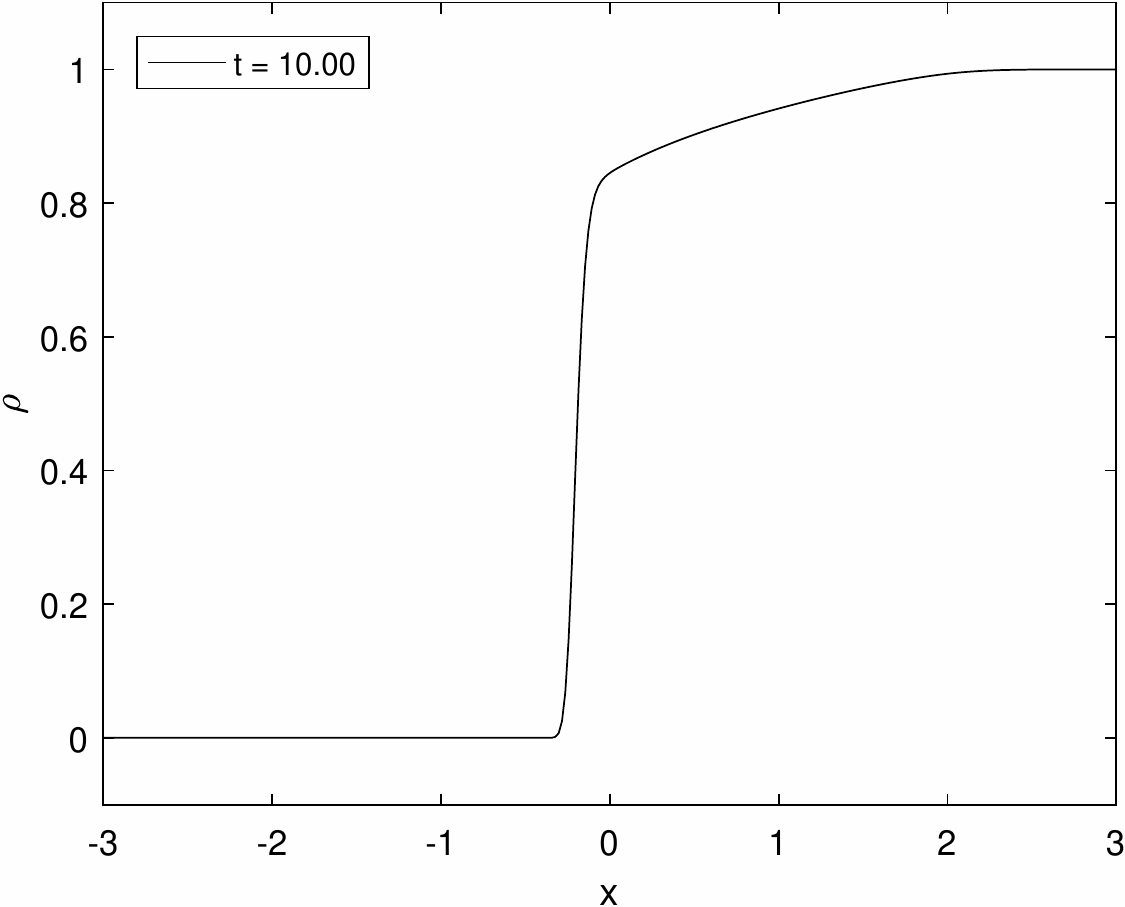}	
		\caption{Forward horizon case with 512 mesh points at times t = 2.50, 5.01, 7.50, 10.00. The expected difference in the two solutions due to the position of the support of $\eta$ was correctly captured by our method.}
		\label{figlwr1}
	\end{figure}

	Parameters $x_1$ and $x_2$ are the 
	\textit{horizon} of each driver, in the sense 
	that a driver situated at $x$ adjusts his speed 
	according to the average vehicular density he 
	sees on interval
	$[x-x_1,x+x_2]$. We followed the exact same 
	numerical approximation of the convolution 
	integral presented in \cite{amorim2015}. And 
	we selected two situations (as in 
	\cite{amorim2015}): (1) the drivers 
	look forward or (2) backward $
	(x_1,x_2) = (0,0.25)$
	\text{and} $
	(x_1,x_2) = (0.25,0).$ The initial condition is given by
	\begin{equation}
	\rho_{0}(x) = \begin{cases}
	\frac{1}{2}, \quad \quad -2.8 \leq x \leq -1.8; \quad\quad\\
	\frac34, \quad \quad -1.2 \leq x \leq -0.2; \\
	\frac34, \quad \quad 0.6 \leq x \leq 1.0;\\
	1, \quad \quad  1.5 \leq x < \infty; \\
	0, \quad \quad  \text{ otherwise.}
	\end{cases} 
	\end{equation}
	which represents three groups of vehicles lining 
	up in a queue. From \cite{amorim2015}, for 
	any $\rho^{,} \in L^{1}(\mathbb{R};[0,1])$, 
	the Cauchy problem with initial datum $\rho_{0}$
	allows a unique solution $\rho = \rho(t,x)$, 
	reaching values in $[0,1]$. The qualitative 
	behaviors of the solution are rather different 
	in the two situations in \eqref{eq:3.2}. The 
	expected big oscillations in the vehicular 
	density caused by the backward-looking case can be 
        seen in Figure \ref{figlwr1} (as opposed to the  
        far more reasonable behavior 
	in the forward-looking scenario in Figure \ref{figlwr2}).
The structure of the numerical 
	solutions presented here are in particularly good 
	agreement with \cite{amorim2015}. We will also present in Table \ref{nonlocerror} an error analysis, so that it is possible to observe that our method presents first-order accuracy behavior.
	
	\begin{figure}[ht]
		\centering
		\includegraphics[width=0.3\textwidth]{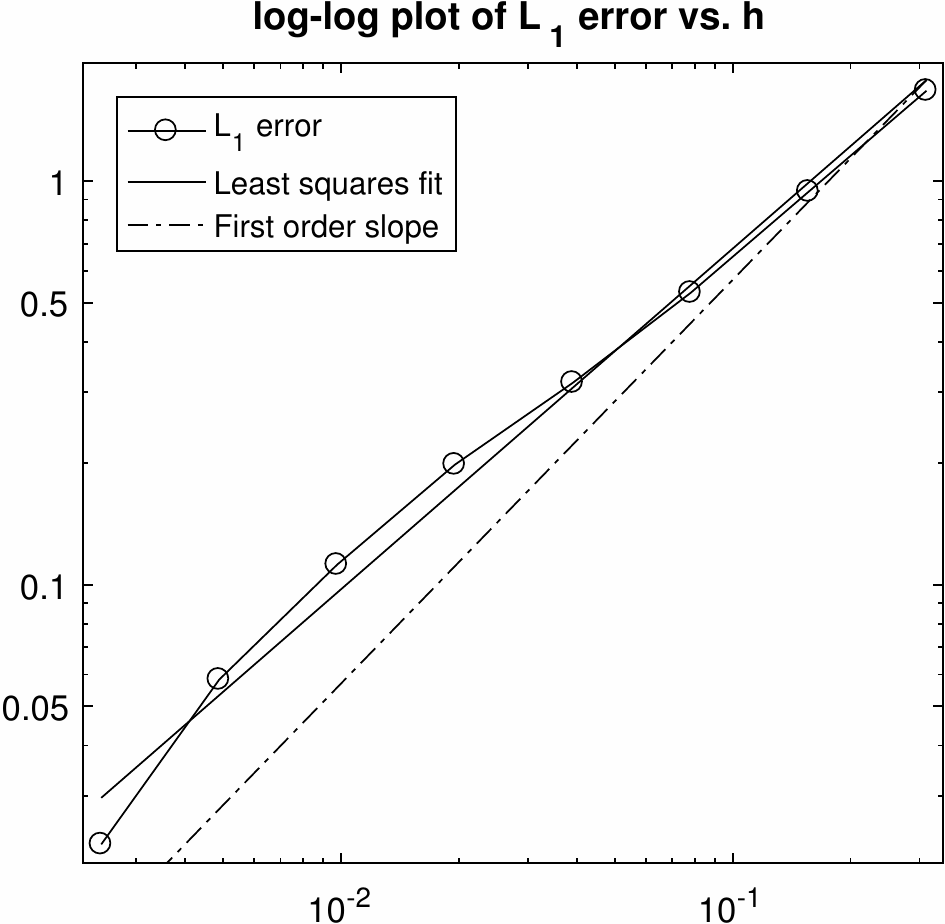}\qquad\qquad
		\includegraphics[width=0.3\textwidth]{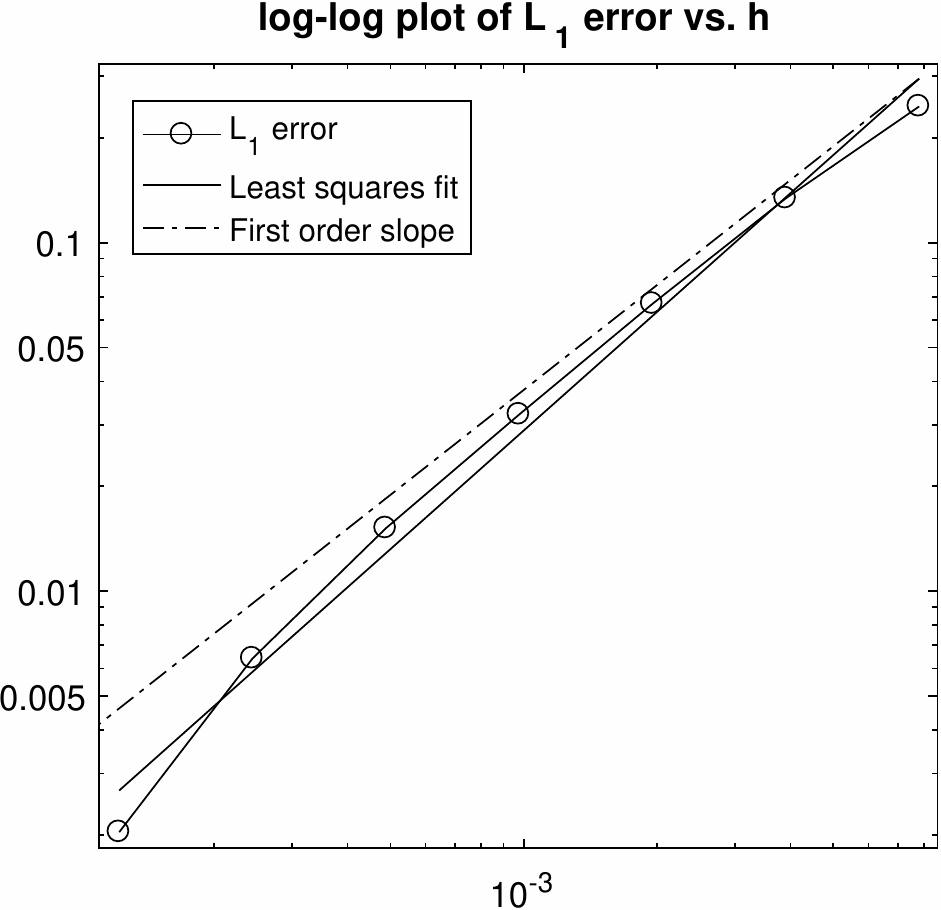}
		\caption{Log-log plots for norm $L^1$ 
            of the error versus the cell sizes, for the traffic problem (\ref{eq:3.2}), at time $T = 0.5$ with backward horizon (left) and problem \ref{eq:22kkr} 
			at time $T = 0.5$ (right).  We can see first-order behavior of accuracy in the numerical 
			solutions.}
		\label{loglogerror}
	\end{figure}

	\begin{table}[ht]
		\begin{center}
			\begin{tabular}{c c c }
				\hline
				Cells &  $h$  &   $ \|u - U\|_{L_h^{1}}$  \\
				\hline 
				$  64$ & \quad $0.15625$ & \quad $9.4\times10^{-1} $ \\
				$ 128$ & \quad $0.07813$ & \quad $5.29\times10^{-1} $ \\
				$ 256$ & \quad $0.03906$ & \quad $3.17\times10^{-1} $ \\
				$ 512$ & \quad $0.01953$ & \quad $1.99\times10^{-1} $ \\
				$1024$ & \quad $0.00976$ & \quad $1.12\times10^{-1} $ \\
				$2048$ & \quad $0.00488$ & \quad $5.82\times10^{-2} $ \\
				$4096$ & \quad $0.00244$ & \quad $2.28\times10^{-2} $ \\
				\hline
				LSF $E(h)$ & \quad & \quad $5.034 \times h^{0.856} $ \\
				\hline
			\end{tabular}\qquad
			\begin{tabular}{c c c }
				\hline
				Cells &  $h$  &   $ \|u - U\|_{L_h^{1}}$  \\
				\hline 
				$  512$ &  $3.90 \times 10^{-3}$ &  $2.46 \times 10^{-1}$\\
				$ 1024$ &  $1.95 \times 10^{-3}$ &  $1.34 \times 10^{-1}$\\
				$ 2048$ &  $9.76 \times 10^{-4}$ &  $6.68 \times 10^{-2}$\\
				$ 4096$ &  $4.88 \times 10^{-4}$ &  $3.21 \times 10^{-2}$\\
				$ 8192$ &  $2.44 \times 10^{-4}$ &  $1.52 \times 10^{-2}$\\
				$16384$ &  $1.22 \times 10^{-4}$ &  $6.4 \times 10^{-3}$\\
				$32768$ &  $6.10 \times 10^{-5}$ &  $2.04 \times 10^{-3}$\\
				\hline
				LSF $E(h)$ &  &  $71.161  \times h^{1.13037} $  \\
				\hline
			\end{tabular}
		\end{center}
		\caption{Left: Corresponding errors 
			between the numerical approximations ($U$) and a reference
			solution ($u$) with 8192 mesh points for the nonlocal problem.
			Right: Corresponding errors 
			between the numerical approximations ($U$) and a reference
			solution ($u$) for the variable $u_1$ with 65536 mesh 
			points for the Keyfitz-Kranzer system problem. The bottom row in both tables presents least square fits for the error profiles.}
        \label{nonlocerror}
	\end{table}

	\subsection{The $2\times2$ symmetric Keyfitz--Kranzer system}
	
\begin{figure}[ht]
		\centering
		\includegraphics[width=0.35\textwidth]{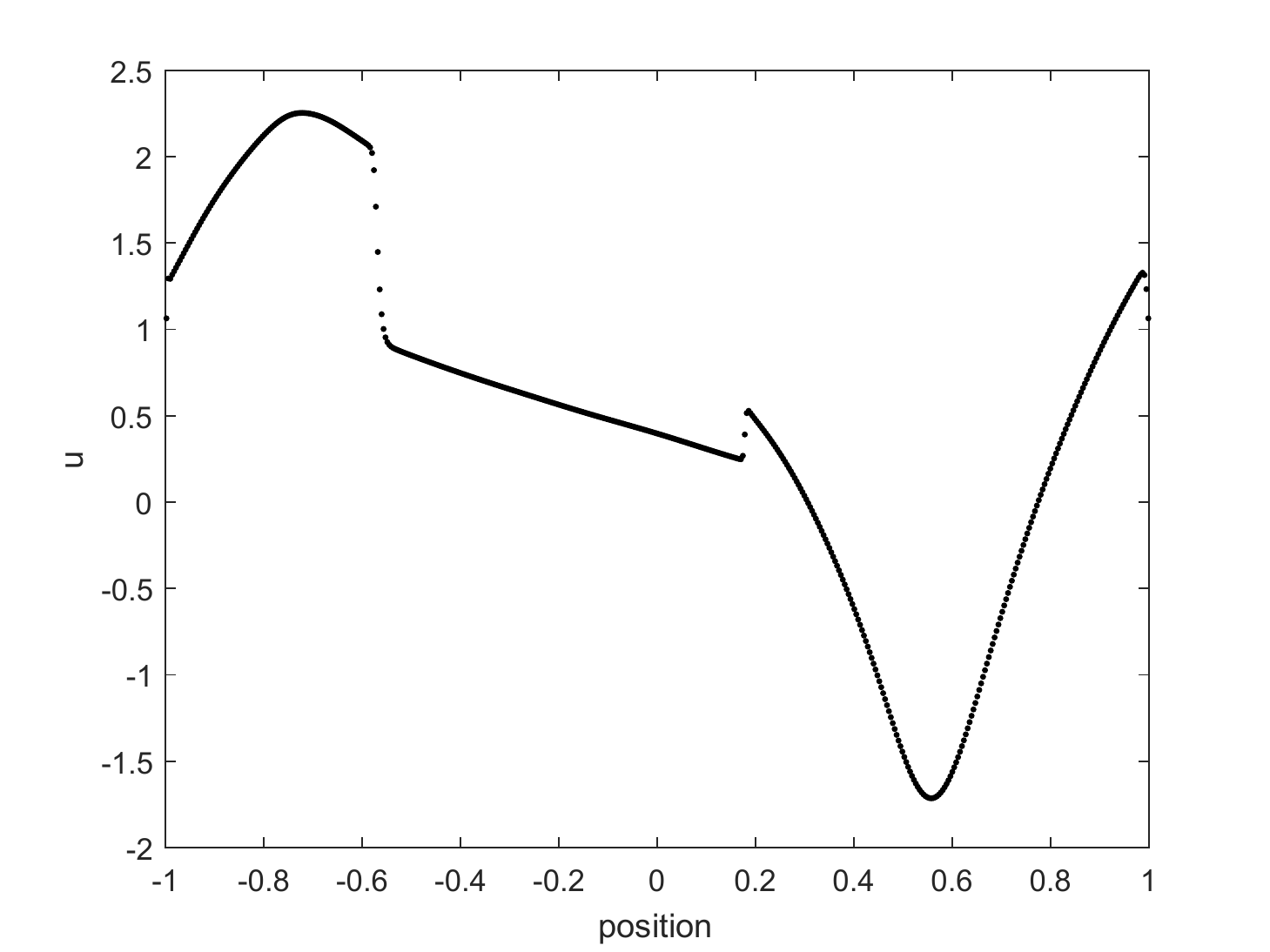}
		\includegraphics[width=0.35\textwidth]{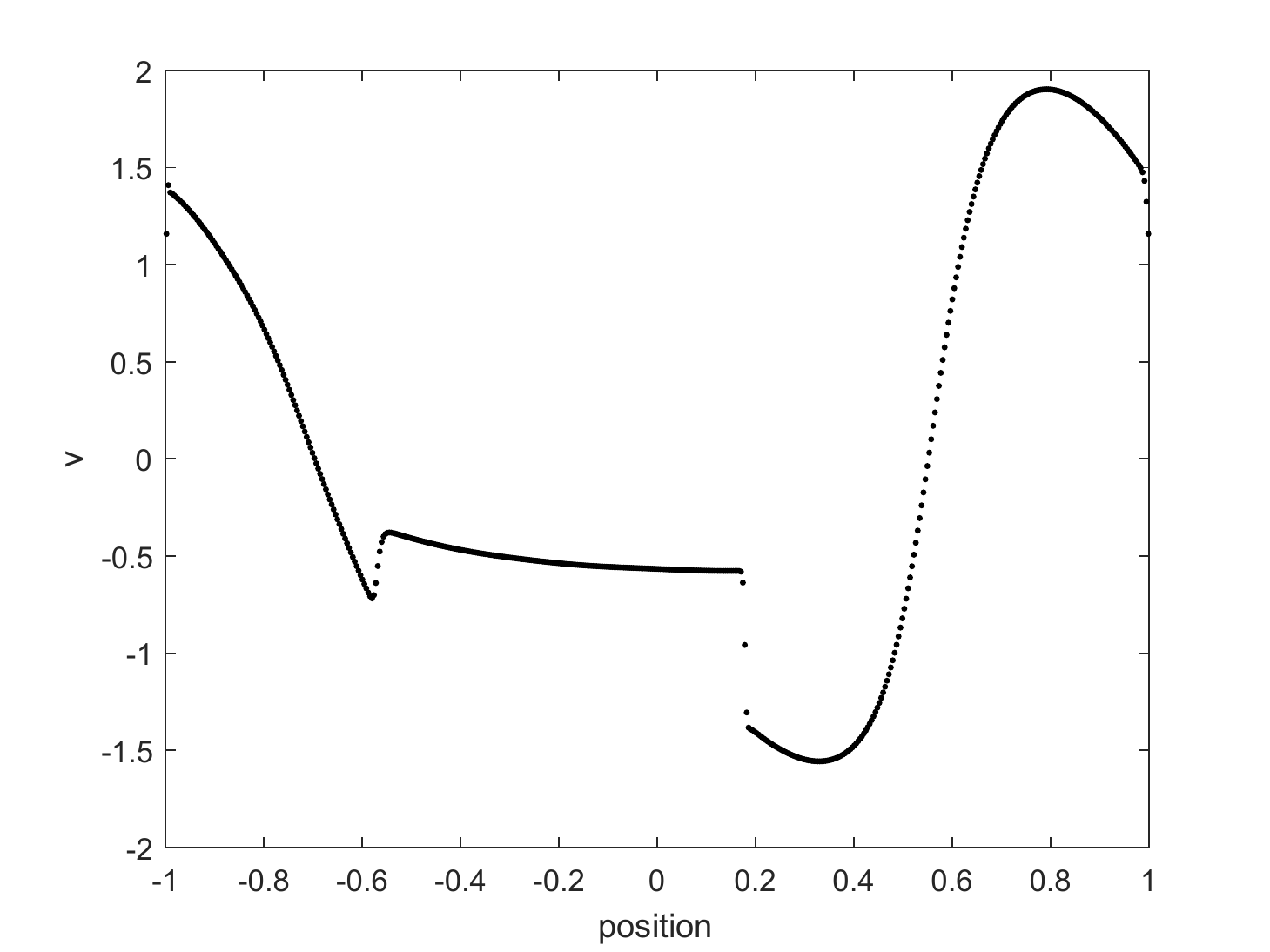}
		\caption{Symmetric Keyfitz--Kranzer system with 512 points at time t = 0.50, $u_1$ (left) and $u_2$ (right).}
		\label{figkk}
	\end{figure}

	We consider the Cauchy problem for the $2\times2$ 
	 Keyfitz--Kranzer system as in \cite{risebro2013},
	
	\begin{equation}
	\begin{cases}
	u_t + (u \phi(|u|))_x  = 0, &\quad (x,t) \in \mathbb{R}\times(0,T),  \\
	u(x,0) = u_0(x), &\quad x \in \mathbb{R},
	\end{cases} 
	\label{eq:22kk}
	\end{equation}
	where $T>0$ is the final simulation time, $u$ is the unknown solution, such that $u = (u_1,u_2) : \mathbb{R} \times (0, T) \mapsto \mathbb{R}^n$ with $|u| := \sqrt{u_1^2+u_2^2}$. The initial datum is given by $u_0 =(u_1^0,u_2^0) \in L^\infty(\mathbb{R}, \mathbb{R}^n)$ and $\phi(r)$ is a scalar function such that $\phi(r) \in C^1(\mathbb{R}^+)$, where $r\phi(r)\to 0$ as $r\to0$. This system 
	is a prototype of a non strictly hyperbolic system of 
	conservation laws, serving as a model for the elastic 
	string (see \cite{keyfitz1980}). Nevertheless, such model appears in magnetohydrodynamics, where it has been used, for example, 
	to explain certain features of the solar wind, such as in
	\cite{cohen1974}. We can rewrite Eq. 
	(\ref{eq:22kk}) in a more explicit form by introducing 
	variable $r = |u|$ to approximate its 
	strong generalized entropy solution (as in \cite{risebro2013}):
	\begin{equation}
	\begin{cases}
	r_t + (r \phi(r))_x  = 0, &\quad (x,t) \in \mathbb{R}\times(0,T),  \\
	u_t + (u \phi(r))_x  = 0, &\quad (x,t) \in \mathbb{R}\times(0,T),  \\
	u(x,0) = u_0(x), \; r(x,0) = r_0(x) = |u_0(x)|, &\quad x \in \mathbb{R}.
	\end{cases} 
	\label{eq:22kkr}
	\end{equation}
	%
	Here, $\phi(r) = r^2 - 4r + 5.5$. This 
	function has a minimum at $r = 2$; hence, the
	ordering of the eigenvalues changes,
        a nonconvex flux function. We tested the method 
	with the following initial data
	$r_0 = \sin(\pi x) + 1.5$, $v_0 = (\sin(\pi x), 
	\cos(\pi x))$,\linebreak$ x \in [-1, 1]$ with periodic boundary 
	conditions. The solution to the Riemann problem 
	with left state $u_L$ and right state $u_R$
	consists of left, right, and middle states separated 
	by shocks, rarefaction waves, or contact discontinuities 
	along which only $r$ changes and by contact 
	discontinuities along which only $v= u/r$ changes. 
	Figure \ref{figkk} illustrates the precise structure 
	of the solution, which is also in good agreement with 
	\cite{risebro2013}.
	
		\section{Concluding remarks}\label{ConRem}
		In this paper, we constructed a Lagrangian--Eulerian 
		framework as a novel tool for balancing discretization 
		in order to deal with nonlinear wave formations and rarefaction 
		interactions in several applications. By implementing the weak 
		asymptotic method, we used $L^1$-norm as well as a comparison between numerical studies and the W1 distance. As a result, we concluded that we were effectively 
		computing the expected approximate 
		solutions linked to problems exhibiting intricate 
		nonlinear wave formations; for instance, Burgers' equation (with initial data containing jumps that also exhibit non-linear rarefaction waves with sonic points and non-linear wave interaction in shock-wave focusing process in time-space), a nonlocal Lighthill-Whitham-Richards model for vehicular traffic model, and a $2 \times 2$ symmetric Keyfitz--Kranzer system. The weak asymptotic 
		solutions we computed with our novel 
		Lagrangian--Eulerian framework have been shown to coincide 
		with classical (regular) solutions and weak Kruzhkov 
		entropy solutions. Our 
		scheme is promising, and it has shown it is a 
		suitable foundation to develop novel 
		constructive methods in abstract as well as practical 
		computational mathematics settings.
	 \FloatBarrier
	\appendix 
	\section*{Appendix}
\numberwithin{equation}{section}

		\section{Proof that the reconstructions are Lipschitz continuous}
	\label{lip}
	In Section $\ref{convergencee}$, to prove the 
	convergence of our numerical scheme, we assumed that the 
	reconstructions were Lipschitz continuous. The 
	reconstructions were obtained from the slope limiters 
	discussed in Section $\ref{slope}$. Here, we prove, 
	first, that each slope limiter (Eqs. $(\ref{mm2})$, 
	$(\ref{mm3})$ and $(\ref{uno})$) is a Lipschitz function, and second, that the reconstruction 
	for each case is also Lipschitz continuous.
	\subsection*{\textbf{The function $MM_2$ is Lipschitz continuous}}
	Let $\vec{x}_1=(x_1,y_1)$ and $\vec{x}_2=(x_2,y_2)$, 
	thus, from Eq. $(\ref{mm})$, we obtain
	\begin{equation*}
	||\MM_2(\vec{x}_1)-\MM_2(\vec{x}_2)||=\left|{\frac{1}{2}\left[\sgn(x_1)+\sgn(y_1)\right]\min(|x_1|,|y_1|)-
		{\frac{1}{2}\left[\sgn(x_2)+\sgn(y_2)\right]\min(|x_2|,|y_2|)}}\right|.
	\end{equation*}
	We then divide our analysis into two possibilities:
	
	\textit{(i)} - First,  we assume that $x_1\geq 0$, 
	$y_1\geq 0$, $x_2\geq0$ and $y_2\geq0$ (the case in which all the negatives are similar). Thus 
	$\frac12\left[\sgn(x_1)+\sgn(y_1)\right]\min(|x_1|,|y_1|)= 
	\hat{x}_1$, where $\hat{x}_1=x_1$ or 
	$\hat{x}_1=y_1$, we also have that \linebreak $\frac12
	\left[\sgn(x_2)+\sgn(y_2)\right]\min(|x_2|,|y_2|)= 
	\hat{x}_2$, where $\hat{x}_2=x_2$ or 
	$\hat{x}_2=y_2$. Thus
$	||\MM_2(\vec{x}_1)-\MM_2(\vec{x}_2)||=|\hat{x}_1-\hat{x}_2|.
$ 
	Notice that, if $\hat{x}_1=x_1$ and $\hat{x}_2=x_2$, 
	then $|\hat{x}_1-\hat{x}_2|=|x_1-x_2|$, and we have that
	\begin{equation*}
	||\MM_2(\vec{x}_1)-\MM_2(\vec{x}_2)||=|\hat{x}_1-
	\hat{x}_2|=|x_1-x_2|\leq |x_1-x_2|+|y_1-y_2|=\norm{\vec{x}_1-\vec{x}_2}.
	\end{equation*}
  
	On the other hand, if $\hat{x}_1=y_1$ and $\hat{x}_2=x_2$, 
	then $|\hat{x}_1-\hat{x}_2|=|y_1-x_2|$. In this case,  
	we have that $x_1>y_1$.  If $y_1>x_2$, then we have 
	that $|y_1-x_2|=y_1-x_2<x_1-x_2=|x_1-x_2|\leq |x_1-x_2|+|y_1-y_2|$. If
	$y_2>x_2>y_1$, then $|y_1-x_2|=x_2-y_1<y_2-y_1=|y_1-y_2|\leq |x_1-x_2|+|y_1-y_2|$, thus
$	||\MM_2(\vec{x}_1)-\MM_2(\vec{x}_2)||=|\hat{x}_1-\hat{x}_2|
	=\leq |x_1-x_2|+|y_1-y_2|=\norm{\vec{x}_1-\vec{x}_2}.
$ 
	The case in which $\hat{x}_1=y_1$, $\hat{x}_2=y_2$,	$\hat{x}_1=x_1$, and $\hat{x}_2=y_2$ is analogous to the previous one.
	
	\textit{(ii)} - Now, we assume that $x_1>0$, $y_1<0$, 
	$x_2>0$, and $y_2>0$ (the other cases for which we 
	have one pair with different signals and another pair 
	with equal signals are similar to this case). 
	Notice that, $\MM_2(\vec{x}_1)=0$, thus
$	||\MM_2(\vec{x}_1)-\MM_2(\vec{x}_2)||=|\hat{x}_2|=\hat{x}_2,
$
	where $\hat{x}_2=x_2$ or $\hat{x}_2=y_2$. If 
	$\hat{x}_2=x_2$, then we have that $\hat{x}_2=
	x_2<x_2-y_1<y_2-y_1\leq |x_2-x_1|+|y_2-y_1|$, where 
	$x_2<x_2-y_1$ because $y_1$ is negative. On the other 
	hand, if $\hat{x}_2=y_2$, then we have that 
	$\hat{x}_2=y_2<y_2-y_1\leq |x_2-x_1|+|y_2-y_1|$. 
	In any case
$	||\MM_2(\vec{x}_1)-\MM_2(\vec{x}_2)||=|\hat{x}_2|
	=\hat{x}_2<||\vec{x}_1-\vec{x}_2||.
$
	Then, the function $\MM_2$ is Lipschitz continuous and its Lipschitz 
	constant equals 1.
	\subsection*{\textbf{The function $\MM_3$ is Lipschitz continuous}}
	Let $\vec{x}_1=(x_1,y_1,z_1)$, $\vec{x}_2=
	(x_2,y_2,z_2)$ and $\MM_2(x,y)$ be Lipschitz continuous 
	with a constant equal to 1
	\begin{align}
	&\norm{\MM_3(\vec{x}_1)-\MM_3(\vec{x}_2)}=\norm{\MM_2
		\left(\MM_2(x_1,y_1),z_1\right)-\MM_2\left(\MM_2(x_2,y_2),z_2\right)}\notag\\
	&\leq |\MM_2(x_1,y_1)-\MM_2(x_2,y_2)|+|z_1-z_2|\leq |x_1-x_2|+|y_1-y_2|+|z_1-z_2|=\norm{\vec{x}_1-
		\vec{x}_2}.\label{ettbs}
	\end{align}
	From inequality $(\ref{ettbs})$, we prove that 
	$\MM_3$ is Lipschitz continuous with a constant equal to 1.  
	\subsection*{\textbf{The functions measuring variations are Lipschitz continuous}}
	The function $U^\prime_j$ that measures variations is 
	defined in Eqs. $(\ref{mm2})$, $(\ref{mm3})$, and $(\ref{uno})$. 
	
	{For Eq. $(\ref{mm2})$}, we defined
$	U^\prime_j=\MM_2\left(F_1(u_{j+1},u_j,u_{j-1})\right),$
    where $F_1(x,y,z)=(x-y,y-z)$. If we prove that 
	$F_2(x,y,z)$ is a Lipschitz function, then, since 
	$U^\prime_j$ is a composition of two Lipschitz continuous, $U^\prime$ is also Lipschitz continuous. 
	Let $\vec{x}_1=(x_1,y_1,z_1)$ and $\vec{x}_2=(x_2,y_2,z_2)$, thus
	\begin{align*}
	\norm{F_1(\vec{x}_1)-F_1(\vec{x}_2)}&=
	\norm{(x_1-x_2-(y_1-y_2),y_1-y_2-(z_1-z_2))}\notag\\
	&\leq |x_1-x_2|+2|y_1-y_2|+|z_1-z_2|\leq 2\norm{\vec{x}_1-\vec{x}_2}.
	\end{align*}    
	
	For Eq. $(\ref{mm3})$, we defined
$	U^\prime_j=\MM_3\left(F_2(u_{j-1},u_j,u_{j+1})\right), $
	where $F_2(x,y,z)=(\alpha(x-y),x-z,\alpha(y-z))$ and 
	$\alpha$ is a nonnegative number. We have now proven that 
	$F_2$ is a Lipschitz function, thus 
	$U^\prime_j$ is also Lipschitz continuous. Let $\vec{x}_1=
	(x_1,y_1,z_1)$ and $\vec{x}_2=(x_2,y_2,z_2)$, thus
	\begin{align*}
	\norm{F_2(\vec{x}_1)-F_2(\vec{x}_2)}&=
	\norm{\left(\alpha(x_1-x_2-(y_1-y_2)),x_1-x_2-(z_1-z_2),
		\alpha(y_1-y_2-(z_1-z_2))\right)}\notag\\
	&\leq \alpha|x_1-x_2|+\alpha|y_1-y_2|+|x_1-x_2|+
	|z_1-z_2|+\alpha|y_1-y_2|+\alpha|z_1-z_2|\notag\\
	&\leq (\alpha+1)|x_1-x_2|+2\alpha|y_1-y_2|+
	(\alpha+1)|z_1-z_2|
\leq (\alpha+2)\norm{\vec{x}_1-\vec{x}_2}.
	\end{align*}  
	
	For Eq. $(\ref{uno})$, we defined 
$	U^\prime_j=\MM_2\left(H(u_{j+2},u_{j+1},u_j,u_{j-1},u_{j-2})\right). $
	Here, $H$ is a more complex function we 
	defined from other auxiliary functions.
	We defined 
	\begin{equation*}
	F_3(x,y,z,w)=(x-2y+z,y-2z+w)
	\quad \text{ and } \quad
	F_4(x,y,z,w)=\frac{1}{2}\MM_2(F_3(x,y,z,w)).
	\end{equation*}

	Using similar calculations, we can prove that 
	$F_3(x,y,z,w)$ is a Lipschitz function. 
	Since $F_4(x,y,z,w)$ is a composition of Lipschitz function, it is also a Lipschitz function. 
	Notice that the function $\delta^2$ that appears in 
	Eq. $(\ref{uno})$ can be written as
$	\delta^2(u_{j+2},u_{j+1},u_{j},u_{j-1})=F_4(u_{j+2},u_{j+1},u_{j},u_{j-1}).
$
	We defined function $H(x,y,z,w,u)$ as
$	H(x,y,z,w,u)=(y-z-F_4(x,y,z,w),z-w+F(y,z,w,u)).$
	Notice that $H(x,y,z,w,u)$ is obtained as a sum 
	of Lipschitz function; thus, it is 
	Lipschitz as well.
	
	\begin{remark}
Since the reconstructions were obtained from linear combinations of slope limiters and function measuring variations, these reconstructions are also Lipschitz continuous.
\end{remark}
\section{The extension of derivatives of $f^+$ and $f^-$}
\label{apendicef}
Sometimes, we are interested in using results for which a continuous derivative of $f^+$ and $f^-$ is necessary. Since these functions are well defined, and their derivative is not well defined only on some points for which $f$ changes their signal, then we can extend the derivative of $f^+$ and $f^-$ in a continuous (but not smooth) way.
First, we assume that $f(x)$ has only a finite number of zeros for which $f$ changes their signal (we are disregarding the zeros for which $f$ does not change their signal); for instance, we denote these zeros as $x_0<x_1\cdots <x_n$. For the sake of simplicity, we assume that $n=2k$ for some $k\in \mathbb{N})$ (if $n=2k+1$, we use a similar argument).

We assume that $f$ satisfies
\begin{equation*}
f(x)=\begin{cases}
f(x)>0,& \text{ if } x\in(-\infty,x_0)\\
f(x)<0,& \text{if } x\in (x_0,x_1)\\
\vdots\\
f(x)>0,& \text{ if } x\in(x_{n-1},x_n)\\
f(x)<0,& \text{if } x\in (x_n,\infty)\\
f(x)=0,& \text{if } x=\{x_0,\cdots,x_n\}
\end{cases} \Rightarrow  
f(x)=\begin{cases}
f(x)>0,& \text{ if } x\in(-\infty,x_0)\bigcup \left[\bigcup\limits_{i=0}^{N/2-1}(x_{2i+1},x_{2i+2})\right]\\
f(x)=0,& \text{if } x=\{x_0,\cdots,x_n\}\\
f(x)<0,& \text{ if } x\in \left[\bigcup\limits_{i=0}^{N/2-1}(x_{2i},x_{2i+1})\right]\bigcup(x_n,\infty)
\end{cases}
\end{equation*}
The derivatives of $f^+$ and $f^-$ are not defined in $x_i$ for $i=0,\cdots,n$. To obtain a continuous extension, we define the derivative of $f^+$, denoted as $(\hat{f}^+)^\prime$, as
\begin{equation}
(\hat{f}^+)^\prime=\begin{cases}
f^\prime(x),& \text{ if } x\in(-\infty,x_0)\bigcup\left[\bigcup\limits_{i=0}^{N/2-1}(x_{2i+1},x_{2i+2})\right].\\
f^\prime(x_{2i})\left(\frac{x_{2i}+\delta-x}{\delta}\right),& \text{ if } x\in\left[\bigcup\limits_{i=0}^{N/2-1}[x_{2i},x_{2i}+\delta]\right]\bigcup[x_n,x_n+\delta].\\
f^\prime(x_{2i+1})\left(\frac{x-x_{2i+1}+\delta}{\delta}\right),& \text{ if } x\in\bigcup\limits_{i=0}^{N/2-1}(x_{2i+1}-
\delta,x_{2i+1}).\\
0,& \text{ if } x\in \left[\bigcup\limits_{i=0}^{N/2-1}(x_{2i}+\delta,x_{2i+1}-\delta]\right]\bigcup(x_n+\delta,\infty).
\end{cases}\label{fdefinem}
\end{equation}
In Eq. $(\ref{fdefinem})$, $\delta$ is arbitrary. For instance, we can take $\delta$, thus satisfying
$\delta=\min_{i=\{0,\cdots,n-1\}}\left(\frac{x_{i+1}-x_{i}}{3}\right).
$ Using $f=f^+-f^-$ and therefore $f^-=f^+-f$, we propose a continuous extension for the derivative of $f^-$, denoted as $
(\hat{f}^-)^\prime$, as
$
(\hat{f}^-)^\prime=(\hat{f}^+)^\prime-f^\prime.
$
\section{The pre-compactness of sequence $u(x,t,\epsilon)$}
\label{precompact}
To prove that the sequence $u(x,t,\epsilon)$ is pre-compact, we used the results in another paper \cite{ACP16}. 
The first result we need is Lemma 1 in \cite{ACP16}:

\textbf{Lemma 1}. \textit{Suppose that} $u(x)\in L^1(\mathbb{T}^n)$, $h>0$. \textit{Then}
\begin{equation*}
\int_{\mathbb{T}^n}|u(x)(\sgn u)^h(x)-|u(x)||dx\leq 2\omega^x(h),
\end{equation*}
\textit{where} 
$\omega^x(h)=\sup_{|\Delta x|\leq h}\int_{\mathbb{T}^n}|u(x+\Delta x)-u(x)|dx,$ \textit{is the continuity modulus of} $u(x)$ in $L^1(\mathbb{T}^n)$.

Here, $\mathbb{T}^n$ is the $n$-dimensional torus. In this study, we are interested in a one-dimensional problem. For $n=1$, $\mathbb{T}^n$ reduces to ${\mathbb{S}^1}$.
Since the proof of the previous Lemma does not depend on the scheme, we refer to \cite{ACP16}. Notice that $\omega^x(h)$ is a measure of TVNI\textsubscript{$\epsilon$}, as described in Eq. $(\ref{tv})$. Thus, under the same hypothesis of Proposition \ref{Prop5}, we can prove the following Corollary:
\begin{corollary}
\label{cor1}
Let us assume that $u(x,t,\epsilon)$ is given by scheme $(\ref{ODE})$ and satisfies the hypothesis of Proposition $\ref{Prop5}$. Then, for all $t>0$, $\Delta x \in \mathbb{R}$, we have that
\begin{equation*}
\int_{\mathbb{S}^1}|u(x+\Delta x,t,\epsilon)-u(x,t,\epsilon)|dx\leq \int_{\mathbb{S}^1}|u_0(x+\Delta x,t,\epsilon)-u_0(x,t,\epsilon)|\leq 
\omega^x(|\Delta x|),
\end{equation*}
where 
$\omega^x(|\Delta x|)\leq \sup_{|\Delta x|\leq h}\int_{\mathbb{S}^1}|u_0(x+\Delta x,t,\epsilon)-u_0(x,t,\epsilon)|$ is the continuity modulus of the initial data $u_0(x)$ in ${\mathbb{S}^1}$.
\end{corollary}

The proof of Corollary $\ref{cor1}$ follows from Proposition $\ref{Prop5}$. 
Now, we prove the result to obtain the pre-compactness of sequence $u(x,t,\epsilon)$. The first useful result, similar to that obtained in \cite{ACP16}, is
\begin{lemma}
\label{lema2} Let us assume that $\phi(x)\in C^1({\mathbb{S}^1})$. Then $\forall \Delta t>0$,
\begin{equation}
\int_{\am{\mathbb{S}^1}} (u(c,t+\Delta t,\epsilon)-u(x,t,\epsilon)\phi(x)dx\leq N ||\nabla \phi||_{\infty} \Delta t\mu(\mathbb{S}^1).\label{dtt}
\end{equation}
\end{lemma}
Here, $\mu(\mathbb{S}^1)$ is the measure of $\mathbb{S}^1$ and
$
N=\max_{|u|\leq \bar{M}}(|u| (|f^+(\hat{u})|+|f^-(\hat{u})|)
\quad \text{ and }\quad 
\bar{M}=||u_0||_\infty.
$

\textbf{Proof}. Let us denote $I(t)=\int_{{\mathbb{S}^1}}u(x,t,\epsilon)\phi(x)$. Differentiating $I(t)$ from $t$ and using $(\ref{ODE})$, we have that
\begin{align}
I^\prime(t)&=\frac{1}{\epsilon}\int_{{\mathbb{S}^1}}
\left(u_{_{\epsilon-1}}f^+\left(\hat{u}_{_{\epsilon-1/2}}\right)-
u_{_{\epsilon}}f^+\left(\hat{u}_{_{\epsilon+1/2}}\right)-
u_{_{\epsilon}}f^-\left(\hat{u}_{_{\epsilon-1/2}}\right)+
u_{_{\epsilon+1}}f^-\left(\hat{u}_{_{\epsilon+1/2}}\right)\right)\phi(x)dx\notag\\
&=\int_{{\mathbb{S}^1}} u_{_{\epsilon}}f^+\left(\hat{u}_{_{\epsilon+1/2}}\right)\frac{\phi(x+\epsilon)-\phi(x)}{\epsilon}dx-
\int_{{\mathbb{S}^1}}u_{_{\epsilon}}f^-\left(\hat{u}_{_{\epsilon+1/2}}\right)\frac{\phi(x+\epsilon)-\phi(x)}{\epsilon}dx.\label{eqt}
\end{align}
Since $I^\prime(t)=G(t)$ implies that $|I(t+\Delta t)-I(t)|\leq \max G(t)\Delta t$, we can estimate RHS of Eq. $(\ref{eqt})$ as
\begin{equation*}
|RHS|\leq \int_{{\mathbb{S}^1}} |u_{_{\epsilon}}||f^+\left(\hat{u}_{_{\epsilon+1/2}}\right)|\left|\frac{\phi(x+\epsilon)-\phi(x)}{\epsilon}\right|dx+
\int_{{\mathbb{S}^1}}|u_{_{\epsilon}}||f^-\left(\hat{u}_{_{\epsilon+1/2}}\right)|\left|\frac{\phi(x+\epsilon)-\phi(x)}{\epsilon}\right|dx.
\end{equation*}
Using 
$\left|\frac{\phi(x\pm\epsilon)-\phi(x)}{\epsilon}\right|dx\leq ||\nabla \phi||_{\infty}$ and $|u(x,t,\epsilon)|\leq \bar{M}$, Eq. $(\ref{dtt})$. follows $\quad \quad \square$

Since we obtained similar estimates in \cite{ACP16}, we used Lemma 3 in reference \cite{ACP16}.

\textbf{Lemma 3}. \textit{For every} $t\geq 0$, $\Delta t>0$
\begin{equation*}
\int_{{\mathbb{S}^1}}|u(x,t+\Delta t,\epsilon)-u(x,t,\epsilon)|dx \leq \omega^t(\Delta t),
\end{equation*}
\textit{where} $\omega^t(\Delta t)=\inf_{h>0}(4\omega^x(h)+cN\Delta t/h)$, \textit{and $c$ is a universal constant}.

Note that, in $\omega^t(\Delta t)$, since this parameter is the infimum, $\omega^t(\Delta t)$ for fixed $\Delta t$ reduces to $\inf_{h>0}(4\omega^x(h))$.

Moreover, since $\omega^x(h)\longrightarrow 0$ as $h\longrightarrow 0$ and does not depend on $\epsilon$ (based on previous results), family $u(x,t,\epsilon)$ is uniformly bounded and equicontinuous in $L^1({\mathbb{S}^1}\times [0,T])$ for every $T>0$. Thus, $u(x,t,\epsilon)$ is a pre-compact sequence in $L^1({\mathbb{S}^1}\times [0,T])$, which implies that we can extract a sequence $\epsilon_k\longrightarrow 0$ such that $u_k(x,t)=u(x,t,\epsilon_k)\longrightarrow u(x,t)$ as $k\longrightarrow\infty$ in $L_{loc}^1({\mathbb{S}^1}\times [0,\infty])$.

\printendnotes


\end{document}